\documentclass[a4,10pt]{article}
\usepackage{amsmath,amscd,amssymb}
\usepackage[margin=.8in]{geometry}

\newcommand{\nbiga}{\mathcal{A}}

\newcommand{\nbigh}{\mathcal{H}}

\newcommand{\nbigm}{\mathcal{M}}

\newcommand{\nbigo}{\mathcal{O}}

\newcommand{\nbigu}{\mathcal{U}}

\newcommand{\seisuu}{{\mathbb Z}}

\newcommand{\cnum}{{\mathbb C}}
\newcommand{\real}{{\mathbb R}}

\newcommand{\gbigr}{\mathfrak R}

\newcommand{\vecv}{{\boldsymbol v}}

\newcommand{\lrarr}{\longrightarrow}

\newcommand{\pf}{{\bf Proof}\hspace{.1in}}
\newcommand{\qed}{\mbox{\rule{1.2mm}{3mm}}}

\def\Hom{\mathop{\rm Hom}\nolimits}

\def\End{\mathop{\rm End}\nolimits}

\def\Image{\mathop{\rm Im}\nolimits}

\def\Re{\mathop{\rm Re}\nolimits}

\def\SL{\mathop{\rm SL}\nolimits}

\def\rank{\mathop{\rm rank}\nolimits}

\def\sym{\mathop{\rm sym}\nolimits}
\def\ad{\mathop{\rm ad}\nolimits}
\def\Res{\mathop{\rm Res}\nolimits}

\def\tr{\mathop{\rm tr}\nolimits}
\def\Tr{\mathop{\rm Tr}\nolimits}

\def\dvol{\mathop{\rm dvol}\nolimits}

\def\id{\mathop{\rm id}\nolimits}

\newcommand{\del}{\partial}
\newcommand{\delbar}{\overline{\del}}

\newcommand{\barz}{\overline{z}}
\newcommand{\zbar}{\barz}

\newcommand{\Hbar}{\overline{H}}

\newcommand{\Abar}{\overline{A}}

\newcommand{\lefttop}[1]{{}^{#1}\!}

\def\Def{\mathop{\rm Def}\nolimits}

\newcommand{\openclosed}[2]{]#1,#2]}

\newcommand{\htilde}{\widetilde{h}}

\newcommand{\atilde}{\widetilde{a}}

\newcommand{\Dtilde}{\widetilde{D}}

\def\asym{\mathop{\rm asym}\nolimits}
\def\hor{\mathop{\rm hor}\nolimits}
\def\ver{\mathop{\rm ver}\nolimits}
\def\semiflat{\mathop{\rm sf}\nolimits}
\def\Pic{\mathop{\rm Pic}\nolimits}
\def\aux{\mathop{\rm aux}\nolimits}

\newcommand{\Ptilde}{\widetilde{P}}

\newcommand{\betabar}{\overline{\beta}}

\newcommand{\taubar}{\overline{\tau}}

\newcommand{\Gbar}{\overline{G}}

\newcommand{\gammatilde}{\widetilde{\gamma}}

\newcommand{\Sigmatilde}{\widetilde{\Sigma}}

\newcommand{\Bbar}{\overline{B}}

\newcommand{\ttI}{{\tt I}}

\newcommand{\nubar}{\overline{\nu}}
\newcommand{\bbar}{\overline{b}}
\newcommand{\btilde}{\widetilde{b}}
\newcommand{\ttV}{{\tt V}}
\newcommand{\ttH}{{\tt H}}
\newcommand{\gammabar}{\overline{\gamma}}
\newcommand{\ttHtilde}{\widetilde{\ttH}}

\newtheorem{thm}{Theorem}[section]
\newtheorem{cor}[thm]{Corollary}

\newtheorem{rem}[thm]{Remark}
\newtheorem{lem}[thm]{Lemma}
\newtheorem{prop}[thm]{Proposition}

\newtheorem{assumption}[thm]{Assumption}

\begin{document}

\title{Comparison of the Hitchin metric
and the semi-flat metric \\
in the rank two case}

\author{Takuro Mochizuki\thanks{Research Institute for Mathematical Sciences, Kyoto University, Kyoto 606-8502, Japan, takuro@kurims.kyoto-u.ac.jp}}
\date{}
\maketitle

\begin{abstract}
Let $(E,\theta)$ be a Higgs bundle of rank $2$ and degree $0$
on a compact Riemann surface $X$ whose spectral curve is smooth.
The tangent space of the moduli space of Higgs bundles
at $(E,\theta)$ 
is equipped with two natural metrics
called the Hitchin metric and the semi-flat metric.
It is known that the difference between two metrics along the curve
$(E,t\theta)$ $(t\geq 1)$ decays in an exponential way.
In this paper, we shall study how the exponential rate is
improved.

\vspace{.1in}
\end{abstract}

\section{Introduction}

\subsection{Hitchin metric and semi-flat metric
of the moduli space of Higgs bundles}

Let $X$ be a compact Riemann surface.
Let $\nbigm_H$ denote the moduli space of
stable Higgs bundles of rank $2$ and degree $0$ on $X$.
According to Hitchin \cite{Hitchin-self-duality},
$\nbigm_H$ is equipped with the hyperk\"ahler metric $g_{H,\real}$
called the Hitchin metric.

There also exists the Hitchin fibration
$\Phi_H:\nbigm_H\to \nbiga_H=H^0(X,K_X^{\otimes 2})$,
which admits a section $\nbiga_H\to \nbigm_H$
called a Hitchin section
depending only on the choice of $K_X^{1/2}$.
Let $\nbiga_H'\subset\nbiga_H$ denote the set of regular values.
We set $\nbigm_H'=\Phi_H^{-1}(\nbiga_H')$.
The restriction of the Hitchin fibration
$\nbigm_H'\to \nbiga_H'$ is an algebraic integrable system,
and in this situation 
$\nbigm_H'$ is equipped with the canonically defined hyperk\"ahler metric
$g_{\semiflat,\real}$ \cite{Freed}, called the semi-flat metric.

Let $(E,\theta)\in \nbigm_H'$.
We obtain the curve $(E,t\theta)$ $(t\geq 1)$ in $\nbigm_H'$.
Let $g_H$ and $g_{\semiflat}$ denote the K\"ahler metrics
of $\nbigm_H$ with the natural complex structure
induced by $g_{H,\real}$ and $g_{\semiflat,\real}$.
It was predicted in \cite{Gaiotto-Moore-Neitzke, Gaiotto-Moore-Neitzke2}
and proved in \cite{Fredrickson2} that
there exists $\epsilon>0$ such that
\begin{equation}
\label{eq;23.10.3.2}
 g_{H|(E,t\theta)}
 -g_{\semiflat|(E,t\theta)}
 =O(e^{-\epsilon t})
\end{equation}
as $t\to\infty$.
Such an estimate was obtained in \cite{Fredrickson2}
under a mild assumption on $(E,\theta)$,
and in \cite{Mochizuki-Asymptotic-Hitchin-metric}
without the additional assumption
by a different method.
However, the original prediction in
\cite{Gaiotto-Moore-Neitzke, Gaiotto-Moore-Neitzke2}
suggests an explicit range of $\epsilon$.
In this paper,
we shall partially verify it in the rank two case
by using the method in \cite{Mochizuki-Asymptotic-Hitchin-metric}.

Before explaining our main result more explicitly,
let us recall more detail about the related objects.

\subsubsection{Hitchin metric}

Let $(E,\theta)$ be a Higgs bundle
on a compact Riemann surface $X$
of degree $0$ and $\rank 2$.
Assume that $\tr\theta=0$ and $\det(E)=\nbigo_X$,
for simplicity.
We also assume $(E,\theta)\in\nbigm'_H$,
which is equivalent to the condition that
the spectral curve $\Sigma_{\theta}$ of $(E,\theta)$
is smooth.

According to the fundamental theorem of Hitchin and Simpson,
under the condition,
$(E,\theta)$ has a unique harmonic metric $h$
such that $\det(h)=1$,
i.e.,
the Hitchin equation
$R(h)+[\theta,\theta^{\dagger}_h]=0$ holds,
where
$R(h)$ denotes the curvature
of the Chern connection of $(E,h)$,
and $\theta^{\dagger}_h$ denotes the adjoint of $\theta$
with respect to $h$.
We remark that the induced harmonic metric
on $(\End(E),\ad\theta)$ is also denoted by $h$.

We consider the Deformation complex
$\Def(E,\theta)$
given by
\[
 \begin{CD}
  \End(E)
  @>{\ad\theta}>>
  \End(E)\otimes K_X,
 \end{CD}
\]
where the first term sits at the degree $0$.
It is well known that
\begin{equation}
\label{eq;23.10.19.30}
 T_{(E,\theta)}\nbigm'_H
 \simeq
 H^1(X,\Def(E,\theta)).
\end{equation}
We recall that
an $\End(E)$-valued $1$-form $\sigma$
is called harmonic
if
\[
 (\delbar_E+\ad\theta)\sigma
=(\delbar_E+\ad\theta)^{\ast}_{h,g_X}\sigma=0
\]
where $\delbar_E$ denotes the holomorphic structure of $E$,
and $g_X$ denotes a compact K\"ahler metric $g_X$ of $X$.
It is equivalent to the condition that
\[
 (\delbar_E+\ad\theta)\sigma
=(\del_{E,h}+\ad\theta^{\dagger}_h)\sigma=0,
\]
where $\del_{E,h}$ denotes the $(1,0)$-part of
the Chern connection of $(E,h)$.
Then, (\ref{eq;23.10.19.30})
is also isomorphic to the space
$\nbigh^1(\End(E),\ad\theta,h)$
of harmonic $1$-forms of $(\End(E),\ad\theta,h)$,

We consider the standard $L^2$-Hermitian product
of $\End(E)$-valued $1$-forms $\sigma_i$ $(i=1,2)$:
\[
 (\sigma_1,\sigma_2)_{L^2,h}
 =2\sqrt{-1}\int_X
 \Tr\Bigl(
 \sigma_1^{1,0}
 \cdot
 (\sigma_2^{1,0})^{\dagger}_h
-\sigma_1^{0,1}
 \cdot
 (\sigma_2^{0,1})^{\dagger}_h
 \Bigr).
\]
Then, the Hitchin metric $g_{H|(E,\theta)}$
is obtained as the restriction of
$(\cdot,\cdot)_{L^2,h}$
to $\nbigh^1(\End(E),\ad\theta,h)
\simeq T_{(E,\theta)}\nbigm_H'$.

\subsubsection{Semi-flat metric}
\label{subsection;24.7.6.10}

Note that $\Phi_H:\nbigm_H'\to\nbiga_H'$
is a locally principal torus bundle.
It is a Lagrangian fibration,
and the fiber
$\Phi_H^{-1}(\Phi_H(E,\theta))$
is isomorphic to
$\Pic_{2g(X)-2}(\Sigma_{\theta})$.
There exists a Hitchin section $\nbiga_H'\to\nbigm_H'$.
There uniquely exists an integrable connection of
$\Phi_H$
such that
any Hitchin section is horizontal.
We obtain the decomposition
of $T_{(E,\theta)}\nbigm_H'$
into the vertical part
$(T_{(E,\theta)}\nbigm_H')^{\ver}$
and the horizontal part
$(T_{(E,\theta)}\nbigm_H')^{\hor}$.
Because the fiber is identified with
a component of the Picard variety of $\Sigma_{\theta}$,
there exists an isomorphism
\[
 (T_{(E,\theta)}\nbigm_H')^{\ver}
 \simeq
 H^1(\Sigma_{\theta},\nbigo_{\Sigma_{\theta}}).
\]
Because the vertical part is a Lagrangian subspace,
$(T_{(E,\theta)}\nbigm_H')^{\hor}$
is the dual space of
$(T_{(E,\theta)}\nbigm_H')^{\ver}$,
and hence
\[
(T_{(E,\theta)}\nbigm_H')^{\hor}
\simeq
 H^0(\Sigma_{\theta},K_{\Sigma_{\theta}}).
\]
Therefore, we obtain
\[
 T_{(E,\theta)}\nbigm_H'
 \simeq
 H^1(\Sigma_{\theta},\nbigo_{\Sigma_{\theta}})
 \oplus
 H^0(\Sigma_{\theta},K_{\Sigma_{\theta}}).
\]
It is also identified with
the space
$\nbigh^1(\Sigma_{\theta})$
of the harmonic $1$-forms of $\Sigma_{\theta}$.

We consider the standard $L^2$-Hermitian product
of $1$-forms $\mu_i$ $(i=1,2)$ on $\Sigma_{\theta}$:
\[
 (\mu_1,\mu_2)_{L^2,\Sigma_{\theta}}
 =2\sqrt{-1}\int_{\Sigma_{\theta}}
 \Bigl(
 \mu_1^{1,0}\wedge
 \overline{\mu_2^{1,0}}
-\mu_1^{0,1}\wedge
 \overline{\mu_2^{0,1}}
 \Bigr).
\]

Then, the semi-flat metric
$g_{\semiflat|(E,\theta)}$
is obtained as the restriction of
$(\cdot,\cdot)_{L^2,\Sigma_{\theta}}$
to $\nbigh^1(\Sigma)\simeq T_{(E,\theta)}\nbigm'_H$.

\subsection{An auxiliary Hermitian product}

\subsubsection{Decompositions}

Let $D_{\theta}$ denote the set of the critical values of
$\pi:\Sigma_{\theta}\to X$.
Let $Q\in X\setminus D_{\theta}$.
There exists a neighbourhood $X_Q$ of $Q$
and a decomposition
\[
 (E,\theta)_{|X_Q}
=(E_{Q,1},\theta_{Q,1})
 \oplus
 (E_{Q,2},\theta_{Q,2}).
\]
By setting
$\End(E)_{|X_Q}^{\circ}
=\End(E_{Q,1})\oplus \End(E_{Q,2})$
and
$\End(E)_{|X_Q}^{\bot}
=\Hom(E_{Q,1},E_{Q,2})\oplus
 \Hom(E_{Q,2},E_{Q,1})$,
we obtain
the decomposition 
\[
 \End(E)_{|X_Q}
=\End(E)_{|X_Q}^{\circ}\oplus\End(E)_{|X_Q}^{\bot}.
\] 
By varying $Q\in X\setminus D_{\theta}$,
we obtain the decomposition
\[
 \End(E)_{|X\setminus D_{\theta}}
=\End(E)^{\circ}
 \oplus
 \End(E)^{\bot}.
\]
For any local section $s$ of
$\End(E)$ on an open set in $X\setminus D_{\theta}$,
we obtain the decomposition
$s=s^{\circ}+s^{\bot}$.
We have similar decompositions
for any $\End(E)$-valued $1$-forms
on an open subset of $X\setminus D_{\theta}$.

\subsubsection{Functions and $1$-forms on the spectral curves}

We recall there exists a holomorphic line bundle $L$
on $\Sigma_{\theta}$
such that $\pi_{\ast}(L)\simeq(E,\theta)$.
Here, $\pi_{\ast}(L)$ is equipped with
the Higgs field induced by the tautological $1$-form
of $T^{\ast}X$,
and the isomorphism preserves the Higgs fields.

Let $U\subset X\setminus D_{\theta}$ be an open subset.
Let $\alpha$ be a function
on $\Sigma_{\theta|U}=\pi^{-1}(U)\subset\Sigma_{\theta}$,
then the multiplication of $\alpha$ on $L$
induces
an endomorphism $F_{\alpha}$ of $E_{|U}$
such that $[\theta,F_{\alpha}]=0$.
Similarly,
a $1$-form $\tau$ on $\Sigma_{\theta|U}$
induces
an $\End(E)$-valued $1$-form $F_{\tau}$ on $U$.
If $\alpha$ (resp. $\tau$) is holomorphic,
then $F_{\alpha}$ (resp. $F_{\tau}$)
is also holomorphic.
Note that even if $U\cap D_{\theta}\neq\emptyset$,
a holomorphic function $\alpha$ on $\Sigma_{\theta|U}$
induces a holomorphic endomorphism $F_{\alpha}$
of $E_{|U}$ such that $[\theta,F_{\alpha}]=0$.
A holomorphic $1$-form $\nu$ on $\Sigma_{\theta|U}$
induces a holomorphic section of
$\End(E)\otimes K_X(D_{\theta})$ on $U$.

\subsubsection{Holomorphic $1$-forms on the spectral curve
induced by harmonic metrics}

Let $h_{\infty}$ be the limiting configuration of $(E,\theta)$
in the sense of \cite{MSWW}
(see \S\ref{subsection;23.10.19.40}).
We consider the following
$\End(E)$-valued $1$-form on $X\setminus D_{\theta}$:
\[
 \Psi_h=(\del_{E,h_{\infty}}+\theta^{\dagger}_{h_{\infty}})
 -(\del_{E,h}+\theta^{\dagger}_{h}).
\]
Because
both
$\del_{E,h_{\infty}}+\theta^{\dagger}_{h_{\infty}}$
and 
$\del_{E,h}+\theta^{\dagger}_{h}$
commute with $\delbar_E+\theta$,
we obtain $(\delbar_E+\ad\theta)\Psi_h=0$.
There exists a holomorphic $1$-form $\eta_{(E,\theta)}$
on $\Sigma_{\theta|X\setminus D_{\theta}}$
such that
\[
 (\Psi_h^{1,0})^{\circ}=F_{\eta_{(E,\theta)}}.
\]
\begin{lem}[Lemma
\ref{lem;23.10.19.1}]
$\eta_{(E,\theta)}$ is holomorphic on $\Sigma_{\theta}$. 
\end{lem}

\subsubsection{An auxiliary Hermitian product}

Let $\phi^{1/2}$ denote the holomorphic $1$-form
obtained as the canonical root of
$\pi^{\ast}(-\det\theta)$ on $\Sigma_{\theta}$.
Let $\Dtilde_{\theta}$ be the set-theoretic inverse image
$\pi^{-1}(D_{\theta})$.
For $\nu_i\in H^0(\Sigma_{\theta},K_{\Sigma_{\theta}})$,
we obtain
\[
  \nu_1\cdot\nu_2\cdot\eta_{(E,\theta)}
 \cdot
 (\phi^{1/2})^{-1}
 \in H^0\Bigl(
 \Sigma_{\theta},
 (K_{\Sigma_{\theta}}(\Dtilde_{\theta}))^{2}
 \Bigr).
\]
Note that for each $P\in D_{\theta}$,
there exists a natural map
\[
 \Res^{(2)}_{\pi^{-1}(P)}:
 H^0\bigl(\Sigma_{\theta},
 (K_{\Sigma_{\theta}}(\Dtilde))^2
 \bigr)
 \lrarr\cnum
\]
obtained as the residue.
We define
\[
 \langle \nu_1,\nu_2\rangle^{\aux}_{(E,\theta)}
 =-4\pi
 \sum_{P\in D_{\theta}}
 \Res^{(2)}_{\pi^{-1}(P)}
 \Bigl(
  \nu_1\cdot\nu_2\cdot\eta_{(E,\theta)}
 \cdot
 (\phi^{1/2})^{-1}
  \Bigr).
\]

By using the isomorphisms
$\iota^{\hor}:H^0(\Sigma_{\theta},K_{\Sigma_{\theta}})
\simeq
 (T_{(E,\theta)}\nbigm_H')^{\hor}$
 and
 $\iota^{\ver}:H^1(\Sigma_{\theta},\nbigo_{\Sigma_{\theta}})
\simeq
 (T_{(E,\theta)}\nbigm_H')^{\ver}$,
we define the Hermitian product $g^{\aux}_{(E,\theta)}$
on $T_{(E,\theta)}\nbigm_H'$ by the following conditions.
\begin{itemize}
 \item $g^{\aux}_{(E,\theta)}(u_1,u_2)=0$
       if $u_i\in (T_{(E,\theta)}\nbigm_H')^{\ver}$ $(i=1,2)$
       or $u_i\in (T_{(E,\theta)}\nbigm_H')^{\hor}$ $(i=1,2)$.
 \item For $\nu\in H^0(\Sigma_{\theta},K_{\Sigma_{\theta}})$
       and $\tau\in
       H^1(\Sigma_{\theta},\nbigo_{\Sigma_{\theta}})
       =
       H^0(\Sigma_{\theta}^{\dagger},
       K_{\Sigma_{\theta}^{\dagger}})$,
       we set
       \[
       g^{\aux}_{|(E,\theta)}
       (\iota^{\hor}(\nu),\iota^{\ver}(\tau))
       =\langle
       \nu,\taubar
       \rangle^{\aux}_{(E,\theta)}.
       \]
\end{itemize}

\subsection{Main result}

We set $\phi=-\det(\theta)\in H^0(X,K_X^{\otimes 2})$,
which induces the distance $d_{\phi}$ on $X$.
According to \cite[\S3.2]{Dumas-Neitzke},
a saddle connection of $\phi$
is a geodesic segment on $(X,d_{\phi})$
such that
(i) the starting point and the ending point are
contained in $D_{\phi}$,
(ii) any other point is contained in $X\setminus D_{\phi}$.
Note that the starting point and the ending point
are not necessarily distinct.
Let $M(\phi)$ denote the minimum length
of the saddle connections of $\phi$,
which is called the threshold.

\begin{prop}[Corollary
\ref{cor;23.10.19.11}]
$g^{\aux}_{|(E,t\theta)}=O(e^{-2a t})$
with respect to $g_{\semiflat|(E,t\theta)}$
for any $0<a<M(\phi)$.
As a result,
there exists $t_0>0$ such that
 $g_{\semiflat|(E,t\theta)}
+g^{\aux}_{|(E,t\theta)}$
are positive definite Hermitian products 
on $T_{(E,t\theta)}\nbigm_H'$
for any $t\geq t_0$.
\end{prop}

The following theorem is the main result of this paper.

\begin{thm}[Theorem
\ref{thm;23.10.8.10}]
\label{thm;23.10.19.100}
For any $(E,\theta)\in\nbigm_H'$,
we obtain
\begin{equation}
 g_{H|(E,t\theta)}
 -\bigl(
 g_{\semiflat|(E,t\theta)}
+g^{\aux}_{|(E,t\theta)}
\bigr)
 =O\bigl(e^{-4at}\bigr)
 \quad
 (0<a<M(\phi)).
 \end{equation}
 \end{thm}

\subsubsection{Symmetric case}

A non-degenerate symmetric product of $(E,\theta)$
is a non-degenerate symmetric product
$C:E\otimes E\to \nbigo_X$
such that
$\theta$ is self-adjoint with respect to $C$.
Note that if $(E,\theta)$ is contained in a Hitchin section,
then $(E,\theta)$ has a non-degenerate symmetric product.

\begin{cor}[Lemma
\ref{lem;23.10.19.100}]
\label{cor;23.10.19.101}
If $(E,\theta)$ is equipped with a non-degenerate symmetric product,
then we obtain $\eta_{(E,\theta)}=0$.
As a result, we obtain
\begin{equation}
 \label{eq;23.10.3.3}
g_{H|(E,t\theta)}-g_{\semiflat|(E,t\theta)}
=O(e^{-4at}),
\quad
(0<a<M(\phi)).
\end{equation}
 \end{cor}

\subsubsection{Previous works}

Originally, Gaiotto-Moore-Neitzke
\cite{Gaiotto-Moore-Neitzke, Gaiotto-Moore-Neitzke2}
predicted (\ref{eq;23.10.3.3}).
(See also \cite{Dumas-Neitzke, Fredrickson2, FMSW, Neitzke}.)
In \cite{Dumas-Neitzke},
Dumas and Neitzke studied the issues
in the case where $(E,\theta)$ is contained in
a Hitchin section,
and they proved
\begin{equation}
\label{eq;23.10.3.4}
g_{H|(E,t\theta)}-g_{\semiflat|(E,t\theta)}
=O(e^{-4at}),
\quad
(0<a<M(\phi)/2)
\end{equation}
on the subspace
$T_{t^2\phi}\nbiga_H'\subset
T_{(E,t\theta)}\nbigm_H'$
induced by the Hitchin section.
As we mentioned,
Fredrickson obtained the estimate (\ref{eq;23.10.3.2})
on the basis of \cite{Fredrickson-Generic-Ends, MSWW, MSWW2}.
It was further generalized to the parabolic case
by Fredrickson, Mazzeo, Swoboda and Weiss in \cite{FMSW},
where they suggested that they could also obtain 
the estimate (\ref{eq;23.10.3.4}).
In \cite{Holdt}, Holdt studied the original approach of
\cite{Gaiotto-Moore-Neitzke, Gaiotto-Moore-Neitzke2}
in the case that a Higgs bundle has a non-trivial parabolic structure.
More recently, Chen and Li studied
more singular case in \cite{Chen-Li}.

In \cite{Mochizuki-Asymptotic-Hitchin-metric},
we obtained (\ref{eq;23.10.3.2})
for higher rank case
without any assumption
on the ramification of $\Sigma_{\theta}\to X$
in a way different from \cite{Fredrickson2}.
A naive hope was to obtain (\ref{eq;23.10.3.3})
by the method in \cite{Mochizuki-Asymptotic-Hitchin-metric}.
Corollary \ref{cor;23.10.19.101} is a satisfactory solution
in the symmetric case.
In the general case,
Theorem \ref{thm;23.10.19.100}
implies (\ref{eq;23.10.3.4}),
but it is not a complete answer to
the prediction (\ref{eq;23.10.3.3})
of Gaiotto-Moore-Neitzke
because it is not clear whether
$g^{\aux}_{(E,t\theta)}
=O(e^{-4at})$ $(0<a<M(\phi))$
or not
with respect to $g_{\semiflat|(E,\theta)}$.

\subsection{Outline in the symmetric case}

It would be instructive to explain an outline
to obtain (\ref{eq;23.10.3.3}) directly 
in the case
where $(E,\theta)$ is equipped with
a non-degenerate symmetric pairing $C$
such that $\det(C)=1$,
the arguments can be simplified because some issues disappear.
It is close to the case of \cite{Dumas-Neitzke}.
Note that if $(E,\theta)$ is contained in
a Hitchin section,
then $(E,\theta)$ is equipped with
a non-degenerate symmetric pairing.

Let $D_{\theta}$ denote the zero set of $\phi=-\det(\theta)$.
Let $\Sigma_{\theta}\subset T^{\ast}X$
denote the spectral curve.
Because $(E,\theta)\in\nbigm_H'$,
$\Sigma_{\theta}$ is smooth,
and the projection $\pi:\Sigma_{\theta}\to X$
is a ramified covering of degree $2$.
By \cite{Beauville-Narasimhan-Ramanan, Hitchin-self-duality},
there exists a line bundle $L$ on $\Sigma_{\theta}$
such that
$\pi_{\ast}(L)\simeq (E,\theta)$,
where the Higgs field of $\pi_{\ast}(L)$
is induced by the tautological $1$-form on $T^{\ast}X$.

\subsubsection{The deformation complexes}

Let $\End(E)^{\sym}$ denote the sheaf of
endomorphisms of $E$ which are symmetric with respect to $C$.
Similarly, let $\End(E)^{\asym}$ denote the sheaf of
endomorphisms of $E$ which are anti-symmetric with respect to $C$.
We obtain the following subcomplex of $\Def(E,\theta)$:
\[
 \Def(E,\theta,C)=
 \Bigl(
 \End(E)^{\asym}
 \to
 \End(E)^{\sym}\otimes K_X
 \Bigr).
\]
We may also consider the following subcomplex of
$\Def(E,\theta)$:
\[
\Def'(E,\theta,C)
=\Bigl(
 \End(E)^{\sym}
 \to
 \End(E)^{\asym}\otimes K_X
\Bigr).
\]
We have the decomposition of the complexes:
\begin{equation}
\label{eq;23.10.7.20}
 \Def(E,\theta)
 =\Def(E,\theta,C)
 \oplus
 \Def'(E,\theta,C).
\end{equation}

The natural action $\nbigo_{\Sigma}$ on $L$
induces
$\pi_{\ast}(\nbigo_{\Sigma})
\to
\End(E)$.
It induces a morphism of complexes
$\pi_{\ast}(\nbigo_{\Sigma})
\lrarr
\Def'(E,\theta,C)
\subset\Def(E,\theta)$,
and the isomorphism of sheaves
\[
 \pi_{\ast}(\nbigo_{\Sigma_{\theta}})
 \simeq
 \nbigh^0(\Def'(E,\theta,C))
 \simeq
 \nbigh^0(\Def(E,\theta)).
\]
We also have $\nbigh^1(\Def'(E,\theta,C))=0$.
We obtain the following injection:
\begin{equation}
\label{eq;23.10.5.1}
 H^0(\Sigma_{\theta},\nbigo_{\Sigma_{\theta}})
\simeq
 H^1(X,\Def'(E,\theta,C))
 \lrarr
H^1(X,\Def(E,\theta)).
\end{equation}

As explained in \cite{Mochizuki-Asymptotic-Hitchin-metric},
there exist natural isomorphisms
\[
 \nbigh^1(\Def(E,\theta))\simeq
 \nbigh^1(\Def(E,\theta,C))\simeq
 \pi_{\ast}(K_{\Sigma_{\theta}}).
\]
We also have
$\nbigh^0(\Def(E,\theta,C))=0$.
We obtain the following injection:
\begin{equation}
\label{eq;23.10.5.2}
 H^0(\Sigma_{\theta},K_{\Sigma_{\theta}})
 \simeq
 H^1(X,\Def(E,\theta,C))
 \lrarr
 H^1(X,\Def(E,\theta)).
\end{equation}

\subsubsection{The tangent space, the vertical direction
and the horizontal direction}

The restriction of the Hitchin fibration
$\nbigm_H'\to \nbiga_H'$
is locally a principal torus bundle
equipped with an integrable connection.
There exists the decomposition
into the horizontal direction and the vertical direction:
\begin{equation}
\label{eq;23.10.7.30}
 T_{(E,\theta)}\nbigm_H
=(T_{(E,\theta)}\nbigm_H)^{\hor}
 \oplus
 (T_{(E,\theta)}\nbigm_H)^{\ver}.
\end{equation}

We identify
$T_{(E,\theta)}\nbigm_H
 \simeq
 H^1(X,\Def(E,\theta))$.
The vertical direction
$(T_{(E,\theta)}\nbigm_H')^{\ver}$
equals
$T_{(E,\theta)}\Phi_H^{-1}(\Phi_H(E,\theta))$,
which is identified with
the image of (\ref{eq;23.10.5.1}).
The horizontal direction equals
the image of (\ref{eq;23.10.5.2})
as explained in \cite{Mochizuki-Asymptotic-Hitchin-metric}.
In other words,
the decomposition (\ref{eq;23.10.7.30})
is induced by the decomposition (\ref{eq;23.10.7.20})
in the symmetric case.

\subsubsection{The Hitchin metric and the semi-flat metric}

Let $\tau\in H^1(\Sigma_{\theta},\nbigo_{\Sigma_{\theta}})$
and $\nu\in H^0(\Sigma_{\theta},K_{\Sigma_{\theta}})$.
We regard $\tau$ as an anti-holomorphic $1$-form
on $\Sigma_{\theta}$,
i.e.,
$\tau\in H^0(\Sigma_{\theta}^{\dagger},K_{\Sigma_{\theta}^{\dagger}})$.
Because $\Sigma_{\theta}=\Sigma_{t\theta}$,
there exist the corresponding elements
\[
 \iota_{t}^{\ver}(\tau),
 \iota_t^{\hor}(\nu)
 \in
 H^1(X,\Def(E,t\theta))
 =T_{(E,t\theta)}\nbigm_H.
\]

We have the following:
\begin{equation}
\label{eq;23.10.5.3}
 g_{\semiflat}(\iota_t^{\ver}(\tau),\iota_t^{\ver}(\tau))
=-2\sqrt{-1}\int_{\Sigma_{\theta}}
 \tau\wedge \overline{\tau}
=\|\tau\|_{L^2}^2,
\end{equation}
\begin{equation}
\label{eq;23.10.5.4}
 g_{\semiflat}(\iota_t^{\hor}(\nu),\iota_t^{\hor}(\nu))
=2\sqrt{-1}\int_{\Sigma_{\theta}}
\nu\wedge \overline{\nu}
=\|\nu\|_{L^2}^2,
\end{equation}
\begin{equation}
\label{eq;23.10.5.5}
  g_{\semiflat}(\iota_t^{\hor}(\nu),\iota_t^{\ver}(\tau))=0.
\end{equation}

Let $h_t$ be the harmonic metric of $(E,t\theta)$,
which are compatible with $C$.
The Hitchin metrics $g_{H}$ at $(E,t\theta)$
are induced by
the natural $L^2$-metric on
the space of $\End(E)$-valued $1$-forms.
Namely, for any $L^2$-sections
$\rho_i=\rho_i^{1,0}+\rho_i^{0,1}$ $(i=1,2)$
of $\End(E)\otimes\Omega^{1,0}\oplus\End(E)\otimes\Omega^{0,1}$
on $W\subset X$,
we set
\[
 (\rho_1,\rho_2)_{L^2,h_t,W}
=2\sqrt{-1}\int_W\Tr\Bigl(
 \rho_1^{1,0}\cdot
 (\rho_{2}^{1,0})^{\dagger}_{h_t}
- \rho_1^{0,1}\cdot
 (\rho_{2}^{0,1})^{\dagger}_{h_t}
 \Bigr).
\]
If $W=X$,
$(\rho_1,\rho_2)_{L^2,h_t,X}$
is denoted by
$(\rho_1,\rho_2)_{L^2,h_t}$.
Let $\ttH_t(\nu)$ and $\ttV_t(\tau)$
$(t\geq 1)$ denote the harmonic $1$-forms
of $(\End(E),\ad\theta,h_t)$
corresponding to $\iota^{\hor}_t(\nu)$
and $\iota^{\ver}_t(\tau)$,
respectively.
We have the following:
\begin{equation}
\label{eq;23.10.5.6}
 g_H(\iota_t^{\ver}(\tau),\iota_t^{\ver}(\tau))
 =(\ttV_t(\tau),\ttV_t(\tau))_{L^2,h_t},
\end{equation}
\begin{equation}
\label{eq;23.10.5.7}
 g_H(\iota_t^{\hor}(\nu),\iota_t^{\hor}(\nu))
 =(\ttH_t(\nu),\ttH_t(\nu))_{L^2,h_t},
\end{equation}
\begin{equation}
 g_H(\iota_t^{\hor}(\nu),\iota_t^{\ver}(\tau))
 =(\ttH_t(\nu),\ttV_t(\tau))_{L^2,h_t}.
\end{equation}

Because $C$ is compatible with $h_t$,
the decomposition
$\End(E)=\End(E)^{\sym}\oplus \End(E)^{\asym}$
is orthogonal with respect to $h$.
We can observe that
$\ttH_t(\nu)$ is a section of
the Dolbeault resolution of $\Def(E,\theta,C)$,
and
$\ttV_t(\nu)$ is a section of
the Dolbeault resolution of $\Def'(E,\theta,C)$.
Hence,
we obtain
\[
  g_H(\iota_t^{\hor}(\nu),\iota_t^{\ver}(\tau))
  =0
  =g_{\semiflat}(\iota_t^{\hor}(\nu),\iota_t^{\ver}(\tau)).
\]
We would like to show the following estimates
for any $0<\kappa<M(\phi)/2$:
\begin{equation}
\label{eq;23.10.7.40}
 g_H(\iota_t^{\hor}(\nu),\iota_t^{\hor}(\nu))
 =g_{\semiflat}(\iota_t^{\hor}(\nu),\iota_t^{\hor}(\nu))
+ O\bigl(e^{-8\kappa t}\|\nu\|_{L^2}^2\bigr).
\end{equation}
\begin{equation}
\label{eq;23.10.7.41}
 g_H(\iota_t^{\ver}(\tau),\iota_t^{\ver}(\tau))
 =g_{\semiflat}(\iota_t^{\ver}(\tau),\iota_t^{\ver}(\tau))
+ O\bigl(e^{-8\kappa t}\|\tau\|_{L^2}^2\bigr).
\end{equation}

\subsubsection{Canonically decoupled harmonic metric}

Let $Q\in X\setminus D_{\theta}$.
Let $(X_Q,z_Q)$ denote
a holomorphic coordinate neighbourhood around $Q$
such that $\phi=dz_Q\,dz_Q$.
There exists the decomposition,
which is orthogonal with respect to $C$:
\[
 (E,\theta)_{|X_Q}
=(E_{Q,1},dz_Q)
 \oplus
 (E_{Q,2},-dz_Q).
\]
Let $v_i$ $(i=1,2)$ be frames of $E_{Q,i}$
such that $C(v_i,v_i)=1$.
We define $h_{\infty,Q}$ by
$h_{\infty,Q}(v_i,v_i)=1$
and
$h_{\infty,Q}(v_1,v_2)=0$,
which is independent of $v_i$.
By varying $Q\in X\setminus D_{\theta}$,
we obtain the canonical decoupled harmonic metric
$h_{\infty}$ of $(E,\theta)_{|X\setminus D_{\theta}}$
compatible with $C$.
It equals the limiting configuration in \cite{MSWW}.

Let $s_t$ be the automorphism of $E_{|X\setminus D_{\theta}}$
determined by
$h_t=h_{\infty}s_t$.
We have $\det(s_t)=1$.

We set $g_{\phi}=|\phi|$ as
the K\"ahler metric of $X\setminus D_{\theta}$,
which induces the distance $d_{\phi}$ on $X$.
We obtain the connection $\nabla_{\infty}$ of
$\End(E)\otimes T^{\ast}(X\setminus D_{\theta})^{\otimes m}$
for any $m\in\seisuu_{\geq 0}$
induced by $h_{\infty}$ and $g_{\phi}$.

For any $r>0$ and $P\in X$,
we set $X_P(r)=\{Q\in X\,|\,d_{\phi}(Q,P)<r\}$.
We set $\kappa_0=M(\phi)/2$.
We have
$X_P(\kappa_0)\cap X_{P'}(\kappa_0)=\emptyset$
if $P,P'$ are two distinct points in $D_{\theta}$.

For $0<\kappa<\kappa_0$,
we set
$\delta(\kappa):=
\min\{(\kappa-\kappa_0)/10,\kappa/10\}$.
For each $P\in D_{\theta}$,
we set
\[
X^{(1)}_{P,\kappa}=X_P(\kappa_0-\delta(\kappa)),
\quad
X^{(2)}_{P,\kappa}=X_P(\kappa_0-2\delta(\kappa)).
\]
For any non-negative integer $j$,
there exist $C_j>0$ such that the following holds
on $X\setminus \bigcup X^{(2)}_{P,\kappa}$
for any $t\geq 1$
(see Corollary \ref{cor;23.9.28.40}):
\begin{equation}
\label{eq;23.10.5.11}
 \bigl|
  \nabla_{\infty}^j(s_t-\id)
  \bigr|_{h_{\infty},g_{\phi}}
  \leq
  C_je^{-4(\kappa+7\delta(\kappa))t}.   
\end{equation}
This is essentially contained in \cite{Dumas-Neitzke}.

The multiplication of
$\nu\in H^0(\Sigma_{\theta},K_{\Sigma_{\theta}})$ on $L$ induces
a section 
$F_{\nu}$ of $(\End(E)\otimes\Omega^{1,0})_{|X\setminus D_{\theta}}$
such that $(\delbar_E+t\ad\theta)F_{\nu}=0$.
Similarly, $\tau$ induces
a section
$F_{\tau}$ of $(\End(E)\otimes\Omega^{0,1})_{|X\setminus D_{\theta}}$
such that $(\delbar_E+t\ad\theta)F_{\tau}=0$.
We can regard
$F_{\nu}$ and $F_{\tau}$ as approximations of
$\ttH_t(\nu)$ and $\ttV_{t}(\tau)$
on $X\setminus\bigcup_PX_{P,\kappa}^{(2)}$.

\subsubsection{Approximation of $\ttH_t(\nu)$ around $P\in D_{\theta}$}

For any $W\subset X$,
we set $\Sigma_{\theta|W}=\Sigma_{\theta}\times_XW$.
There exists a holomorphic function $\alpha_P$
on $\Sigma_{\theta|X_P(\kappa_0)}$
such that $\nu=d\alpha_P$
and that $\alpha_P=0$ at $\pi^{-1}(P)$.
The multiplication of $\alpha_P$ on $L$
induces
a holomorphic endomorphism
$F_{\alpha_P}$ on $E_{|X_P(\kappa_0)}$
such that
$[\theta,F_{\alpha_P}]=0$.
We obtain the following $\End(E)$-valued
$C^{\infty}$ $1$-forms for $t\geq 1$:
\[
 \ttH_{P,t}(\nu)=
 (\del_{E,h_t}+t\ad\theta^{\dagger}_{h_t})
 F_{\alpha_P}.
\]
They are harmonic $1$-forms of
$(\End(E),t\ad\theta,h_t)_{|X_P}$,
i.e.,
\[
(\delbar_E+t\ad\theta)
 \ttH_{P,t}(\nu)=0,\quad
(\del_{E,h_t}+t\ad\theta^{\dagger}_{h_t})
 \ttH_{P,t}(\nu)=0.
\]
We have
$\ttH_{P,t}(\nu)^{\circ}
= F_{\nu}$
on $X_P(\kappa_0)\setminus\{P\}$.
There exists a $C^{\infty}$-section
$\rho_P$ of $\End(E)^{\asym}_{|X_P(\kappa_0)\setminus\{P\}}$
such that
\[
 (\delbar_E+t\ad\theta)\rho_{P,t}=
 F_{\nu}-\ttH_{P,t}(\nu).
\]

\subsubsection{Approximation of $\ttH_t(\nu)$ on $X$}

By patching $F_{\nu}$ and $\ttH_{P,t}(\nu)$
$(P\in D_{\theta})$,
we construct $\ttH'_{\kappa,t}(\nu)$ as follows.
\begin{itemize}
 \item On $X\setminus \bigcup X_{P,\kappa}^{(1)}$,
       we set
       $\ttH'_{\kappa,t}(\nu)=F_{\nu}$.
 \item For $P\in D_{\theta}$,
       let $\chi_{P,\kappa}:X\to [0,1]$ be
       a $C^{\infty}$-function
       such that $\chi_{P,\kappa}=1$ on $X_{P,\kappa}^{(2)}$
       and $\chi_{P,\kappa}=0$ on $X\setminus X_{P,\kappa}^{(1)}$.
       Then, on $X^{(1)}_{P,\kappa}\setminus X^{(2)}_{P,\kappa}$,
       we set
       $\ttH'_{\kappa,t}(\nu)=F_{\nu}-(\delbar_E+t\ad\theta)$
       $(\chi_{P,\kappa}\rho_{P})$.
 \item On $X^{(2)}_{P,\kappa}$,
       we set
       $\ttH'_{\kappa,t}(\nu)=\ttH_{P,t}(\nu)$.
\end{itemize}
By the construction,
we obtain 
$(\delbar_E+t\ad\theta)\ttH_{\kappa,t}'(\nu)=0$.
We can observe that
$\ttH'_{\kappa,t}(\nu)$ is a section of
the Dolbeault resolution of
$\Def(E,\theta,C)$,
and that $\ttH_{\kappa,t}'(\nu)^{\circ}=F_{\nu}$.
As explained in \cite{Mochizuki-Asymptotic-Hitchin-metric},
the cohomology class of
$\ttH_{\kappa,t}'(\nu)$
is $\iota_t^{\hor}(\nu)$.

There exist $C_j>0$ $(j\geq 0)$ such that
the following holds
on $X^{(1)}_{P,\kappa}\setminus X^{(2)}_{P,\kappa}$
by (\ref{eq;23.10.5.11}):
\begin{equation}
\label{eq;23.10.16.11}
\Bigl|
\nabla_{\infty}^j
 (\ttH_{P,t}(\nu)-F_{\nu})
\Bigr|_{h_{\infty},g_{\phi}}
\leq C_je^{-4(\kappa+6\delta(\kappa))t}
\|\nu\|_{L^2}.
\end{equation}
There exist $C_j>0$ $(j\geq 0)$ such that
the following holds
on $X^{(1)}_{P,\kappa}\setminus X^{(2)}_{P,\kappa}$:
\[
 \bigl|\nabla^j_{\infty}\rho_{P}\bigr|_{h_{\infty},g_{\phi}}
 \leq C_je^{-4(\kappa+6\delta(\kappa))t}
 \|\nu\|_{L^2}.
\]
Hence, there exist $C_j>0$ $(j\geq 0)$ such that
on $X^{(1)}_{P,\kappa}\setminus X^{(2)}_{P,\kappa}$:
\[
\bigl|
 \nabla^j_{\infty}(\ttH'_{\kappa,t}(\nu)-F_{\nu})
\bigr|_{h_{\infty},g_{\phi}}
\leq C_je^{-4(\kappa+6\delta(\kappa))t}
\|\nu\|_{L^2}.
\]

Let $g_X$ be a K\"ahler metric.
Note that
$(\del_{E,h_t}+t\ad\theta^{\dagger}_{h_t})
 \ttH'_{\kappa,t}(\nu)=0$
on $X\setminus \bigcup_P X^{(1)}_{P,\kappa}$
and $\bigcup_PX_{P,\kappa}^{(2)}$.
There exists $C>0$ such that
the following holds on $X$:
\[
\bigl|
 (\del_{E,h_t}+t\ad\theta^{\dagger}_{h_t})
 \ttH'_{\kappa,t}(\nu)
 \bigr|_{h_{t},g_{X}}
 \leq Ce^{-4(\kappa+6\delta(\kappa))t}
 \|\nu\|_{L^2}.
\]

There exist unique $C^{\infty}$-sections $\gamma_{\kappa,t}(\nu)$
of $\End(E)^{\asym}$ such that
\[
 (\del_{E,h_t}+t\ad\theta^{\dagger}_{h_t})
 (\delbar_E+t\ad\theta)\gamma_{\kappa,t}(\nu)
=(\del_{E,h_t}+t\ad\theta^{\dagger}_{h_t})
 \ttH'_{\kappa,t}(\nu).
\]
By \cite[Proposition 2.39]{Mochizuki-Asymptotic-Hitchin-metric},
there exists $C>0$ such that
\[
 \|\gamma_{\kappa,t}(\nu)\|_{L^2,h_t,g_X}
+\|(\delbar_E+t\ad\theta)\gamma_{\kappa,t}(\nu)\|_{L^2,h_t,g_X} 
 \leq Ce^{-4(\kappa+5\delta(\kappa)) t}\|\nu\|_{L^2}.
\]

We have
\[
 \ttH_t(\nu)
=\ttH'_{\kappa,t}(\nu)-\gamma_{\kappa,t}(\nu).
\]

\subsubsection{Estimate of the pairings of $\ttH_t(\nu)$}

First, we observe
\begin{equation}
\label{eq;23.10.5.10}
 (\ttH_t(\nu),\ttH_t(\nu))_{L^2,h_t}
=(\ttH'_{\kappa,t}(\nu),\ttH'_{\kappa,t}(\nu))_{L^2,h_t}
+O\bigl(e^{-8\kappa t}\|\nu\|_{L^2}^2\bigr).
\end{equation}
Indeed,
because
$(\delbar+t\ad\theta)^{\ast}_{h_t,g_X}\ttH_t(\nu)=0$,
we obtain
\begin{multline}
 (\ttH_t(\nu),\ttH_t(\nu))_{L^2,h_t}
=(\ttH_t(\nu),\ttH'_{\kappa,t}(\nu))_{L^2,h_t}
 +\bigl(\ttH_t(\nu),
 (\delbar_E+t\ad\theta)\gamma_{\kappa,t}(\nu)\bigr)_{L^2,h_t}
=
\\
(\ttH_t(\nu),\ttH'_{\kappa,t}(\nu))_{L^2,h_t}
+\bigl((\delbar_E+t\ad\theta)^{\ast}_{h_t,g_X}\ttH_t(\nu),
\gamma_{\kappa,t}(\nu)\bigr)_{L^2,h_t,g_X}
=(\ttH_t(\nu),\ttH'_{\kappa,t}(\nu))_{L^2,h_t}. 
\end{multline}
Because
$(\delbar+t\ad\theta)^{\ast}_{h_t,g_X}\ttH'_{\kappa,t}(\nu)=
O\bigl(e^{-4\kappa t}\|\nu\|_{L^2}\bigr)$
and
$\gamma_{\kappa,t}(\nu)
=O(e^{-4\kappa t}\|\nu\|_{L^2})$,
we obtain
\[
\bigl(\ttH'_{\kappa,t}(\nu),
(\delbar_E+t\ad\theta)\gamma_{\kappa,t}(\nu)
\bigr)_{L^2,h_t}
=\bigl((\delbar_E+t\ad\theta)^{\ast}_{h_t,g_X}\ttH'_{\kappa,t}(\nu),
\gamma_{\kappa,t}(\nu)\bigr)_{L^2,h_t,g_X}
=O\bigl(e^{-8\kappa t}\|\nu\|_{L^2}^2\bigr).
\]
Hence, we obtain (\ref{eq;23.10.5.10}).

We have the following estimates on
$W_{\kappa}=X\setminus \bigcup X_{P,\kappa}^{(2)}$:
\begin{equation}
\label{eq;23.10.16.10}
 \ttH_{\kappa,t}'(\nu)^{\circ}
=F_{\nu},
\quad\quad
\Bigl|
\bigl(
\ttH'_{\kappa,t}(\nu)
\bigr)^{\bot}
\Bigr|_{h_t,g_X}
=O\bigl(e^{-4\kappa t}\|\nu\|_{L^2}\bigr),
\quad\quad
 \bigl|
 s_{t}-\id
 \bigr|_{h_t}
=O(e^{-4\kappa t}).
\end{equation}
We can prove the following
by using (\ref{eq;23.10.16.10})
with an elementary consideration
as in \S\ref{subsection;23.10.4.1}:
\begin{multline}
(\ttH'_{\kappa,t}(\nu),\ttH'_{\kappa,t}(\nu))_{L^2,h_t,W_{\kappa}}
=
 2\sqrt{-1}
 \int_{W_{\kappa}}
 \Tr\bigl(
 F_{\nu}
 \cdot
 (F_{\nu})^{\dagger}_{h_{\infty}}
 \bigr)
 +O\bigl(e^{-8\kappa t}\|\nu\|_{L^2}^2\bigr)
\\
 =
2\sqrt{-1}
 \int_{\Sigma_{\theta|W_{\kappa}}}
 \nu\wedge \overline{\nu}
 +O\bigl(e^{-8\kappa t}\|\nu\|_{L^2}^2\bigr).
\end{multline}

On $X^{(2)}_{P,\kappa}$,
we have
$\ttH'_{\kappa,t}(\nu)=\ttH_{P,t}(\nu)$.
Because
$(\del_{E,h_t}+t\theta^{\dagger}_{h_t})^{\ast}_{h_t,g_X}
\ttH_{P,t}(\nu)=0$,
we obtain the following
by using the Stokes formula (Lemma \ref{lem;23.10.7.20}):
\[
 \bigl(
 \ttH'_{\kappa,t}(\nu),
 \ttH'_{\kappa,t}(\nu)
 \bigr)_{L^2,h_t,X^{(2)}_{P,\kappa}}
=
 2\sqrt{-1}
 \int_{\del X_{P,\kappa}^{(2)}}
 \Tr\Bigl(
 F_{\alpha_P}
 \cdot
 (\ttH_{P,t}(\nu)^{1,0})^{\dagger}_{h_t}
\Bigr).
\]
By using the estimates (\ref{eq;23.10.5.11})
and (\ref{eq;23.10.16.11})
with the elementary consideration
in \S\ref{subsection;23.10.4.1},
we obtain
\begin{multline}
 2\sqrt{-1}
 \int_{\del X_{P,\kappa}^{(2)}}
 \Tr\Bigl(
 F_{\alpha_P}
 \cdot
 (\ttH_{P,t}(\nu)^{1,0})^{\dagger}_{h_t}
\Bigr)
= 
2\sqrt{-1}\int_{\del X_{P,\kappa}^{(2)}}
 \Tr\Bigl(
 F_{\alpha_P}
 \cdot
 (F_{\nu})^{\dagger}_{h_{\infty}}
 \Bigr)
 +O\bigl(e^{-8\kappa t} \|\nu\|^2_{L^2}\bigr)
 \\
 =
 2\sqrt{-1}\int_{\del\Sigma_{\theta|X_{P,\kappa}^{(2)}}}
 \alpha_P\cdot \nubar
 +O\bigl(e^{-8\kappa t} \|\nu\|^2_{L^2}\bigr)
=2\sqrt{-1}\int_{\Sigma_{\theta|X_{P,\kappa}^{(2)}}}
 \nu\nubar
 +O\bigl(e^{-8\kappa t} \|\nu\|^2_{L^2}\bigr).
\end{multline}
In all, we obtain
\[
 (\ttH_{t}(\nu),\ttH_t(\nu))_{L^2,h_t}
 =(\nu,\nu) _{L^2}
 +O\bigl(e^{-8\kappa t}\|\nu\|_{L^2}^2\bigr)
 \quad
 (0<\kappa<M(\phi)/2).
\]
This implies (\ref{eq;23.10.7.40}).

\subsubsection{Approximation of $\ttV_t(\tau)$ on $X$}

For each $P\in D_{\theta}$,
there exists an anti-holomorphic function $\beta_P$
on $\Sigma_{\theta|X_P(\kappa_0)}$
such that $\tau=d\beta_P$.
Let $\chi_{P,\kappa}:X\to[0,1]$ be as before.
We set
\[
 \tau^{\circ}_{\kappa}
 =\tau-\sum_{P\in D_{\theta}}
 \delbar(\pi^{\ast}(\chi_{P,\kappa})\cdot\beta_P).
\]
The multiplication of $\tau^{\circ}_{\kappa}$ on $L$
induces
a $C^{\infty}$-section $F_{\tau^{\circ}_{\kappa}}$ of
$\End(E)\otimes\Omega^{0,1}$
such that
$[\theta,F_{\tau^{\circ}_{\kappa}}]=0$.
We obtain
$(\delbar_E+t\ad\theta)F_{\tau^{\circ}_{\kappa}}=0$.
The cohomology class of
$F_{\tau^{\circ}_{\kappa}}$
is $\iota^{\ver}_t(\tau)$.

We may regard $\pi^{\dagger}:\Sigma_{\theta}^{\dagger}\to X^{\dagger}$
as the spectral curve of
the Higgs bundle $(E,\del_{E,h_t},\theta^{\dagger}_{h_t})$
on $X^{\dagger}$.
We may regard $\beta_P$ as a holomorphic function on
$\Sigma^{\dagger}_{\theta|X_P(\kappa_0)}$.
There exists $L_t$ on $\Sigma^{\dagger}_{\theta}$
with an isomorphism
$\pi^{\dagger}_{\ast}(L_t)\simeq (E,\del_{E,h_t},\theta^{\dagger}_{h_t})$.
The multiplication of $\beta_P$ induces
a holomorphic endomorphism $F^{\dagger}_{\beta_P}$
of $(E,\del_{E,h_t})$.
We have $F^{\dagger}_{\beta_P}=(F_{\overline{\beta_P}})^{\dagger}_{h_t}$.
We note that
$(\delbar_E+t\theta)F^{\dagger}_{\beta_P}$
is a harmonic $1$-form of
$(\End(E),t\ad\theta,h_t)$.

We set
\[
 \ttV'_{\kappa,t}(\tau)
=F_{\tau^{\circ}_{\kappa}}
+\sum_{P\in D_{\theta}}
(\delbar_E+t\ad\theta)(\chi_{P,\kappa}F^{\dagger}_{\beta_P}).
\]
By the construction,
we have
$(\delbar_E+t\ad\theta)\ttV'_{\kappa,t}(\tau)=0$,
and its cohomology class equals
$\iota_t^{\ver}(\tau)$.

By (\ref{eq;23.10.5.11}),
there exists $C>0$ such that
\[
\Bigl|
(\del_{E,h_t}+t\ad\theta^{\dagger}_{h_t})
 \ttV'_{\kappa,t}(\tau)
 \Bigr|_{h_t,g_X}
 \leq
 Ce^{-4(\kappa+7\delta(\kappa))t}\|\tau\|_{L^2}.
\]
We also have
$\del_X\tr(\ttV'_{\kappa,t}(\tau))=0$.
There exist unique sections $\gamma_{\kappa,t}(\tau)$
of $\End(E)$ such that
\[
 (\del_{E,h_t}+t\ad\theta^{\dagger}_{h_t})
 (\delbar_E+t\ad\theta)\gamma_{\kappa,t}(\tau)
=(\delbar_{E,h_t}+t\ad\theta^{\dagger}_{h_t})
 \ttV'_{\kappa,t}(\tau).
\]
By \cite[Proposition 2.39]{Mochizuki-Asymptotic-Hitchin-metric},
there exists $C>0$ such that
\[
 \|
 \gamma_{\kappa,t}(\tau)
 \|_{L^2,h_t,g_X}
+\|(\delbar_E+t\ad\theta)
 \gamma_{\kappa,t}(\tau)
 \|_{L^2,h_t,g_X}
 \leq Ce^{-4(\kappa+5\delta(\kappa)) t}\|\tau\|_{L^2}.
\]
We have
\[
 \ttV_t(\tau)=\ttV'_{\kappa,t}(\tau)
 -(\delbar_E+t\ad\theta)\gamma_{\kappa,t}(\tau).
\]

\subsubsection{Estimate of the pairings of $\ttV_t(\tau)$}

As in the case of (\ref{eq;23.10.5.10}),
we obtain
\[
 \bigl(\ttV_t(\tau),\ttV_t(\tau)\bigr)_{L^2,h_t}
=\bigl(\ttV'_{\kappa,t}(\tau),
 \ttV'_{\kappa,t}(\tau)\bigr)_{L^2,h_t}
+O\bigl(e^{-8\kappa t}\|\tau\|_{L^2}^2\bigr).
\]

By using (\ref{eq;23.10.5.11})
with the elementary consideration in \S\ref{subsection;23.10.5.12},
we obtain
\begin{multline}
 \bigl(
 \ttV'_{\kappa,t}(\tau),
 \ttV'_{\kappa,t}(\tau)
 \bigr)_{L^2,h_t,W_{\kappa}}
=
 -2\sqrt{-1}
 \int_{W_{\kappa}}
 \Tr\bigl(
 F_{\tau}
 \cdot
 (F_{\tau})^{\dagger}_{h_{\infty}}
 \bigr)
 +O\bigl(e^{-8\kappa t}  \|\tau\|_{L^2}^2\bigr)
\\
 =-2\sqrt{-1}
 \int_{\Sigma_{\theta|W_{\kappa}}}
 \tau\overline{\tau}
  +O\bigl(e^{-8\kappa t}  \|\tau\|_{L^2}^2\bigr).
\end{multline}

By using the Stokes formula
(Lemma \ref{lem;23.10.7.20}),
and by using (\ref{eq;23.10.5.11})
with the elementary consideration in \S\ref{subsection;23.10.5.12},
we obtain 
\begin{multline}
 \bigl(
 \ttV'_{\kappa,t}(\tau),\ttV'_{\kappa,t}(\tau)
 \bigr)_{L^2,h_t,X^{(2)}_{P,\kappa}}
=
 -2\sqrt{-1}\int_{\del X^{(2)}_{P,\kappa}}
 \Tr\Bigl(
 F^{\dagger}_{\beta_P}
 \cdot
 \bigl(
 \ttV'_t(\tau)^{0,1}
 \bigr)^{\dagger}_{h_t}
 \Bigr)
\\
=-2\sqrt{-1}
  \int_{\del \Sigma_{\theta|X^{(2)}_{P,\kappa}}}
 \beta_P\cdot \taubar
 +O\bigl(e^{-8\kappa t}\|\tau\|_{L^2}^2\bigr)
 =-2\sqrt{-1}
 \int_{\Sigma_{\theta|X^{(2)}_{P,\kappa}}}
 \tau\taubar
  +O\bigl(e^{-8\kappa t}\|\tau\|_{L^2}^2\bigr).
\end{multline}
In all, we obtain
\[
 \bigl(\ttV_t(\tau),\ttV_t(\tau)\bigr)_{L^2,h_t}
 =(\tau,\tau)_{L^2}
 +O\bigl(e^{-8\kappa t}\|\tau\|_{L^2}^2\bigr)
 \quad\quad
 (0<\kappa<M(\phi)/2).
\]
This implies (\ref{eq;23.10.7.41}).

\paragraph{Acknowledgement}

This study grew out of my attempt to understand the work of
David Dumas and Andrew Neitzke \cite{Dumas-Neitzke}.
I thank Laura Fredrickson who asked
whether the method in \cite{Mochizuki-Asymptotic-Hitchin-metric}
is useful to improve the exponential rate.
I thank Rafe Mazzeo, Motohico Mulase and Olivia Dumitrescu
for some discussions.
I thank Qiongling Li and Szilard Szabo
for discussions and collaborations.
This study is partially based on the joint works
\cite{Li-Mochizuki3, Mochizuki-Szabo}.
I thank Laura Fredrickson and Andy Neitzke for encouraging comments.
I am grateful to Gao Chen and Nianzi Li for discussions.
I thank Yoshifumi Tsuchimoto and Akira Ishii
for their constant encouragements.

I prepared the original version of this manuscript
for the workshops
``Complex Lagrangians, Mirror Symmetry, and Quantization''
at Banff International Research Station
in October 2023,
and ``Higgs bundles, non-Abelian Hodge correspondence and related topics 2023''
at Chern Institute of Mathematics
in November 2023.
I revised it for the workshop
``New developments in Kobayashi-Hitchin correspondence and Higgs bundles''
at Osaka Metropolitan University
in August 2024.
I thank the organizers for the opportunities of the talks.

I am partially supported by
the Grant-in-Aid for Scientific Research (A) (No. 21H04429),
the Grant-in-Aid for Scientific Research (A) (No. 22H00094),
the Grant-in-Aid for Scientific Research (A) (No. 23H00083),
and the Grant-in-Aid for Scientific Research (C) (No. 20K03609),
Japan Society for the Promotion of Science.
I am also partially supported by the Research Institute for Mathematical
Sciences, an International Joint Usage/Research Center located in Kyoto
University.

\section{Local estimates}

\subsection{Variations of Simpson's main estimate}

We explain a variant of Simpson's main estimate
in the rank $2$ case.
We also use the functions used in \cite{Dumas-Neitzke}
to improve the estimate.

\subsubsection{A general estimate in the rank $2$ case}

For any $R>0$,
we set $U(R)=\bigl\{z\in\cnum\,\big|\,|z|<R\bigr\}$.

Let $R_0>0$.
Let $(E,\delbar_E,\theta)$ be a Higgs bundle of rank $2$ on $U(R)$
for some $R>2R_0$.
We fix an isomorphism $\kappa:\det(E)\simeq\nbigo_{U(R)}$.
Let $f$ be the endomorphism of $E$ determined by
$\theta=f\,dz$.
We assume that there exits a decomposition
\[
 (E,f)=
 (E_1,1)
 \oplus
 (E_{-1},-1).
\]
Let $\nbigh(E,\theta;\kappa)$
denote the set of harmonic metrics $h$ of
$(E,\theta)$
such that $\det(h)=1$
under the isomorphism $\kappa:\det(E)\simeq\nbigo_X$.

For any $h\in\nbigh(E,\theta;\kappa)$,
let $\pi_1$ denote the projection onto $E_1$.
We obtain the orthogonal projection $\pi_1^h$ of $E$
onto $E_1$.
We regard
$\pi_1$ and $\pi_1^h$ as endomorphisms of $E$
in natural ways.
We set $\rho:=\pi_1-\pi_1^h$.
We have
$|\pi_1|_h^2
=|\pi_1^h|_h^2+|\rho|_h^2$.
We shall prove the following proposition
in \S\ref{subsection;23.9.28.1}.
\begin{prop}
\label{prop;23.9.27.20}
For any $0<\gamma<2\sqrt{2}$,
there exists $C(R_0,\gamma)>0$ depending only on $R_0$ and $\gamma$,
such that the following holds on $U(R-R_0)$
for any $h\in\nbigh(E,\theta;\kappa)$: 
\[
 |\rho|_h
 \leq
 C(R_0,\gamma)e^{-\gamma (R-|z|)}.
\] 
\end{prop}

\subsubsection{Refinement under an additional assumption}

Let $v_i$ $(i=1,2)$ be frames of $E_i$
such that $\kappa(v_1\wedge v_2)=1$.
We shall prove the following proposition in
\S\ref{subsection;24.7.5.1}.

\begin{prop}
\label{prop;23.9.27.10}
Let $h\in \nbigh(E,\theta;\kappa)$.
Suppose that there exist $\gamma_0>2$ and $B(h)>0$ 
such that the following holds on $U(R)$:
\begin{equation}
\label{eq;23.9.27.11}
 h(v_1,v_1)+h(v_2,v_2)
+|\del_z h(v_1,v_1)|
+|\del_z h(v_2,v_2)|
+|\del_{z}h(v_1,v_2)|
+|\del_{\zbar}h(v_1,v_2)|
\leq B(h) e^{-\gamma_0 (R-|z|)}.
\end{equation}
Then, 
for any $0<\gamma<4$,
there exist $C(R_0,\gamma,\gamma_0,B(h))>0$
depending only on $R_0,\gamma,\gamma_0,B(h)$,
such that the following holds on $U(R-R_0)$:
\[
 |\rho|_h\leq
 C\bigl(R_0,\gamma,\gamma_0,B(h)\bigr)
 e^{-\gamma(R-|z|)}.
\]
\end{prop}

\subsubsection{Refinement in the symmetric case}

Let $C$ be a non-degenerate symmetric pairing of $(E,\theta)$
such that $\det(C)=1$ under the isomorphism $\kappa$.
Let $\nbigh(E,\theta;C)$ denote the set of
harmonic metrics of $(E,\theta)$ compatible with $C$.
For any $h\in\nbigh(E,\theta;C)$,
we obtain an estimate as in Proposition \ref{prop;23.9.27.10}
without the assumption (\ref{eq;23.9.27.11}).
We shall prove the following proposition
in \S\ref{subsection;23.9.28.10}.

\begin{prop}
\label{prop;23.9.27.110}
For any $0<\gamma<4$,
there exists $C(R_0,\gamma)>0$
such that the following holds on $U(R-R_0)$
for any $h\in\nbigh(E,\theta;C)$: 
\[
 |\rho|_h
 \leq
 C(R_0,\gamma)e^{-\gamma (R-|z|)}.
\] 
\end{prop}

Let $v_i$ be frames of $E_i$
such that $C(v_i,v_i)=1$.
The following proposition
is a refinement of Proposition \ref{prop;23.9.27.110}.

\begin{prop}
\label{prop;23.9.27.111}
For any $0<\gamma<4$
and $(j,k)\in\seisuu^2_{\geq 0}$,
there exists $C_{j,k}(R_0,\gamma)>0$
such that the following holds
for any $h\in\nbigh(E,\theta;C)$
on $U(R-R_0)$:
\begin{equation}
\label{eq;23.9.28.20}
 |\del_z^j\del_{\zbar}^k(h(v_i,v_i)-1)|
 \leq
 C_{j,k}(R_0,\gamma)e^{-2\gamma(R-|z|)},
\end{equation}
\begin{equation}
\label{eq;23.9.28.21}
 |\del_z^j\del_{\zbar}^kh(v_1,v_2)|
 \leq
 C_{j,k}(R_0,\gamma)e^{-\gamma(R-|z|)}.
\end{equation}
\end{prop}

Let $h^C$ denote the unique decoupled harmonic metric of $(E,\theta)$
compatible with $C$
(see \cite{Li-Mochizuki3}).
We have $h^C(v_i,v_i)=1$ and $h^C(v_1,v_2)=0$.
Let $\nabla$ denote the Chern connection of $E$
determined by $h^C$.
For any $h\in \nbigh(E,\theta;C)$,
let $s(h^C,h)$ be the automorphism of $E$
determined by $h=h^C\cdot s$.
We set
$\End(E)^{\circ}=\End(E_1)\oplus\End(E_2)$
and
$\End(E)^{\bot}=\Hom(E_1,E_2)\oplus\Hom(E_2,E_1)$.
We obtain the decomposition
\[
 \End(E)=\End(E)^{\circ}\oplus\End(E)^{\bot}.
\]
We have the corresponding decomposition
$s(h^C,h)=s(h^C,h)^{\circ}+s(h^C,h)^{\bot}$.

\begin{cor}
\label{cor;23.9.28.40}
For any $0<\gamma<4$
and $(j,k)\in\seisuu^2_{\geq 0}$,
there exists $C'_{j,k}(R_0,\gamma)>0$
such that the following holds
for any $h\in\nbigh(E,\theta;C)$
on $U(R-R_0)$:
\[
 |\nabla_z^j\nabla_{\zbar}^k(s(h^C,h)^{\circ}-\id)|
 \leq
 C'_{j,k}(R_0,\gamma)
 e^{-2\gamma(R-|z|)},
\]
\[
   |\nabla_z^j\nabla_{\zbar}^k(s(h^C,h)^{\bot})|
 \leq
 C'_{j,k}(R_0,\gamma)
 e^{-\gamma(R-|z|)}.
\]
\hfill\qed
\end{cor}

\begin{rem}
Proposition {\rm\ref{prop;23.9.27.110}}
and Proposition {\rm\ref{prop;23.9.27.111}}
are essentially contained in 
{\rm\cite{Dumas-Neitzke}}.
They can be proved independently from
Proposition {\rm\ref{prop;23.9.27.20}}
and Proposition {\rm\ref{prop;23.9.27.10}}.
\hfill\qed 
\end{rem}

\subsubsection{A general inequality}

Let $h\in\nbigh(E,\theta;\kappa)$.
We have the following equality:
\begin{equation}
\label{eq;23.9.27.101}
 -\del_z\del_{\zbar}
 |\pi_1|_{h}^2
=-|\del_{E,h}\pi_1|^2_h
-\bigl|[f^{\dagger}_h,\pi_1]\bigr|_h^2.
\end{equation}

Let $v_i$ $(i=1,2)$ be frames of $E_i$
such that $v_1\wedge v_2=1$.
Let $H$ be the Hermitian matrix valued function
determined by
$H_{i,j}=h(v_i,v_j)$.
To simplify,
we set
$a=h(v_1,v_1)$,
$b=h(v_1,v_2)$
and $c=h(v_2,v_2)$,
i.e.,
\[
 H=\left(
 \begin{array}{cc}
  a & b \\ \bbar & c
 \end{array}
 \right).
\]
Because $\det(h)=1$,
we have $\det(H)=ac-|b|^2=1$.
We also have $|\rho|_h^2=|b|^2$.

We shall use the following lemma implicitly.
\begin{lem}
Let $g$ be an endomorphism of $E$.
Let $g^{\dagger}_h$ be the adjoint of $g$ with respect to $h$.
Let $G$ and $G^{\dagger}_H$ denote the matrix valued functions
determined by
$g\vecv=\vecv G$  and $g^{\dagger}_h\vecv=\vecv G^{\dagger}_H$.
Then, we have
$G^{\dagger}_H
=\Hbar^{-1}\cdot (\lefttop{t}\Gbar)\cdot \Hbar$.
\hfill\qed
\end{lem}

\begin{lem}
\label{lem;23.9.24.1}
We have
$\bigl|[f^{\dagger}_h,\pi_1]\bigr|_h^2
=8(1+|b|^2)|b|^2$.
\end{lem}
\pf
Let $F$ and $F^{\dagger}_{H}$
be the matrix valued functions determined by
$f\vecv=\vecv F$
and
$f^{\dagger}_h\vecv=\vecv\cdot F^{\dagger}_H$.
More explicitly,
\[
F=\left(
\begin{array}{cc}
 1 & 0 \\ 0 & -1
\end{array}
\right),\quad
F^{\dagger}_H
=\left(
 \begin{array}{cc}
  c& -\bbar\\
  -b & a
 \end{array}
 \right)
 \left(
\begin{array}{cc}
 1 & 0 \\ 0 & -1
\end{array}
\right)
\left(
 \begin{array}{cc}
 a & \bbar \\
 b & c
 \end{array}
 \right)
 =\left(
 \begin{array}{cc}
  1+2|b|^2 &
 2\bbar c
 \\
  -2ab & -1-2|b|^2
 \end{array}
 \right).
\]
Let $\Pi_1$ be the matrix valued function determined by
$\pi_1\vecv=\vecv\cdot \Pi_1$,
i.e.,
\[
 \Pi_1=\left(
 \begin{array}{cc}
  1 & 0 \\ 0 & 0
 \end{array}
 \right).
\]

We can check the following:
\[
 [F^{\dagger}_H,\Pi_1]
=\left(
 \begin{array}{cc}
  0 & -2\bbar c\\ -2ab & 0
 \end{array} 
 \right)
 =[F^{\dagger}_H,\Pi_1]^{\dagger}_{H}.
\]
We obtain
\[
  [F^{\dagger}_H,\Pi_1]
  \cdot
  [F^{\dagger}_H,\Pi_1]^{\dagger}_H
=\left(
 \begin{array}{cc}
  4(1+|b|^2) |b|^2 & 0\\
  0 & 4(1+|b|^2)|b|^2
 \end{array}
 \right).
\]
Hence, we obtain
\[
\bigl|
[f^{\dagger}_h,\pi_1]
\bigr|_h^2
=\Tr\bigl(
[F^{\dagger}_H,\Pi_1]^{\dagger}_H
 \cdot
 [F^{\dagger}_H,\Pi_1]
 \bigr)
=8(1+|b|^2)|b|^2.
\]
Thus, we obtain Lemma \ref{lem;23.9.24.1}.
\hfill\qed

\begin{cor}
We obtain
$-\del_z\del_{\zbar}
 |b|_h^2
 \leq -8|b|_h^2$.
\hfill\qed
\end{cor}

\subsubsection{Bessel function of the first kind}

We recall
some functions by following 
\cite[\S5]{Dumas-Neitzke}.
Let $I_0(z):\cnum\to\real_{>0}$
be the function as in \cite[\S5]{Dumas-Neitzke},
i.e.,
$I_0$ denotes the unique positive even $C^{\infty}$-function
such that
(i) $I_0(z)$ depends only on $|z|$,
(ii) $(\del_x^2+\del_y^2)I_0=I_0$,
(iii) $I_0(0)=1$.
We have the convergence
$I_0(z)\cdot
(2\pi|z|)^{1/2}
 e^{-|z|}
 \to 1$
 as $|z|\to\infty$.
In particular, we have the following lemma.
\begin{lem}
For any $0<u<1$,
there exists $C_0(u)>0$
such that the following holds on $\cnum$:
\[
 I_0(z)^{-1}
 \leq
 C_0(u)e^{-u|z|}.
\]
\hfill\qed
\end{lem}

\begin{lem}[\cite{Dumas-Neitzke}]
For any $0<\gamma_1<\gamma_2$,
there exists $A(\gamma_1,\gamma_2)>0$
such that the following holds for any $0<b<a$:
\[
 e^{-\gamma_2(a-b)}
<
 A(\gamma_1,\gamma_2)
 \frac{I_0(\gamma_1b)}{I_0(\gamma_1a)}
\]
\hfill\qed
 \end{lem}

\subsubsection{Proof of Proposition \ref{prop;23.9.27.20}}
\label{subsection;23.9.28.1}

For any $r>0$ and $z_0\in \cnum$,
we set $U(z_0,r)=\bigl\{z\in\cnum\,\big|\,|z-z_0|<r\bigr\}$.
We set $r_0=R_0/10$.
We recall the following general lemma.
\begin{lem}
\label{lem;23.9.27.30}
There exists $C_1(r_0)$ depending only on $r_0$
such that the following holds:
\begin{itemize}
 \item  Let $h$ be any harmonic metric of
	$(E,\theta)_{|U(z_0,r_0)}$.
	Then,
	$|f|_h\leq C_1(r_0)$ holds
	on $U(z_0,r_0/2)$.
	\hfill\qed
\end{itemize}
\end{lem}
\begin{cor}
$|f|_h\leq C_1(r_0)$ holds 
on $U(R-r_0/2)$ for any $h\in\nbigh(E,\theta)$.
\hfill\qed 
\end{cor}

Note that $|f|_h^2=2(1+2|b|^2)$.
\begin{cor}
There exists $C_2(r_0)>0$
such that $|b|\leq C_2(r_0)$ on $U(R-r_0/2)$.
\end{cor}

\begin{lem}
The following hols:
\[
 |b(0)|^2
 \leq
 2C_2(r_0)^2I_0(4\sqrt{2}(R-r_0/2))^{-1}.
\]
\end{lem}
\pf
We set
$J(z)=2C_2(r_0)^2I_0(4\sqrt{2} z)I_0(4\sqrt{2}(R-r_0/2))^{-1}$.
We have
\[
 -\del_z\del_{\zbar}J(z)
 =-8J(z),
 \quad
J(R-r_0/2)=2C_2(r_0)^2.
\]
Let $Z=\{z\in U(R-r_0/2)\,|\,|b(z)|^2>J(z)\}$.
We deduce a contradiction by assuming $Z\neq\emptyset$.
It is relatively compact in $U(R-r_0/2)$.
We have $|b(z)|^2=J(z)$ on $\del Z$.
On $Z$, we have
\[
 -\del_z\del_{\zbar}(|b|^2-J)
 \leq -8(|b|^2-J)
 <0.
\]
Hence, $|b|^2-J$ is subharmonic on $Z$.
By the maximum principle,
we obtain $|b|^2-J\leq 0$, which contradicts with
the definition of $Z$.
Hence, we obtain $Z=\emptyset$.
In particular,
$|b(0)|^2\leq J(0)$.
\hfill\qed

\begin{cor}
\label{cor;23.9.27.50}
For any $0<u<1$, 
\[
 |b(0)|^2
 \leq
 2C_2(r_0)^2C_0(u)
 e^{-4\sqrt{2}(R-r_0/2)u}
\leq 2C_2(r_0)^2C_0(u)
 e^{2\sqrt{2}r_0}
 e^{-4\sqrt{2}uR}.
\]
\hfill\qed
\end{cor}

For $z_0\in U(R-R_0)$,
we have $U(z_0,R-|z_0|)\subset U(R)$.
By applying Corollary \ref{cor;23.9.27.50},
we obtain
\[
 |b(z_0)|^2
 \leq
  2C_2(r_0)^2C_0(u)
 e^{-4\sqrt{2}(R-|z_0|-r_0/2)u}
 \leq 2C_2(r_0)^2
 C_0(u)
 e^{2\sqrt{2}r_0}
 e^{-4\sqrt{2}u(R-|z_0|)}.
\]
It implies the claim of Proposition \ref{prop;23.9.27.20}.
\hfill\qed

\subsubsection{Proof of Proposition \ref{prop;23.9.27.10}}
\label{subsection;24.7.5.1}

\begin{lem}
We obtain
\begin{multline}
\label{eq;23.9.27.100}
\bigl|
 \del_{E,h,z}\pi
 \bigr|_h^2
=(1+|b|^2)^2
 \bigl(|\del_{\zbar}b|^2+|\del_zb|^2\bigr)
+|b|^2(a^2|\del_zc|^2+c^2|\del_za|^2)
 \\
 -2\Re\Bigl(
 ba^2c\del_z\bbar\del_{\zbar}c
+\bbar ac^2\del_zb\del_{\zbar}a
 +b^2\bigl(
 -c\del_z\bbar+\bbar\del_zc
 \bigr)
 \bigl(-\bbar\del_{\zbar}a+a\del_{\zbar}\bbar\bigr)
 \Bigr).
\end{multline}
\end{lem}
\pf
We have
\[
 \Hbar^{-1}\del_z\Hbar
=\left(
 \begin{array}{cc}
 c & -\bbar \\ -b & a
 \end{array}\right)
\left(
 \begin{array}{cc}
 \del_za & \del_z\bbar \\ \del_zb & \del_zc
 \end{array}\right)
 =\left(
\begin{array}{cc}
 c\del_za-\bbar\del_zb&
 c\del_z\bbar-\bbar\del_zc
 \\
 -b\del_za+a\del_zb
  &
  -b\del_z\bbar+a\del_zc
\end{array}
 \right).
\]
Let $A$ be determined by
$(\del_{E,h,z}\pi_1)\vecv=\vecv\cdot A$.
We set
$\beta=-c\del_z\bbar+\bbar\del_zc$
and $\gamma=-b\del_za+a\del_zb$.
We obtain
\[
 A=\del_z\Pi_1
 +[\Hbar^{-1}\del_z\Hbar,\Pi_1]
 =\left(
 \begin{array}{cc}
  0 & \beta
   \\
  \gamma
   & 0
 \end{array}
 \right).
\]
We have
\[
 A\cdot \Hbar^{-1}\lefttop{t}\Abar\cdot \Hbar
 =\left(
 \begin{array}{cc}
  -b^2\beta\gammabar+a^2|\beta|^2&
   -bc\beta\gammabar+a\bbar|\beta|^2
   \\
  bc|\gamma|^2-a\bbar\gamma\betabar
   &
   |\gamma|^2c^2-\gamma\betabar\bbar^2
 \end{array}
 \right).
\]

We have
\[
 \bigl|
 \del_{E,h,z}\pi_1
 \bigr|_h^2
 =
 \Tr\bigl(
 A\cdot \Hbar^{-1}\lefttop{t}\Abar\cdot \Hbar
 \bigr)
 =a^2|\beta|^2+c^2|\gamma|^2
 -2\Re(b^2\beta\gammabar).
\]
It is equal to
\begin{multline}
 a^2\Bigl(
 c^2|\del_z\bbar|^2
 +|b|^2|\del_zc|^2
 -2\Re(bc\del_z\bbar\del_{\zbar}c)
 \Bigr)
 +c^2\Bigl(
 |b|^2|\del_za|^2
 +|a|^2|\del_zb|^2
-2\Re(\bbar a\del_zb\del_{\zbar}a)
 \Bigr)
 \\
 -2\Re\Bigl(
 b^2(-c\del_z\bbar+\bbar\del_zc)
 (-\bbar\del_{\zbar}a+a\del_{\zbar}\bbar)
 \Bigr).
\end{multline}
Then,
by using $ac=1+|b|^2$,
we obtain (\ref{eq;23.9.27.100}).
\hfill\qed

\vspace{.1in}
We set $W=\{z\in U(R)\,|\,b(z)\neq 0\}$.
We take $4<\gamma_1<\gamma_0+\min\{2\sqrt{2},\gamma_0\}$.

\begin{lem}
There exists $C_{20}(B(h),\gamma_1)>0$,
depending only on $B(h)$ and $\gamma_1$
such that the following holds 
on $W$:
\begin{equation}
\label{eq;23.9.27.102}
 -\del_z\del_{\zbar}|b|
 \leq
 -4|b|
 +C_{20}(B(h),\gamma_1) e^{-\gamma_1(R-|z|)}.
\end{equation}
\end{lem}
\pf
On $W$,
we have
\[
 \del_z|b|
=\del_{\zbar}|b|
 \leq
 \frac{1}{2}
 \bigl(|\del_zb|+|\del_{\zbar}b|\bigr).
\]
We obtain
\[
 2\del_z|b|
 \del_{\zbar}|b|
 \leq
 \frac{1}{2}
 \bigl(|\del_zb|+|\del_{\zbar}b|\bigr)^2
 \leq
 |\del_zb|^2+|\del_{\zbar}b|^2.
\]
Then, we obtain 
(\ref{eq;23.9.27.102})
from
(\ref{eq;23.9.27.101})
and (\ref{eq;23.9.27.100}).
\hfill\qed

\vspace{.1in}

Let $0<\gamma<4$.
We set $\eta=(4+\gamma)/2$.
We have
\[
 -\del_z\del_{\zbar}I_0(\eta z)
 =-\frac{1}{4}\eta^2I_0(\eta z).
\]
There exists $C_{21}(B(h),\gamma_1,\gamma)$
such that
$|b|<C_{21}(B(h),\gamma_1,\gamma)$ on $|z|=R-r_0$,
and that the following holds
for        $g=C_{21}(B(h),\gamma_1,\gamma)
       I_0(\eta z)I_0(\eta (R-r_0))^{-1}$
on $U(R-r_0)$:
\[
       -\del_z\del_{\zbar}g
       \geq
       -4g
       +C_{20}(B(h),\gamma_1)e^{-\gamma_1(R-|z|)}.
\]

\begin{lem}
We have $|b|\leq g$ on $U(R-r_0)$.
In particular,
$|b(0)|\leq C_{21}(B(h),\gamma_1,\gamma)I_0(\eta(R-r_0))^{-1}$.
\end{lem}
\pf
Let $Z=\{z\in U(R-r_0)\,|\,|b(z)|>g(z)\}$.
We shall deduce a contradiction by assuming $Z\neq\emptyset$.
We note that $Z\cap \del U(R-r_0)=\emptyset$.
We also have $Z\subset W$.
Hence, we have $|b|=g$ on $\del Z$,
and the following inequality on $Z$:
\[
 -\del_z\del_{\zbar}(|b|-g)
 \leq
 -4(|b|-g).
\]
Then, we obtain $|b|\leq g$ on $Z$
which is a contradiction.
\hfill\qed

\vspace{.1in}
We obtain
\[
 |b(0)|
 \leq
 C_{21}(B(h),\gamma,\gamma_1)
 C_0(\gamma/\eta)
 e^{\gamma r_0}
 e^{-\gamma R}.
\]
By applying the argument to
the restriction of $(E,\theta,h)$
to $U(z_0,R-|z_0|)$ for $z_0\in U(R_1)$,
we obtain
\[
  |b(z_0)|
 \leq
 C_{21}(B(h),\gamma,\gamma_1)
 C_0(\gamma/\eta)
 e^{\gamma r_0}
 e^{-\gamma (R-|z_0|)}.
\]
Thus, we obtain Proposition \ref{prop;23.9.27.10}.
\hfill\qed

\subsubsection{Proof of Proposition \ref{prop;23.9.27.110}
and the $C^0$-part of Proposition \ref{prop;23.9.27.111}}
\label{subsection;23.9.28.10}

Let $v_1,v_2$ be a frame such that
$C(v_i,v_i)=1$ and $C(v_1,v_2)=0$.
Let $h\in\nbigh(E,\theta;C)$.
Let $H$ be the Hermitian matrix valued function
determined by $H_{i,j}=h(v_i,v_j)$.
We can describe $H$ as 
\[
 H=\left(
 \begin{array}{cc}
  a & \sqrt{-1}\btilde \\ -\sqrt{-1}\btilde & a
 \end{array}
 \right),
\]
where $a$ is an $\real_{>0}$-valued function,
and $\btilde$ is an $\real$-valued function
such that $a^2-\btilde^2=1$.

\begin{lem}
The following holds:
\[
 |\del_{E,h,z}\pi_1|^2_h
 =\frac{2(1+2\btilde^2)}{1+\btilde^2}
 |\del_z\btilde|^2.
\]
\end{lem}
\pf
We set 
\[
 \beta=-\sqrt{-1}\btilde\del_za+\sqrt{-1}a\del_z\btilde
=\sqrt{-1}a^{-1}\del_z\btilde.
\]
We have
\[
\Hbar^{-1}\del_z\Hbar
=
\left(
\begin{array}{cc}
 0& -\beta\\
\beta 
  & 0
\end{array}
\right).
\]
Then, we obtain
\[
 \bigl[
 \del_z+\Hbar^{-1}\del_z\Hbar,
 \Pi_1
 \bigr]
 =\left(
 \begin{array}{cc}
  0 & \beta \\ \beta & 0
 \end{array}
 \right).
\]
We obtain
\begin{multline}
 \bigl|
 \del_{E,h,z}\pi_1
 \bigr|^2_h
 =\tr\left(
 \left(
 \begin{array}{cc}
  0 & \beta \\ \beta & 0
 \end{array}
 \right)
 \left(
 \begin{array}{cc}
  a & \sqrt{-1}\btilde \\ -\sqrt{-1}\btilde & a
 \end{array}
 \right)
 \left(
 \begin{array}{cc}
  0 & \betabar \\ \betabar & 0
 \end{array}
 \right)
 \left(
 \begin{array}{cc}
  a & -\sqrt{-1}\btilde \\ \sqrt{-1}\btilde & a
 \end{array}
 \right)
 \right) 
\\
 =2(1+2\btilde^2)|\beta|^2
 =\frac{2(1+2\btilde^2)}{1+\btilde^2}
 |\del_z\btilde|^2.
\end{multline}
\hfill\qed

\vspace{.1in}
We set $W=\{z\in U(R)\,|\,\btilde(z)\neq 0\}$.

\begin{lem}
On $W$, the following holds:
\[
-\del_z\del_{\zbar}|\btilde|
\leq -4|\btilde|.
\]
\end{lem}
\pf
We have
\[
-\del_z\del_{\zbar}\btilde^2
=-2|\btilde|\del_z\del_{\zbar}|\btilde|
-2(\del_z|\btilde|)\cdot(\del_{\zbar}|\btilde|).
\]
Because $\btilde$ is real valued,
we obtain
\[
 \del_z|\btilde|
=\del_{\zbar}|\btilde|
\leq |\del_z\btilde|.
\]
We have
\[
-2\frac{1+2\btilde^2}{1+\btilde^2}
 |\del_z\btilde|^2
 +2\del_z|\btilde|\del_{\zbar}|\btilde|
\leq 0.
\]
Then, we obtain
$-2|\btilde|\del_z\del_{\zbar}|\btilde|
 \leq
 -8|\btilde|^2$,
which implies the claim of the lemma.
\hfill\qed

\vspace{.1in}

\begin{prop}
\label{prop;23.9.28.3}
For any $0<\gamma<4$,
there exists $C'(R_0,\gamma)>0$
such that the following holds
for any $h\in\nbigh(E,\theta;C)$
on $U(R-R_0)$:
\[
 |\btilde(z)|\leq C'(R_0,\gamma)e^{-\gamma (R-|z|)}.
\]
It implies
$|a-1|\leq C'(R_0,\gamma)e^{-2\gamma (R-|z|)}$.
\end{prop}
\pf
We use the notation in \S\ref{subsection;23.9.28.1}.
\begin{lem}
The following holds:
\[
 |\btilde(0)|
 \leq
 2C_2(r_0)
 I_0(4(R-r_0/2))^{-1}.
\]
\end{lem}
\pf
We set $J(z)=2C_2(r_0)I_0(4z)I_0(4(R-r_0/2))^{-1}$.
\[
 -\del_z\del_{\zbar}J(z)
 =-4J(z),
 \quad
 J(z)_{|z|=R-r_0/2}=2C_2(r_0).
\]
Let $Z=\{z\in U(R-r_0/2)\,|\,|\btilde(z)|>J(z)\}$.
We deduce a contradiction by assuming $Z\neq\emptyset$.
Clearly, $Z$ is contained in $W$,
and relatively compact in $U(R-r_0/2)$.
We have $|\btilde(z)|=J(z)$ on $\del Z$,
and 
\[
 -\del_z\del_{\zbar}(|\btilde|-J)
 \leq -4(|\btilde|^2-J)
 <0
\]
on $Z$.
Hence, $|\btilde|-J$ is subharmonic on $Z$.
By the maximum principle,
we obtain $|\btilde|-J\leq 0$, which contradicts with
the definition of $Z$.
Hence, we obtain $Z=\emptyset$.
In particular,
$|\btilde(0)|\leq J(0)$.
\hfill\qed

\begin{cor}
\label{cor;23.9.28.2}
For any $0<u<1$, we have
\[
  |\btilde(0)|
 \leq
 2C_2(r_0)C_0(u)
 e^{-4(R-r_0/2)u}
\leq 2C_2(r_0)C_0(u)
 e^{2r_0}
 e^{-4uR}.
\]
\hfill\qed
\end{cor}

By applying Corollary \ref{cor;23.9.28.2}
to the restriction to $U(z_0,R-|z_0|)$,
we obtain
\[
 |\btilde(z_0)|
 \leq
 2C_2(r_0)C_0(u)
 e^{-4(R-r_0/2)u}
\leq 2C_2(r_0)C_0(u)
 e^{2r_0}
 e^{-4u(R-|z_0|)}.
\]
It implies the claim of Proposition \ref{prop;23.9.28.3}.
\hfill\qed

\vspace{.1in}
Proposition \ref{prop;23.9.28.3}
implies
the claim of Proposition \ref{prop;23.9.27.111}
in the case $(j,k)=(0,0)$.
It implies Proposition \ref{prop;23.9.27.110}.

\subsubsection{Proof of Proposition \ref{prop;23.9.27.111}}

We continue to use the notation in \S\ref{subsection;23.9.28.10}.
We have
\[
 \del_{\zbar}
 \Bigl(
 \Hbar^{-1}\del_z\Hbar
 \Bigr)
= \left(
 \begin{array}{cc}
  0 & -\sqrt{-1}\del_{\zbar}(a^{-1}\del_z\btilde) \\
 \sqrt{-1}\del_{\zbar}(a^{-1}\del_z\btilde) & 0
 \end{array}
 \right).
\]
We have 
\[
 F^{\dagger}_H=
\Hbar^{-1}F
 \Hbar
 =\left(
 \begin{array}{cc}
  1+ \btilde^2 & -2\sqrt{-1}a\btilde \\
  -2\sqrt{-1}a\btilde & -1-2\btilde^2
 \end{array}
 \right).
\]
We obtain
\[
[F,F^{\dagger}_H]
=\left(
 \begin{array}{cc}
  0 & -4\sqrt{-1}a\btilde \\
  4\sqrt{-1}a\btilde & 0
 \end{array}
\right).
\] 
Because
$\del_{\zbar}
 \left(
 \Hbar^{-1}\del_{z}\Hbar
 \right)
-[F,F^{\dagger}_H]=0$,
we obtain
\begin{equation}
\label{eq;23.9.28.11}
 \del_{\zbar}\bigl(
 a^{-1}\del_{z}\btilde
 \bigr)
-4a\btilde=0.
\end{equation}
By using $a^2-\btilde^2=1$,
we obtain
\begin{equation}
\label{eq;23.9.28.10}
 \del_{\zbar}\del_z\btilde
=\frac{1}{1+\btilde^2}
 \btilde|\del_z\btilde|^2
 +4\btilde.
\end{equation}

For $z_0\in U(R_1)$,
let $\chi_{z_0}:\cnum\to [0,1]$ be a $C^{\infty}$-function
such that
$\chi_{z_0}(z)=1$ $(|z-z_0|<r_0/2)$
and
$\chi_{z_0}(z)=0$ $(|z-z_0|\geq 2r_0/3)$.
We also assume that
$\chi_{z_0}^{-1/2}d\chi_{z_0}$ induces a $C^{\infty}$-function.
By (\ref{eq;23.9.28.10}),
we obtain
\[
0\leq
 \int \chi_{z_0}\btilde\del_{\zbar}\del_z\btilde
=-\int \del_{\zbar}(\chi_{z_0}) \btilde\del_z\btilde
-\int \chi_{z_0}|\del_{z}\btilde|^2.
\]

We obtain
\[
 \int \chi_{z_0}|\del_z\btilde|^2
 \leq
 \left(
 \int |\chi_{z_0}^{-1/2}d\chi_{z_0}|
 \cdot |\btilde|^2
 \right)^{1/2}
\left(
 \int \chi_{z_0}|\del_z\btilde|^2
 \right)^{1/2}.
\]
Hence,
\[
  \int \chi_{z_0}|\del_z\btilde|^2
 \leq
 \int |\chi_{z_0}^{-1/2}d\chi_{z_0}|\cdot
 |\btilde|^2.
\]
For any $0<\gamma<4$,
we obtain
\[
\left(
 \int_{|z-z_0|<r_0/2}
 |\del_z\btilde|^2
 \right)^{1/2}
 \leq
 \pi^{1/2} (r_0/2)
 C'(R,R_1+r_0,\gamma)
 e^{\gamma r_0}
 e^{-\gamma(R-|z_0|)}.
\]
By using (\ref{eq;23.9.28.11}),
we obtain the estimate for the sup norms of
$\del_z\btilde$ on $\{|z-z_0|<r_0/3\}$.
By a standard bootstrapping argument,
we obtain the estimate (\ref{eq;23.9.28.21}).
Because $a^2=1+\btilde^2$,
it implies the estimate (\ref{eq;23.9.28.20}).
\hfill\qed

\subsection{An estimate for derivatives}

Fix $R_0>0$ and $C_0>0$.
Let $(E,\theta)$ be a Higgs bundle of rank $r$ on $U(R)$
for $R>2R_0$.
Let $f$ be the endomorphism determined by $\theta=f\,dz$.
Assume that
$|\tr(f^j)|<C_0$ for $j=1,\ldots,r$.
Let $h_0$ be a harmonic metric.
Let $g_0$ be the Euclidean metric $dz\,d\zbar$.
Let $\nabla$ denote the Chern connection of $(E,h_0)$.
Let $\gbigr$ denote the endomorphism of $E$
such that $\gbigr\,dz\,d\zbar$ is the curvature of $\nabla$.

\begin{lem}
\label{lem;23.9.30.20}
There exist $C(j,k)>0$ depending only on
$R_0$, $C_0$, $r$ and $(j,k)$
such that the following holds on $U(R-R_0)$
for any harmonic metric $h$ of $(E,\theta)$:
\[
 |\nabla_z^jf|_{h_0}
+|\nabla_{\zbar}^kf^{\dagger}_{h_0}|_{h_0}
 +|\nabla_z^j\nabla_{\zbar}^k\gbigr|_{h_0}\leq
 C(j,k). 
\] 
\end{lem}
\pf
For $z_0\in U(R-R_0)$,
we have $U(z_0,R_0)\subset U(R)$.
By applying \cite[Proposition 4.1]{Mochizuki-Szabo}
to $(E,\theta)_{|U(z_0,R_0)}$,
we obtain the estimate on $U(z_0,R_0/2)$.
Hence, we obtain the desired estimate on $U(R-R_0)$.
\hfill\qed

\vspace{.1in}

For positive constants $C_1,\epsilon_0,\epsilon_1$,
let $\nbigh(E,\theta;h_0,C_1,\epsilon_0,\epsilon_1)$
denote the set of harmonic metrics $h$ of
$(E,\theta)$
such that the following holds on $U(R)$.
\begin{itemize}
 \item
      $|s(h_0,h)|_{h_0}+|s(h_0,h)^{-1}|_{h_0}\leq C_1$.
 \item $\int \bigl|\delbar s(h_0,h)\bigr|_{h_0}^2
       +\int \bigl|[\theta,s(h_0,h)]\bigr|_{h_0}^2\leq \epsilon_1^2$.
 \item $\Bigl|
       \delbar\bigl(s(h_0,h)^{-1}\del_{h_0}s(h_0,h)\bigr)
       \Bigr|_{h_0,g_0}\leq \epsilon_0$.
\end{itemize}
We set $\epsilon_2=\max\{\epsilon_0,\epsilon_1\}$.

\begin{lem}
\label{lem;23.9.28.110}
There exists
$C_2$ depending only on 
$R_0,C_0,C_1,r,\epsilon_0,\epsilon_1$,
such that the following holds on $U(R-R_0)$
for any $h\in \nbigh(E,\theta;h_0,C_1,\epsilon_0,\epsilon_1)$:
\[
 \bigl|\nabla s(h_0,h)\bigr|_{h_0,g_0}
+\bigl|[\theta,s(h_0,h)]\bigr|_{h_0,g_0}
 \leq C_2 \cdot \epsilon_2. 
\] 
\end{lem}
\pf
To simplify the description,
we denote $s(h_0,h)$ by $s$.
We obtain the estimate for
$\sup_{U(R-R_0)}\bigl|
\delbar s
\bigr|_{h_0,g_0}$
by using the Sobolev embedding theorem
and Uhlenbeck's theorem on the existence of Coulomb gauge.
(See the proof of \cite[Proposition 4.1]{Mochizuki-Szabo},
for example.)
Because $s$ is self-adjoint with respect to $h_0$,
we obtain the estimate for
$\sup_{U(R-R_0)}\bigl|
\nabla s
\bigr|_{h_0,g_0}$.
Because 
$\delbar[\theta,s]
=-[\theta,\delbar s]$,
we also obtain the estimate for
$\sup_{U(R-R_0)}\bigl|
[\theta,s]
\bigr|_{h_0,g_0}$.
\hfill\qed

\begin{lem}
\label{lem;23.9.28.111}
For any $(j,k)\in\seisuu_{\geq 0}^2\setminus\{(0,0)\}$,
there exist positive constants
$C_3(j,k)$ depending only on
$R_0$, $C_0$, $C_1$,
$r$, $\epsilon_0$ and $\epsilon_1,(j,k)$,
such that the following holds on $U(R-R_0)$:
\[
\bigl|\nabla_z^j\nabla_{\zbar}^k s(h_0,h)\bigr|_{h_0}
\leq C_3(j,k) \cdot \epsilon_2.
\] 
\end{lem}
\pf
We have
\[
 \delbar\bigl(
 s(h_0,h)^{-1}
 \del_{E,h_0}s(h_0,h)
 \bigr)
 +\bigl[
 \theta,
 s^{-1}[\theta^{\dagger}_{h_0},s]
 \bigr]
=0.
\]
Then, we obtain the desired estimate
by using a standard bootstrapping argument,
the Sobolev embedding theorem
and Uhlenbeck's theorem on the existence of Coulomb gauge.
\hfill\qed

\begin{cor}
\label{cor;24.7.6.11}
For any $(j,k)\in\seisuu_{\geq 0}^2$,
there exist positive constants
$C_4(j,k)$ depending only on
$R_0$, $C_0$, $C_1$,
$r$, $\epsilon_0$ and $\epsilon_1,(j,k)$,
such that the following holds on $U(R-R_0)$:
\[
 \bigl|\nabla_z^j\nabla_{\zbar}^k
 [\theta,s(h_0,h)]\bigr|_{h_0}
+\bigl|\nabla_z^j\nabla_{\zbar}^k
 [\theta^{\dagger}_{h_0},s(h_0,h)]\bigr|_{h_0}
 \leq C_4(j,k) \cdot \epsilon_2.
\] 
\hfill\qed
\end{cor}

\section{Global estimates}

\subsection{Flat metrics induced by quadratic differentials}
\label{subsection;24.7.6.1}

\subsubsection{Holomorphic local coordinates}
\label{subsection;23.9.30.3}

Let $X$ be a compact Riemann surface.
Let $\phi$ be a non-zero quadratic differential on $X$
such that each zero of $\phi$ is simple.
Let $D_{\phi}$ denote the set of zeroes of $\phi$.
Let $g_{\phi}$ be the K\"ahler metric of $X\setminus D_{\phi}$
obtained as $g_{\phi}=|\phi|$.
Let $d_{\phi}$ denote the distance of $X$ induced by $g_{\phi}$.
For any $r>0$ and $P\in X$,
we set $X_P(r)=\{Q\in X\,|\,d_{\phi}(P,Q)<r\}$.

For each $P\in X\setminus D_{\phi}$,
there exists a coordinate neighbourhood $(X_P,z_P)$ around $P$
such that $\phi_{|X_P}=dz_P\,dz_P$
and $z_P(P)=0$.
We obtain $g_{\phi|X_P}=dz_P\,d\zbar_P$.
Such a coordinate $z_P$ is well defined up to
the multiplication of $\pm 1$.
There exists $r_0>0$ such that 
the coordinate $z_P$
induces the following holomorphic isomorphism
for any $r<r_0$:
\[
 (X_P(r),\phi)
 \simeq
 (U(r),(dz_P)^2).
\]

For each $P\in D_{\phi}$,
there exists a coordinate neighbourhood
$(X_P,z_P)$ such that
$\phi=\left(\frac{3}{2}\right)^2z_P\,dz_P\,dz_P$.
We have
$g_{\phi|X_P}=\left(\frac{3}{2}\right)^{2}|z_P|\,dz_P\,d\zbar_P$.
Such a coordinate $z_P$ is well defined
up to the multiplication of cubic roots of $1$.
We consider the quadratic differential
$\phi_0=\left(\frac{3}{2}\right)^2z\,(dz)^2$ on $\cnum$.
It induces the distance $d_{\phi_0}$ on $\cnum$.
We set
\[
 \nbigu_P(r)=\bigl\{
 z\in\cnum\,\big|\,
 d_{\phi_0}(0,z)<r
 \bigr\}
 =\bigl\{z\in\cnum\,\big|\,
 |z|<r^{2/3}
 \bigr\}.
\]
There exists $r_0>0$ such that 
the coordinate $z_P$ induces the following holomorphic isomorphism
for any $r<r_0$:
\[
 (X_P(r),\phi)
 \simeq
 \bigl(
 \nbigu_P(r),\phi_0
 \bigr).
\]

\subsubsection{Threshold}

According to \cite[\S3.2]{Dumas-Neitzke},
a saddle connection of $\phi$
is a geodesic segment on $(X,d_{\phi})$
such that
(i) the starting point and the ending point are
contained in $D_{\phi}$,
(ii) any other point is contained in $X\setminus D_{\phi}$.
Note that the starting point and the ending point
are not necessarily distinct.
Let $M(\phi)$ denote the minimum length
of the saddle connections of $\phi$,
which is called the threshold.
The following lemma is obvious.
\begin{lem}
For $P\in D_{\phi}$,
we have
$X_P(M(\phi))\cap D_{\phi}=\{P\}$.
\hfill\qed
\end{lem}

We shall prove the following proposition in
\S\ref{subsection;23.9.30.11}.
\begin{prop}
\label{prop;23.9.30.10}
Let $P\in D_{\phi}$.
\begin{itemize}
 \item 
A local holomorphic coordinate $z_P$ around $P$
as in {\rm\S\ref{subsection;23.9.30.3}}
is analytically continued to
a holomorphic coordinate on $X_P(M(\phi)/2)$,
and it induces a holomorphic isomorphism
\[
       (X_P(\kappa),
       \phi)
       \simeq
       (\nbigu_P(\kappa),\phi_0)
\]
for any $0<\kappa\leq M(\phi)/2$.
       In particular,
       $\phi_{|X_P(\kappa)}=(\frac{3}{2})^2z_P(dz_P)^2$.
 \item The inverse function
       $z_P^{-1}:\nbigu_P(M(\phi)/2)\simeq X_P(M(\phi)/2)$
       is analytically continued to
       a locally bi-holomorphic function
\[
\psi_P:\nbigu_P(M(\phi))\lrarr X.
\]
       We have
       $\psi_P^{\ast}(\phi)=\phi_0$,
       $\psi_P^{-1}(D_{\theta})=\{0\}$,
       and
       $\psi_P(\nbigu_P(\kappa))\subset X_P(\kappa)$
       for any $0<\kappa\leq M(\phi)$.
\end{itemize}
\end{prop}

The following lemma is easier
which is also proved in \S\ref{subsection;23.9.30.11}.

\begin{lem}
\label{lem;23.9.30.13}
Let $P\in X\setminus D_{\phi}$. 
Let $z_P$ be a holomorphic local coordinate around $P$
as in {\rm\S\ref{subsection;23.9.30.3}}.
Then, $z_P^{-1}$ is analytically continued to
a locally bi-holomorphic function 
\[
 \psi_P:
 U\bigl(d_{\phi}(P,D_{\phi})\bigr)
 \lrarr
 X.
\]
       We have
       $\psi_P^{\ast}(\phi)=(dz)^2$,
       $\psi_P^{-1}(D_{\phi})=\emptyset$,
       and
       $\psi_P(\nbigu_P(\kappa))\subset X_P(\kappa)$
       for any $0<\kappa\leq d_{\phi}(P,D_{\phi})$.
\end{lem}

\subsubsection{Metrics induced by holomorphic $1$-forms}

Let $Y$ be a simply connected Riemann surface isomorphic to a disc.
Let $a$ be a holomorphic $1$-form on $Y$.
We set $Z_a=\{Q\in Y\,|\,a(Q)=0\}$.

We obtain the K\"ahler metric $|a|^2$ on $Y\setminus Z_a$
whose curvature is flat.
As explained in \cite{Strebel},
it induces the distance $d_{a}$ on $Y$.
For a curve $\gamma$, the length of $\gamma$
is denoted by $L_{a}(\gamma)$.
In this subsection,
any geodesic is assumed to be parameterized by the arc length.
We assume the following.
\begin{assumption}
$(Y,d_a)$ is a complete metric space. 
\hfill\qed
\end{assumption}
Under the assumption, the following holds.
See \cite[Theorem 14.2.1, Theorem 16.1, Corollary 18.2]{Strebel}.

\begin{prop}
For any two points $P_1,P_2\in Y$,
there is a unique geodesic $\gamma_{P_1,P_2}$
starting from $P_1$ and ending at $P_2$. 
We have $L_{a}(\gamma_{P_1,P_2})=d_{a}(P_1,P_2)$.
If the length of
a curve connecting $P_1$ and $P_2$ is $d_{a}(P_1,P_2)$,
the curve equals to the geodesic $\gamma_{P_1,P_2}$
after a re-parametrization of the parameter.
\hfill\qed
\end{prop}

We fix a base point $P\in Y$.
We set
$M_{a,P}:=\inf_{Q\in Z_a\setminus\{P\}}d_{a}(P,Q)$.
Let $Z_{a,P}$ denote the set of $Q\in Z_a\setminus\{P\}$
such that $\gamma_{P,Q}(s)\not\in Z_a$
for any $0<s<d_a(P,Q)$.
It is easy to see
$M_{a,P}=\inf_{Q\in Z_{a,P}}d_{a}(P,Q)$.

Let $n\in\seisuu_{>0}$ denote the order of zero of $a$ at $P$.
Let $f_{a,P}$ be the holomorphic function on $Y$
such that $df_{a,P}=a$ and $f_{a,P}(P)=0$.
There exists a holomorphic function
$f_{a,P}^{\frac{1}{n+1}}$ defined on a neighbourhood $\nbigu$ of $P$
such that
$(f_{a,P}^{\frac{1}{n+1}})^{n+1}=f_{a,P}$
and that $f_{a,P}^{\frac{1}{n+1}}$ has simple zero at $P$.

For any $r>0$,
let $Y_P(r):=\{Q\in Y\,|\,d_{a}(P,Q)<r\}$
and $Y_P(r)^{\ast}=Y_P(r)\setminus\{P\}$.
We also set 
$U(r):=\{w\in\cnum\,|\,|w|<r\}$
and $U(r)^{\ast}=U(r)\setminus\{0\}$.

\begin{prop}
\label{prop;23.9.30.2}
Let $r\leq M_{a,P}$.
\begin{itemize}
 \item $f_{a,P}$ induces a proper map
       $Y_P(r)\to U(r)$,
       and the induced map
       $Y_P(r)^{\ast}\to U(r)^{\ast}$
       is a covering map of degree $n+1$.
 \item $f_{a,P}^{\frac{1}{n+1}}$
       is analytically continued
       to a holomorphic function on $Y_P(r)$,
       and it induces a bi-holomorphic map
       $Y_P(r)\simeq
       U(r^{\frac{1}{n+1}})$.
       We have
       $(f_{a,P}^{\frac{1}{n+1}})^{n+1}=f_{a,P}$
       on $Y_P(r)$.
\end{itemize} 
\end{prop}
\pf
We shall give some lemmas.

\begin{lem}
\label{lem;23.9.30.1}
For any $Q\in Y_P(M_{a,P})^{\ast}$,
$\gamma_{P,Q}(s)$ $(0<s\leq d_a(P,Q))$
are not contained in $Z_a$.
In particular, 
$Y_P(M_{a,P})^{\ast}\cap Z_a=\emptyset$.
\end{lem}
\pf
Let $S$ denote the set of $0<s\leq d_{a}(P,Q)$
such that $\gamma_{P,Q}(s)\in Z_a$.
Assume that $S$ is not empty.
Let $s_0=\min S$.
Then, $\gamma_{P,Q}(s_0)\in Z_{a,P}$,
and $d_a(P,\gamma_{P,Q}(s_0))\leq
d_a(P,Q)<M_{a,P}$,
which contradicts the definition of $M_{a,P}$.
\hfill\qed

\vspace{.1in}
For $Q\in Y_P(M_{a,P})^{\ast}$,
by Lemma \ref{lem;23.9.30.1},
we obtain
\[
f_{a,P}\circ\gamma_{P,Q}(s)=s\cdot |f_{a,P}(Q)|^{-1}f_{a,P}(Q),
\]
and $d_{a}(P,Q)=|f_{a,P}(Q)|$.
In particular, we obtain
$f_{a,P}(Y_P(r))\subset U(r)$
for any $0<r<M_{a,P}$.

\vspace{.1in}

Let $w\in\cnum\setminus\{0\}$.
There exist
$(n+1)$-distinct points 
$w_{n+1}^{(i)}\in\cnum$ $(i=1,\ldots,n+1)$
such that $(w_{n+1}^{(i)})^{n+1}=w$.
There exist $\epsilon>0$ and
$C^{\infty}$-maps
$\eta^{(i)}_{w,\epsilon}:[0,\epsilon]\to
 \nbigu\subset Y$ $(i=1,\ldots,n+1)$
such that
$\eta^{(i)}_{w,\epsilon}(0)=P$
and $f_{a,P}^{\frac{1}{n+1}}\circ\eta^{(i)}_{w,\epsilon}(s)
=s^{\frac{1}{n+1}}w_{n+1}^{(i)}$
for $0\leq s\leq \epsilon$.

\begin{lem}
\label{lem;23.9.29.1}
Suppose $|w|<M_{a,P}$.
\begin{itemize}
 \item 
There exist $C^{\infty}$-maps
$\eta^{(i)}_w:[0,1]\to Y$
such that
(i) $\eta^{(i)}_w(s)=\eta^{(i)}_{w,\epsilon}(s)$ $(0\leq s\leq\epsilon)$,
(ii) $f_{a,P}\circ\eta^{(i)}_w(s)=sw$ $(0\leq s\leq 1)$.
 \item For $i\neq j$,
       we have
       $\eta^{(i)}_w(s)\neq \eta^{(j)}_w(s)$
       for any $0<s\leq 1$.
\end{itemize}
\end{lem}
\pf
Fix $i\in\{1,\ldots,n+1\}$.
Let $I$ be the set of $\epsilon\leq \nu\leq 1$
such that 
there exists
$\eta^{(i)}_{w,\nu}:[0,\nu]\to Y$
such that
(i) $\eta^{(i)}_{w,\nu}(s)
=\eta^{(i)}_{w,\epsilon}(s)$ $(0\leq s\leq\epsilon)$,
(ii) $f_{a,P}\circ\eta^{(i)}_{w,\nu}(s)=sw$ $(0\leq s\leq \nu)$.
Note that $\epsilon\in I$, and hence $I\neq\emptyset$.
We set $\nu_0=\sup I$.
Suppose $\nu_0<1$.
For any $\epsilon\leq \nu<\nu_0$,
we have $\eta^{(i)}_{w,\nu}:[0,\nu]\to Y$
satisfying the above conditions.
Note that
\[
d_{a}(P,\eta^{(i)}_{w,\nu}(s))
\leq
\bigl|\eta^{(i)}_{w,\nu}(s)\bigr|
=s|w|<M_{a,P}.
\]
Because $\eta^{(i)}_{w,\nu}(s)\not\in Z_a$
for any $\epsilon\leq \nu<\nu_0$ and $0<s\leq \nu$,
we have $\eta^{(i)}_{w,\nu_1}(s)=\eta^{(i)}_{w,\nu_2}(s)$
for any $\epsilon\leq \nu_1\leq \nu_2$
and $0\leq s\leq \nu_1$.
Let us consider the sequence $\eta^{(i)}_{w,\nu}(\nu)$
$(\nu\in I)$.
We note that
$d_{a}(\eta^{(i)}_{w,\nu_1}(\nu_1),\eta^{(i)}_{w,\nu_2}(\nu_2))
=|\nu_1-\nu_2|\cdot|w|$.
Because $d_{a}$ is complete,
there exists the limit
$z_{\infty}=\lim_{\nu\to\nu_0} \eta^{(i)}_{w,\nu}(\nu)\in Y$.
By the continuity of $f_a$,
we have
$f_a(z_{\infty})=\nu_0 w$.
Because $f_a$ is bi-holomorphic around $z_{\infty}$ and $\nu_0w$,
we can construct
$\eta^{(i)}_{w,\nu_0+\delta}$ satisfying the above conditions
for some $\delta>0$.
It contradicts the definition of $\nu_0$.
Thus, we obtain the first claim.

Fix a pair $1\leq i,j\leq n+1$ with $i\neq j$.
Let $J$ be the set of
$0<s\leq 1$
such that
$\eta^{(i)}_{w}(s)=\eta^{(j)}_w(s)$.
Because the condition is closed,
$J$ is a closed subset of $\openclosed{0}{1}$.
Because $f_a$ is locally bi-holomorphic,
$J$ is an open subset of $\openclosed{0}{1}$.
Because $J\subset\{\epsilon\leq s\leq 1\}$,
we obtain $J\neq\openclosed{0}{1}$,
and hence $J=\emptyset$.
Thus, we obtain the second claim.
\hfill\qed

\begin{lem}
\label{lem;23.9.29.10}
For any $w\in U(r)^{\ast}$,
we have
\[
 \bigl\{
 Q\in Y_P(M_{a,P})\setminus\{P\}\,\big|\,
 f_{a,P}(Q)=w
 \bigr\}
 =\bigl\{
 \eta^{(i)}_{w}(1)\,\big|\,
 i=1,\ldots,n+1
 \bigr\}.
\]
\end{lem}
\pf
The right hand side is contained in the left hand side.
Let $Q\in Y_P(M_{a,P})$ such that $f_{a,P}(Q)=P$.
Because 
$f_{a,P}\circ \gamma_{P,Q}(|w|s)=s\cdot w$
for $0\leq s\leq 1$,
$\gamma_{P,Q}(|w|\cdot s)$
equals one of $\eta^{(i)}_{w}$.
Hence, the left hand side is contained in the right hand side.
\hfill\qed

\vspace{.1in}
Let $0<r\leq M_{a,P}$.
Let us consider the induced map
$f_{a|Y_P(r)^{\ast}}:Y_P(r)^{\ast}\to U(r)^{\ast}$.
It is locally bi-holomorphic.
By Lemma \ref{lem;23.9.29.10},
it is a covering map.
Then, $f_a^{\frac{1}{n+1}}$
is analytically continued on $Y_P(r)$.
The other claims of Proposition \ref{prop;23.9.30.2}
also follow.
\hfill\qed

\subsubsection{Proof of Proposition \ref{prop;23.9.30.10}
and Lemma \ref{lem;23.9.30.13}}

\label{subsection;23.9.30.11}

Let $\Sigma\subset T^{\ast}X$ be the spectral curve of $\phi$,
which is also a compact Riemann surface.
Let $\pi:\Sigma\to X$ denote the projection.
There exists a holomorphic $1$-form $a$ on $\Sigma$
such that $\pi^{\ast}(\phi)=a^2$.

Let $\varphi:\Sigmatilde\to \Sigma$ be a universal covering.
We obtain $\atilde=\varphi^{\ast}(a)$.
The Riemann surface $\Sigmatilde$ is bi-holomorphic to a disc.
According to \cite[Lemma 18.2]{Strebel},
the metric space
$(\Sigmatilde,d_{\atilde})$ is complete.

Let $P\in D_{\phi}$.
There exists a unique point
$\Ptilde_1\in \Sigma$ such that
$\pi(\Ptilde_1)=P$.
Let $\Ptilde_2\in \Sigmatilde$
be a point such that
$\varphi(\Ptilde_2)=\Ptilde_1$.
By definitions, we obtain
\[
M(\phi)\leq M_{\atilde,\Ptilde_2}.
\]

We obtain
$f_{\atilde,\Ptilde_2}:\Sigmatilde\to\cnum$
such that $df_{\atilde,\Ptilde_2}=\atilde$
and $f_{\atilde,\Ptilde_2}(\Ptilde_2)=0$.
As explained in Proposition \ref{prop;23.9.30.2},
there exists a holomorphic function
$f_{\atilde,\Ptilde_2}^{\frac{1}{3}}$
on $\Sigmatilde_{\Ptilde_2}(M_{\atilde,\Ptilde_2})$
such that
$(f_{\atilde,\Ptilde_2}^{\frac{1}{3}})^3=f_{\atilde,\Ptilde_2}$,
and $f_{\atilde,\Ptilde_2}^{\frac{1}{3}}$ induces
a bi-holomorphic isomorphism
\[
 \Sigmatilde_{\Ptilde_2}(r)
 \simeq
 \{|\zeta|<r^{1/3}\}
\]
for any $r<M_{\atilde,\Ptilde_2}$.
We have $\atilde=d\zeta^3$.

Let $d_{a}$ be the distance on $\Sigma$ induced by $|a|^2$.
For $r>0$,
let $\Sigma_{\Ptilde_1}(r)
=\{Q\in \Sigma\,|\,d_{a}(\Ptilde_1,Q)<r\}$
and
$\Sigma_{\Ptilde_1}(r)^{\ast}
=\Sigma_{\Ptilde_1}(r)\setminus\{\Ptilde_1\}$.

\begin{lem}
\label{lem;23.9.29.20}
Let $r\leq M_{\atilde,\Ptilde_2}/2$.
\begin{itemize}
 \item $\varphi$ induces an isomorphism
       $\Sigmatilde_{\Ptilde_2}(r)
       \simeq \Sigma_{\Ptilde_1}(r)$.
 \item $\pi$ induces a ramified covering map
       $\Sigma_{\Ptilde_1}(r)\to X_P(r)$
       of degree $2$.
       It induces
       a covering map
       $\Sigma_{\Ptilde_1}(r)^{\ast}
       \to X_P(r)^{\ast}$.
\end{itemize}
\end{lem}
\pf
We have
$\varphi(\Sigmatilde_{\Ptilde_2}(r))
\subset\Sigma_{\Ptilde_1}(r)$ for any $r>0$.
Suppose $r\leq M_{\atilde,\Ptilde_2}/2$.
Let $Q_1,Q_2\in \Sigmatilde_{\Ptilde_2}(r)$
such that $\varphi(Q_1)=\varphi(Q_2)$ and $Q_1\neq Q_2$.
Then, there exists
$\Ptilde'_2\in\varphi^{-1}(\Ptilde_1)\setminus\{\Ptilde_2\}$
such that
$d_{\atilde}(\Ptilde_2',Q_1)<r$.
We obtain
$d_{\atilde}(\Ptilde_2,\Ptilde_2')<2r
\leq M_{\atilde,\Ptilde_2}$,
which contradicts the condition of $M_{\atilde,\Ptilde_2}$.
Hence,
the restriction of $\varphi$ induces
a bi-holomorphic isomorphism
between $\Sigmatilde_{\Ptilde_2}(r)$
and its image
$\varphi(\Sigmatilde_{\Ptilde_2}(r))$.

Let $Q\in \Sigma_{\Ptilde_1}(r)$.
There exists a curve $\gamma^Q_0$
starting from $\Ptilde_1$
and ending at $Q$
whose length is strictly smaller than $r$.
According to \cite[Theorem 18.2.2]{Strebel},
there exists a unique geodesic $\gamma^Q$
from $\Ptilde_1$ to $Q$
such that $\gamma^Q$ is homotopic to $\gamma^Q_0$.
The length $L$ of $\gamma^Q$ is strictly less than $r$.
Let $\gammatilde^Q$ be the lift of $\gamma^Q$
to $\Sigmatilde$
starting from $\Ptilde_2$.
Then, $\gammatilde^Q(L)$ is contained in
$\Sigmatilde_{\Ptilde_2}(r)$.
We obtain
$\Sigma_{\Ptilde_1}(r)
=\varphi(\Sigmatilde_{\Ptilde_2}(r))$.
Thus, we obtain the first claim.

We have
$\pi(\Sigma_{\Ptilde_1}(r))\subset X_P(r)$
for any $r>0$.
Suppose that $r<M_{\atilde,\Ptilde_2}/2$.
The induced map
$\Sigma_{\Ptilde_1}(r)^{\ast} \to X_P(r)^{\ast}$
is locally bi-holomorphic.
Let $Q\in X_P(r)$.
There exists a curve $\gamma^Q$
starting from $P$ and ending at $Q$
whose length is strictly smaller than $r$.
According to \cite[Theorem 18.2.2]{Strebel},
there exists a unique geodesic $\gamma^Q$
from $P$ to $Q$
such that $\gamma^Q$ is homotopic to $\gamma^Q_0$.
The length $L$ of $\gamma^Q$ is strictly less than $r$.
Because $L<M_{\atilde,\Ptilde_2}/2$,
$\gamma^Q(s)$ is not contained in $D_{\phi}$
for any $0<s\leq L$.
Hence, there exists exactly two lifts
$\gamma^{Q(i)}$ $(i=1,2)$ of $\gamma^Q$
to $\Sigma_{\Ptilde_1}(r)$.
Hence,
$\Sigma_{\Ptilde_1}(r)^{\ast}
\to X_P(r)^{\ast}$
is a covering map of degree $2$.
Thus, we obtain the second claim.
\hfill\qed

\vspace{.1in}
Let $\sigma$ be the involution of $\Sigma$
induced by the multiplication of $-1$ on $T^{\ast}X$.
It induces an involution $\sigma$
on $\Sigma_{\Ptilde_1}(r)$ for any $r$.
By Lemma \ref{lem;23.9.29.20},
we obtain the involution $\sigma$
on $\Sigmatilde_{\Ptilde_2}(r)$
for $r<M_{\atilde,\Ptilde_2}/2$.
\begin{lem}
The induced involution $\sigma$
on $\Sigmatilde_{\Ptilde_2}(M_{\atilde,\Ptilde_2})$
is described as $\sigma(\zeta)=-\zeta$.
\end{lem}
\pf
Because
$\atilde=d(\zeta^3)$,
$\sigma^{\ast}\atilde=-\atilde$,
$\sigma(\Ptilde_2)=\Ptilde_2$
and $\sigma^2=\id$,
we obtain $\sigma^{\ast}(\zeta)=-\zeta$.
\hfill\qed

\vspace{.1in}

There exists a bi-holomorphic map
\[
 \Phi_1:
 U\bigl(
 (M_{\atilde,\Ptilde_2}/2)^{2/3}
 \bigr)
 \lrarr
 X_P(M_{\atilde,\Ptilde_2}/2)
\]
such that
$\pi\circ\varphi
=\Phi_1\circ f^{\frac{1}{3}}_{\atilde,\Ptilde_2}$
on $\Sigmatilde_{\Ptilde_2}(M_{\atilde,\Ptilde_2}/2)$.
For each $0<r<M_{\atilde,\Ptilde_2}/2$,
$\Phi_1$ induces an isomorphism
$U(r^{2/3})\simeq X_P(r)$.
Because $M(\phi)\leq M_{\atilde,\Ptilde_2}$,
this implies the first claim of Proposition \ref{prop;23.9.30.10}.
We obtain the second claim of Proposition \ref{prop;23.9.30.10}
from the following lemma.

 \begin{lem}
\label{lem;23.9.30.12}
  $\Phi_1$ is analytically continued to
a locally bi-holomorphic map
\[
 \Phi_2:
 U\bigl(
 M_{\atilde,\Ptilde_2}^{2/3}
 \bigr)
 \lrarr
 X_P(M_{\atilde,\Ptilde_2})
\]
such that 
$\pi\circ\varphi
=\Phi_2\circ f^{\frac{1}{3}}_{\atilde,\Ptilde_2}$
on $\Sigmatilde_{\Ptilde_2}(M_{\atilde,\Ptilde_2})$.
 \end{lem}
\pf
For $\theta_1<\theta_2$ such that $\theta_2-\theta_1<2\pi$,
we consider the following regions:
\[
 S_0(\theta_1,\theta_2)
 =\bigl\{
 \xi\in\cnum\,\big|\,
 (M_{\atilde,\Ptilde_2}/3)^{2/3}<|\xi|<
 (M_{\atilde,\Ptilde_2}/2)^{2/3},\,\,
 \theta_1<\arg(\xi)<\theta_2
 \bigr\},
\]
\[
 S_1(\theta_1,\theta_2)
 =\bigl\{
 \xi\in\cnum\,\big|\,
 (M_{\atilde,\Ptilde_2}/3)^{2/3}<|\xi|
 <M_{\atilde,\Ptilde_2}^{2/3},\,\,
 \theta_1<\arg(\xi)<\theta_2
 \bigr\}.
\]
Let $\rho:\cnum\to \cnum$ be the map defined by
$\rho(\zeta)=\zeta^2$.
For $i=0,1$,
we have the decomposition
into the connected components
\[
 \rho^{-1}(S_i(\theta_1,\theta_2))
=\rho^{-1}(S_i(\theta_1,\theta_2))_1
 \sqcup
 \rho^{-1}(S_i(\theta_1,\theta_2))_2.
\]
We obtain
\[
 (\rho^{-1})_j:
 S_1(\theta_1,\theta_2)
 \simeq
 \rho^{-1}(S_1(\theta_1,\theta_2))_j
\]
obtained as the inverse of $\rho$.
We have the equalities 
\[
\pi\circ\varphi\circ(f_{\atilde,\Ptilde_2}^{\frac{1}{3}})^{-1}\circ
 (\rho^{-1})_1
=\pi\circ\varphi\circ(f_{\atilde,\Ptilde_2}^{\frac{1}{3}})^{-1}\circ
 (\rho^{-1})_2
\]
on $S_1(\theta_1,\theta_2)$,
because their restriction to 
$S_0(\theta_1,\theta_2)$ are equal.
Then, we obtain the claim of Lemma \ref{lem;23.9.30.12},
which implies the second claim of Proposition \ref{prop;23.9.30.10}.
\hfill\qed

\vspace{.1in}
Let $P\in X\setminus D_{\phi}$.
Let $\Ptilde\in (\pi\circ\varphi)^{-1}(P)$.
By definitions,
we have
\[
 M_{\atilde,\Ptilde}
 =d_{\phi}(P,D_{\phi}).
\]
We obtain the locally bi-holomorphic map
$U(M_{\atilde,\Ptilde})\to X$
as the composition of
\[
\begin{CD}
 U(M_{\atilde,\Ptilde})
 @>{(f_{\atilde,\Ptilde})^{-1}}>{\simeq}>
 \Sigmatilde_{\Ptilde}(M_{\atilde,\Ptilde})
 @>{\pi\circ\varphi}>>
 X.
\end{CD}
\]
This gives a desired analytic continuation of $z_P^{-1}$.
Thus, we obtain Lemma \ref{lem;23.9.30.13}.
\hfill\qed

\subsection{Higgs bundles of rank $2$ with smooth spectral curve}

Let $X$ be a compact Riemann surface.
Let $(E,\theta)$ be a Higgs bundle
of degree $0$ and rank $2$ on $X$.
We assume that $\det(E)=\nbigo_X$
and $\tr\theta=0$.
We assume that the spectral curve $\Sigma_{\theta}$ is smooth, reduced
and irreducible.
Let $D_{\theta}\subset X$ denote
the set of the critical values of
the projection $\pi:\Sigma_{\theta}\to X$.
At each point of $D_{\theta}$,
$\phi:=-\det(\theta)$ has a simple zero.
For any $W\subset X$,
we set $\Sigma_{\theta|W}=\Sigma_{\theta}\times_XW$.

We obtain the flat K\"ahler metric
$g_{\phi}=|\phi|$ of $X\setminus D_{\theta}$.
Let $d_{\phi}$ denote the distance of $X$ induced by $g_{\phi}$.
For any $P\in X$ and $r>0$,
we set $X_P(r)=\{Q\in X\,|\,d_{\phi}(P,Q)<r\}$.
We set
\[
 \kappa_0=M(\phi)/2.
\]

\subsubsection{Decomposition into the diagonal part
and the off-diagonal part}
\label{subsection;23.9.28.31}

Let $P\in X\setminus D_{\theta}$.
There exists a holomorphic coordinate neighbourhood
$(X_P,z_P)$ such that $\phi=(dz_P)^2$
and $z_P(P)=0$.
We set $E_P=E_{|X_P}$.
Let $f_P$ be the endomorphism of 
$E_P$
determined by
$\theta_{|X_P}=f_P\,dz_P$.

We have the decomposition
\[
 (E,\theta)_{|X_P}
 =(E_{P,1},dz_P)\oplus (E_{P,2},-dz_P).
\]
We set
\[
 \End(E)^{\circ}_P
 =\bigoplus\End(E_{P,i}),
 \quad
 \quad
 \End(E)^{\bot}_P
 =\Hom(E_{P,1},E_{P,2})
 \oplus
 \Hom(E_{P,2},E_{P,1}).
\]
The subbundles
$\End(E)^{\circ}_P$
and $\End(E^{\bot})_P$
of $\End(E)_{|X_P}$
are well defined.
By varying $P\in X\setminus D_{\theta}$,
we obtain the subbundles
$\End(E)^{\circ}$
and $\End(E)^{\bot}$
of $\End(E)_{|X\setminus D_{\theta}}$,
and the decomposition
\begin{equation}
\label{eq;23.9.28.30}
 \End(E)_{|X\setminus D_{\theta}}
 =\End(E)^{\circ}
 \oplus
 \End(E)^{\bot}.
\end{equation}

\subsection{Approximation of large solutions of Hitchin equation}

\subsubsection{Limiting configuration}
\label{subsection;23.10.19.40}

Let us recall the notion of the limiting configuration
in \cite{MSWW}.
There exists a line bundle $L$ on $\Sigma_{\theta}$
with an isomorphism
$\pi_{\ast}\Sigma\simeq (E,\theta)$,
where the Higgs field of $\pi_{\ast}L$
is induced by the tautological $1$-form of $T^{\ast}X$.
There exists a singular flat Hermitian metric $h_L$ of $L$
with the following property.
\begin{itemize}
 \item Let $Q\in\Sigma_{\theta}$
       be any point of $\pi^{-1}(D_{\theta})$.
       Then, there exists a holomorphic local
       coordinate $\zeta_Q$ of $\Sigma_{\theta}$
       and a frame $v_Q$ of $L$
       around $Q$
       such that
       $h_L(v_Q,v_Q)=|\zeta_Q|$.
\end{itemize}
We obtain the Hermitian metric
$h_{\infty}=\pi_{\ast}h_L$ of $E_{|X\setminus D_{\theta}}$.
Because $\det(h_{\infty})$ induces a $C^{\infty}$ Hermitian metric
of $\det(E)$,
we may adjust $h_{\infty}$ so that $\det(h_{\infty})=1$.

Let $\nabla_{\infty}$ denote the flat unitary connection of
$(E_{|X\setminus D_{\theta}},h_{\infty})$
obtained as the Chern connection of $h_{\infty}$.
We consider the flat connection of
$T^{\ast}(X\setminus D_{\theta})$
obtained as the Levi-Civita connection of $g_{\phi}$.
For any $m\in\seisuu_{\geq 0}$,
the induced connections on
$\End(E)\otimes(T^{\ast}(X\setminus D_{\theta}))^{\otimes m}$
are also denoted by $\nabla_{\infty}$.
We obtain
$\nabla_{\infty}^j:
\End(E)\otimes (T^{\ast}X)^{\otimes k}
\to
\End(E)\otimes (T^{\ast}X)^{\otimes (k+j)}$
for any $k,j\geq 0$.

\subsubsection{Model Higgs bundle}

We set $E_0=\nbigo_{\cnum}u_1\oplus \nbigo_{\cnum}u_2$.
We consider the Higgs field $\theta_0$ of $E_0$
defined by
\[
 \theta_0(u_1)=\frac{3}{2}u_2\,dz,
 \quad
 \theta_0(u_2)=\frac{3}{2}u_1\,z\,dz.
\]
Let $C_0$ be the non-degenerate symmetric pairing
such that
$C_0(u_i,u_i)=0$ $(i=1,2)$ and $C_0(u_1,u_2)=1$.
There exists a unique harmonic metric $h_0$
of $(E_0,\theta_0)$ compatible with $C_0$
(see \cite{Li-Mochizuki3}).
Let $\rho_t:\cnum\to\cnum$ be given by
$\rho_t(z)=t^{2/3}z$.
By setting
\[
 \rho_t^{\ast}(u_1)=t^{1/6}u_1,
 \quad
 \rho_t^{\ast}(u_2)=t^{-1/6}u_2,
\]
we obtain the isomorphism
$\rho_t^{\ast}(E_0,\theta_0,C_0)\simeq
(E_0,t\theta_0,C_0)$.
We obtain harmonic metrics $h_{0,t}=\rho_t^{\ast}h_0$
of $(E_0,t\theta_0)$
which are compatible with $C_0$.

\subsubsection{Harmonic metrics around $P\in D_{\theta}$}
\label{subsection;23.10.16.5}

Let $P\in D_{\theta}$.
We obtain the map
$\psi_{P,2\kappa_0}:
\nbigu_P(2\kappa_0)\to X$
as in \S\ref{subsection;24.7.6.1}.
We obtain the Higgs bundle
\[
(E_P,\theta_P)
=\psi_{P,2\kappa_0}^{\ast}(E,\theta)
\]
on $\nbigu_P(2\kappa_0)$.
There exists the decomposition
\begin{equation}
\label{eq;23.9.28.32}
 \End(E_P)_{|\nbigu_P(2\kappa_0)\setminus\{0\}}
=\End(E_P)^{\circ}
\oplus
 \End(E_P)^{\bot}
\end{equation}
as in \S\ref{subsection;23.9.28.31}.

We set
$h_{P,\infty}:=
\psi_{P,2\kappa_0}^{\ast}(h_{\infty})$
on $\nbigu_P(2\kappa_0)\setminus\{0\}$.
There exists a non-degenerate symmetric pairing $C_P$
of $(E_P,\theta_P)$
such that
$C_P$ is compatible with $h_{P,\infty}$
on $\nbigu_P(2\kappa_0)\setminus\{0\}$
(see \cite{Mochizuki-Szabo}).
For any $m\in\seisuu_{\geq 0}$,
we obtain the connections
$\nabla_{P,\infty}$ of
$\End(E_P)\otimes
T^{\ast}(\nbigu_P(2\kappa_0)\setminus\{0\})^{\otimes m}$
induced by the Chern connection of $h_{P,\infty}$
and the Levi-Civita connection of $|\phi_0|$.

There exists a natural isomorphism
\begin{equation}
\label{eq;23.10.7.100}
 (E_P,\theta_P,C_P)
 \simeq
 (E_0,\theta_0,C_0)_{|\nbigu_P(2\kappa_0)}.
\end{equation}
For any $t>0$,
let $h_{P,t}$ be the harmonic metric of
$(E_P,t\theta_P)$ obtained as the restriction of $h_{0,t}$
under the isomorphism (\ref{eq;23.10.7.100}).
It is compatible with $C_P$.
We define the endomorphism
$s_{P,t}$ of
$E_{P|\nbigu_P(2\kappa_0)\setminus\{P\}}$
by $h_{P,t}=h_{P,\infty}\cdot s_{P,t}$.
We obtain the decomposition
\[
 s_{P,t}
 =s_{P,t}^{\circ}
 +s_{P,t}^{\bot}
\]
on $\nbigu_{P}(2\kappa_0)\setminus\{0\}$
corresponding to the decomposition (\ref{eq;23.9.28.32}).

Let $0<\kappa<\kappa_0$.
We set $\delta(\kappa):=\min\{(\kappa_0-\kappa)/10,\kappa/20\}$.
We set
\[
 \nbigu^{(1)}_{P,\kappa}
 =\nbigu_P(\kappa_0-\delta(\kappa)),
 \quad
 \nbigu^{(2)}_{P,\kappa}
 =\nbigu_P(\kappa_0-2\delta(\kappa)).
\]

\begin{lem}
\label{lem;23.9.28.50}
For any $j\in\seisuu_{\geq 0}$,
there exists $C_j>0$
such that the following holds on
$\nbigu^{(1)}_{P,\kappa}\setminus
 \nbigu^{(2)}_{P,\kappa}$:
\begin{equation}
\label{eq;23.9.30.20}
 |\nabla_{\infty}^j(s_{P,t}^{\circ}-\id)|_{h_{\infty},g_{\phi}}
 \leq C_j
 e^{-8(\kappa+7\delta(\kappa)) t},
 \quad\quad\quad
 |\nabla_{\infty}^j(s_{P,t}^{\bot})|_{h_{\infty},g_{\phi}}
 \leq C_j
 e^{-4(\kappa+7\delta(\kappa)) t}.  
\end{equation}
\end{lem}
\pf
Let $Q\in \nbigu^{(1)}_{P,\kappa}\setminus \nbigu^{(2)}_{P,\kappa}$.
There exists a holomorphic embedding
$\psi_Q:
U(\kappa_0-2\delta(\kappa))\to
\nbigu^{(0)}_{P,\kappa}$
such that
$\psi_Q(0)=Q$,
$\psi_Q^{\ast}(\phi_0)=(dz)^2$,
and $0\not\in\Image(\psi_Q)$.
Take $\kappa+7\delta(\kappa)<\kappa'<
\kappa+8\delta(\kappa)=\kappa_0-2\delta(\kappa)$.

For any $t>0$,
let $\rho_t:\cnum\to \cnum$ be the map
defined by $\rho_t(z)=t^{-1}z$.
We put $\psi_{Q,t}=\psi_Q\circ\rho_t$.
Let $g_0$ be the standard Euclidean metric on $\cnum$.
We apply Corollary \ref{cor;23.9.28.40}
to $\psi_{Q,t}^{\ast}(E_P,t\theta_P)$
on $U(t(\kappa_0-2\delta(\kappa)))$.
For any $j\geq 0$,
there exists $C_j'>0$ such that
the following holds around
$0\in U\bigl(t(\kappa_0-2\delta(\kappa))\bigr)$:
\[
  \Bigl|
  \psi_{Q,t}^{\ast}\bigl(\nabla_{\infty}^j(s_{P,t}^{\circ}-\id)\bigr)
  \Bigr|_{\psi_{Q,t}^{\ast}h_{\infty},g_{0}}
 \leq C'_j
 e^{-8\kappa' t},
 \quad\quad\quad
 \Bigl|
 \psi_{Q,t}^{\ast}\bigl(
 \nabla_{\infty}^j(s_{P,t}^{\bot})
 \bigr)
 \Bigr|_{\psi_{Q,t}^{\ast}h_{\infty},g_{0}}
 \leq C_j'
 e^{-4\kappa' t}.
\]
Hence, the following holds around
$0\in U\bigl(\kappa_0-2\delta(\kappa)\bigr)$:
\[
   \Bigl|
  \psi_{Q}^{\ast}\bigl(\nabla_{\infty}^j(s_P^{\circ}-\id)\bigr)
  \Bigr|_{\psi_{Q}^{\ast}h_{\infty},g_{0}}
 \leq C'_j
 t^je^{-8\kappa' t},
 \quad\quad\quad
 \Bigl|
 \psi_{Q}^{\ast}\bigl(
 \nabla_{\infty}^j(s_P^{\bot})
 \bigr)
 \Bigr|_{\psi_Q^{\ast}h_{\infty},g_0}
 \leq C_j'
 t^je^{-4\kappa' t}.
\]
Then,
we obtain (\ref{eq;23.9.30.20})
around $Q$.
\hfill\qed

\vspace{.1in}

For any $\kappa'<\kappa_0$,
we identify
$\nbigu_P(\kappa')$ and $X_P(\kappa')$.
We also identify
$(E_P,\theta_P)_{|\nbigu_P(\kappa)}$
and $(E,\theta)_{|X_P(\kappa)}$.
We obtain
$X^{(i)}_{P,\kappa}\subset X_{P}(\kappa_0)$ $(i=1,2)$
corresponding to
$\nbigu^{(i)}_{P,\kappa}$.

\subsubsection{Approximation of harmonic metrics}

For $0<\kappa<\kappa_0$,
by patching
$h_{\infty}$ and $h_{P,t}$ $(P\in D_{\theta})$,
we construct Hermitian metrics $\htilde_{\kappa,t}$
such that the following holds.
\begin{itemize}
 \item $\htilde_{\kappa,t}=h_{P,t}$
       on $X^{(2)}_{P,\kappa}$.
 \item $\htilde_{\kappa,t}=h_{\infty}$
       on $X\setminus \bigcup X^{(1)}_{P,\kappa}$.
 \item Let $s(h_{\infty},\htilde_{\kappa,t})$
       be the automorphism of
       $E_{|X\setminus D_{\theta}}$
       determined by
       $\htilde_{\kappa,t}=h_{\infty}\cdot
       s(h_{\infty},\htilde_{\kappa,t})$.
       For any $j\in\seisuu_{\geq 0}$,
       there exists $C_j>0$ such that
       the following holds on
       $\bigcup(X_{P,\kappa}^{(1)}\setminus X^{(2)}_{P,\kappa})$:
\[
       \Bigl|
       \nabla_{\infty}^j
       \bigl(s(h_{\infty},\htilde_{\kappa,t})^{\circ}-\id\bigr)
       \Bigr|_{h_{\infty},g_{\phi}}
       \leq
       C_je^{-8(\kappa+7\delta(\kappa)) t},
\]
\[
        \bigl|
       \nabla_{\infty}^j
       (s(h_{\infty},\htilde_{\kappa,t})^{\bot})
       \bigr|_{h_{\infty},g_{\phi}}
       \leq
       C_je^{-4(\kappa+7\delta(\kappa)) t}.
\]       
\item $\det(\htilde_{\kappa,t})=1$.
\end{itemize}

Let $s(\htilde_{\kappa,t},h_t)$
be the automorphism of $E$
determined by
$h_t=\htilde_{\kappa,t}\cdot s(\htilde_{\kappa,t},h_t)$.

Let $\nabla_{h_t}$ denote the Chern connection of $h_t$.
Let $g_X$ be a K\"ahler metric of $X$.
The Levi-Civita connection and $\nabla_{h_t}$
induce connections on $\End(E)\otimes(T^{\ast}X)^{\otimes m}$,
which are denoted by $\nabla_{h_t,g_X}$.
We obtain the differential operators
$\nabla_{h_t,g_X}^j:\End(E)\otimes(T^{\ast}X)^{\otimes m}
\to \End(E)\otimes(T^{\ast}X)^{\otimes m+j}$.
We shall prove the following theorem
in \S\ref{subsection;23.10.16.1}--\ref{subsection;23.10.16.2}.

\begin{thm}
\label{thm;23.10.1.40}
For any $j\geq 0$,
there exists $C_j>0$ such that 
the following holds on $X$ for any $t\geq 1$:
\[
 \Bigl|
 \nabla_{h_{t},g_X}^j
 \bigl(s(\htilde_{\kappa,t},h_t)-\id\bigr)
 \Bigr|_{h_t,g_X}
+\Bigl|
 \nabla_{h_{t},g_X}^j
 \bigl([\theta,s(\htilde_{\kappa,t},h_t)]\bigr)
 \Bigr|_{h_t,g_X}
+\Bigl|
 \nabla_{h_{t},g_X}^j
 \bigl([\theta^{\dagger}_{h_t},s(\htilde_{\kappa,t},h_t)]\bigr)
 \Bigr|_{h_t,g_X}
 \leq C_je^{-4\kappa t}.
\] 
\end{thm}

To simplify the description,
we set $s_{\kappa,t}=s(\htilde_{\kappa,t},h_t)$.
To simplify the argument,
we use the following result
already proved in \cite{MSWW}.
(See also \cite{Mochizuki-Szabo}.)

\begin{prop}[\cite{MSWW}]
\label{prop;23.9.28.101}
There exists $C>0$ and $\epsilon>0$ such that
$|s_{\kappa,t}-\id|_{\htilde_{\kappa,t}}\leq Ce^{-\epsilon t}$.
Higher derivatives are also dominated.
\hfill\qed
\end{prop}

It is our purpose to improve the rate of the convergence.

\subsubsection{Preliminary}
\label{subsection;23.10.16.1}

Let $R(\htilde_{\kappa,t})$
denote the curvature of
the Chern connection of
$E$ with $\htilde_{\kappa,t}$.
We set
\[
 F(\htilde_{\kappa,t})
 =R(\htilde_{\kappa,t})
 +t^2[\theta,\theta^{\dagger}_{\htilde_{\kappa,t}}].
\]
We obtain the following lemma by the construction of $\htilde_{\kappa,t}$.
\begin{lem}
\label{lem;23.9.28.103}
The support of $F(\htilde_{\kappa,t})$
is contained in 
$\bigcup_P\bigl(
 X^{(1)}_{P,\kappa}\setminus X^{(2)}_{P,\kappa}
 \bigr)$.
Moreover,
there exists $C>0$  such that
\[
 |F(\htilde_{\kappa,t})^{\bot}|_{\htilde_{\kappa,t},g_X}
\leq 
 C e^{-4(\kappa+6\delta(\kappa)) t},
 \quad\quad
  |F(\htilde_{\kappa,t})^{\circ}|_{\htilde_{\kappa,t},g_X}
\leq 
 C e^{-8(\kappa+6\delta(\kappa)) t}.
\]
\hfill\qed
\end{lem}

\subsubsection{Weaker estimate of the $L^2$-norm of the first derivatives}
\label{subsection;23.9.30.50}

Let $s(h_{\infty},h_t)$
be the automorphism of $E_{|X\setminus D_{\theta}}$
determined by
$h_t=h_{\infty}s(h_{\infty},h_t)$.

\begin{lem}
\label{lem;23.9.28.102}
There exists $C>0$ such that
$|s(h_{\infty},h_t)^{\bot}|\leq
 Ce^{-2\sqrt{2}(\kappa+7\delta(\kappa)) t}$
on $X\setminus\bigcup X^{(2)}_{P,\kappa}$.
\end{lem}
\pf
Let $Q\in X\setminus \bigcup X^{(2)}_{P,\kappa}$.
By Lemma \ref{lem;23.9.30.13},
there exists a locally bi-holomorphic map
$\psi_{Q,t}:U\bigl(t(\kappa_0-2\delta(\kappa))\bigr)\to
X\setminus D_{\theta}$
such that $\psi_Q(0)=Q$,
$\psi_Q^{\ast}(t^2\phi)=(dz)^2$.
By applying Proposition \ref{prop;23.9.27.20},
we obtain the claim of Lemma \ref{lem;23.9.28.102}.
\hfill\qed

\begin{cor}
There exists $C>0$ such that
$|s_{\kappa,t}^{\bot}|\leq
 Ce^{-2\sqrt{2}(\kappa+7\delta(\kappa)) t}$
on $X\setminus\bigcup X^{(2)}_{P,\kappa}$.
\hfill\qed
\end{cor}

\begin{cor}
There exists $C>0$ such that
the following holds on $X\setminus D_{\theta}$:
 \begin{equation}
\label{eq;23.9.28.100}
 \Bigl|
 \Tr\bigl(s_{\kappa,t}F(\htilde_{\kappa,t})\bigr)
 \Bigr|\leq Ce^{-(4+2\sqrt{2})(\kappa+6\delta(\kappa)) t}.
 \end{equation}
\end{cor}
\pf
Because
$\Tr\bigl(
 s_{\kappa,t}F(\htilde_{\kappa,t})
 \bigr)
=
 \Tr\bigl(
 s^{\circ}_{\kappa,t}F(\htilde_{\kappa,t})^{\circ}
 \bigr)
+
 \Tr\bigl(
 s^{\bot}_{\kappa,t}F(\htilde_{\kappa,t})^{\bot}
 \bigr)$,
we obtain (\ref{eq;23.9.28.100})
from Proposition \ref{prop;23.9.28.101}
Lemma \ref{lem;23.9.28.103}
and Lemma \ref{lem;23.9.28.102}.
\hfill\qed 

\vspace{.1in}
According to \cite[Lemma 3.1]{s1},
we have
\[
 \Delta_{g_X} \Tr(s_{\kappa,t})
 =-\Tr\bigl(
 F(\htilde_{\kappa,t})s_{\kappa,t}
 \bigr)
-\bigl|
 \delbar_E(s_t)s_t^{-1/2}
 \bigr|^2_{\htilde_{\kappa,t},g_X}
 -\bigl|
 [\theta,s_{\kappa,t}]s_{\kappa,t}^{-1/2}
 \bigr|^2_{\htilde_{\kappa,t},g_X}
\]

Because $\int\Delta_{g_X}\Tr(s_{\kappa,t})=0$,
there exists $C'>0$ such that
\[
 \|\delbar_E(s_{\kappa,t})s_{\kappa,t}^{-1/2}\|
 _{L^2,\htilde_{\kappa,t},g_X}
\leq C'e^{-(2+\sqrt{2})(\kappa+6\delta(\kappa)) t},
 \quad\quad
 \|[\theta,s_{\kappa,t}]s_{\kappa,t}^{-1/2}\|
 _{L^2,\htilde_{\kappa,t},g_X}
\leq C'e^{-(2+\sqrt{2})(\kappa+6\delta(\kappa)) t}.
\]
By Proposition \ref{prop;23.9.28.101},
there exists $C''>0$ such that
\[
  \|\delbar_E(s_{\kappa,t})\|_{L^2,\htilde_{\kappa,t},g_X}
\leq C''e^{-(2+\sqrt{2})(\kappa+6\delta(\kappa)) t},
 \quad\quad
 \|[\theta,s_{\kappa,t}]\|_{L^2,\htilde_{\kappa,t},g_X}
\leq C'' e^{-(2+\sqrt{2})(\kappa+6\delta(\kappa)) t}.
\]
Note that it is equivalent to
\begin{equation}
\label{eq;23.9.28.113}
  \|\delbar_E(s_{\kappa,t})\|_{L^2,\htilde_{\kappa,t},g_{\phi}}
\leq C''e^{-(2+\sqrt{2})(\kappa+6\delta(\kappa)) t},
 \quad\quad
 \|[\theta,s_{\kappa,t}]\|_{L^2,\htilde_{\kappa,t},g_{\phi}}
\leq C'' e^{-(2+\sqrt{2})(\kappa+6\delta(\kappa)) t}.
\end{equation}

\subsubsection{Weaker estimate of higher derivatives}

\begin{lem}
\label{lem;23.10.1.2}
For any $j>0$,
there exists $C_j>0$ such that 
the following holds on 
$X\setminus \bigcup X^{(2)}_{P,\kappa}$:
\begin{equation}
\label{eq;23.9.28.114}
 |\nabla_{\infty}^j(s_{\kappa,t})|_{h_{\infty},g_{\phi}}
 \leq
 C_je^{-2\sqrt{2}(\kappa+6\delta(\kappa)) t}.
\end{equation}
\end{lem}
\pf
By using Lemma \ref{lem;23.9.28.102},
there exists $C'>0$ such that 
\[
\Bigl|
t^2[\theta,\theta^{\dagger}_{h_t}]
\Bigr|_{h_{\infty},g_{\phi}}
\leq C't^2e^{-2\sqrt{2} (\kappa+7\delta(\kappa)) t}
\]
on $X\setminus \bigcup X^{(2)}_{P,\kappa}$.
Let $R(h_{t})$ denote the curvature of
the Chern connection of $h_t$.
Because
$R(h_t)+t^2[\theta,\theta^{\dagger}_{h_t}]=0$,
we obtain
\[
\Bigl|R(h_t)
 \Bigr|_{h_{\infty},g_{\phi}}
 \leq
 C't^2e^{-2\sqrt{2}(\kappa+7\delta(\kappa)) t}
\]
on $X\setminus \bigcup X^{(2)}_{P,\kappa}$.
By using Lemma \ref{lem;23.9.28.110}
and Lemma \ref{lem;23.9.28.111}
with the estimate (\ref{eq;23.9.28.113}),
we obtain the desired estimate (\ref{eq;23.9.28.114}).

\hfill\qed

\vspace{.1in}

We consider the connection
$\nabla_{\htilde_{\kappa,t},g_{\phi}}$
of $\End(E)\otimes T^{\ast}(X\setminus D_{\phi})^{\otimes m}$
induced by the Chern connection of $\htilde_{\kappa,t}$
and the Levi-Civita connection of $g_{\phi}$.
Let $\rho_{P,\kappa}:X_{P,\kappa}^{(2)}\to\real_{\geq 0}$
be defined by
$\rho_{P,\kappa}(Q)=d_{\phi}(Q,X_P(2\delta(\kappa)))$.

\begin{lem}
\label{lem;23.9.30.30}
Let $P\in D_{\theta}$.
For any $j>0$
and $0<\gamma<2\sqrt{2}$,
there exists $C_{j,\gamma}>0$ such that
the following holds on 
$X^{(2)}_{P,\kappa}\setminus X_P(2\delta(\kappa))$:
\[
\Bigl|
 \nabla_{\htilde_{\kappa,t},g_{\phi}}^j
 (s_{\kappa,t})
 \Bigr|_{\htilde_{\kappa,t},g_{\phi}}
 \leq
 C_{j,\gamma}
 e^{-\gamma \cdot\rho_{P,\kappa} t}.
\]
\end{lem}
\pf
Let $Q\in X^{(2)}_{P,\kappa}\setminus X_P(2\delta(\kappa))$.
There exists a locally bi-holomorphic map
\[
 \psi_{Q,t}:
 U\bigl(t \rho_P(Q)\bigr)
 \lrarr
 X
\]
such that $\psi_Q(0)=Q$
and $\psi_Q^{\ast}(t^2\phi)=(dz)^2$.
We apply Proposition \ref{prop;23.9.27.20}
to the Higgs bundle
$(E_{Q,t},\theta_{Q,t})=
\psi_{Q,t}^{\ast}(E,t\theta)$
and harmonic metrics
$h_{t,Q}=\psi_{Q,t}^{\ast}(h_t)$
and
$h_{P,t,Q}=\psi_{Q,t}^{\ast}(h_{P,t})$.
Let $\gamma<\gamma'<2\sqrt{2}$
and $t\geq 10\delta(\kappa)^{-1}$.
There exists $C>0$ such that 
the following holds on
$U(t\rho_{P,\kappa}(Q)-1)$:
\[
\Bigl|
\delbar_{E_{Q,t}}\Bigl(
 \psi_{Q,t}^{\ast}(s_{\kappa,t})^{-1}
 \del_{E_{Q,t},h_{P,t,Q}}
 \psi_{Q,t}^{\ast}(s_{\kappa,t})
\Bigr)
\Bigr|_{h_{P,t,Q},g_0}
 \leq C e^{-\gamma' (t\rho_{P,\kappa}(Q)-|z|)}.
\]
We have
$U(2)\subset U(t\rho_{P,\kappa}(Q))$.
By (\ref{eq;23.9.28.113}),
we obtain
\[
 \int_{U(1)}\bigl|
 \del_{E_{Q,t},h_{P,t,Q}}
 \psi_{Q,t}^{\ast}(s_{\kappa,t})
 \bigr|_{h_{P,t,Q},g_0}^2
 \leq C e^{-(2+\sqrt{2})(\kappa+6\delta(\kappa)) t}.
\]
Note that
$\rho_{P,\kappa}(Q)< \kappa+6\delta(\kappa)$,
and hence
$2\sqrt{2}\rho_{P,\kappa}(Q)
<(2+\sqrt{2})\bigl(
 \kappa+6\delta(\kappa)
 \bigr)$.
By Lemma \ref{lem;23.9.30.20}
and Lemma \ref{lem;23.9.28.110},
we obtain the following on $U(1)$:
\[
\bigl|
\psi_{Q,t}^{\ast}\nabla^{j}_{h_{P,t},g_{\phi}}
 (s_{\kappa,t})
 \bigr|_{\psi_{Q,t}^{\ast}(h_{P,t}),g_0}
 \leq
 Ce^{-\gamma'\rho_{P,\kappa}(Q)t}.
\]
We obtain the following around $Q$:
\[
 \bigl|
 \nabla^{j}_{h_{P,t},g_{\phi}}
 (s_{\kappa,t})
 \bigr|_{h_{P,t},g_{\phi}}
 \leq
 Ct^je^{-\gamma'\rho_{P,\kappa}(Q)t}.
\]
Then, we obtain the claim of the lemma.
\hfill\qed

\begin{lem}
\label{lem;23.9.30.31}
Let $P\in D_{\theta}$.
For any $j\geq 0$
and $0<\gamma<2\sqrt{2}$,
there exists $C_{j,\gamma}>0$ such that
the following holds on 
$X^{(2)}_{P,\kappa}\setminus X_P(2\delta(\kappa))$:
\[
\Bigl|
 \nabla_{\infty}^j
 (s(h_{\infty},h_t)-\id)
 \Bigr|_{h_{\infty},g_{\phi}}
 \leq
 C_{j,\gamma}
 e^{-\gamma \rho_{P,\kappa} t}.
\]
\end{lem}
\pf
There exists $C_{j,\gamma}'>0$ such that
the following holds on 
$X^{(2)}_{P,\kappa}\setminus X_P(2\delta(\kappa))$:
\[
\Bigl|
 \nabla_{\infty}^j
 (s(h_{\infty},\htilde_{\kappa,t})-\id)
 \Bigr|_{h_{\infty},g_{\phi}}
 \leq
 C_{j,\gamma}'
 e^{-\gamma \rho_{P,\kappa} t}.
\]
Together with the estimate in Lemma \ref{lem;23.9.30.30},
we obtain the claim of Lemma \ref{lem;23.9.30.31}.
\hfill\qed

\vspace{.1in}

Let $\rho_{\kappa}:X\to \real_{\geq 0}$
be defined by
$\rho_{\kappa}(Q)=\min d_{\phi}(Q,X_P(2\delta(\kappa)))$.

\begin{cor}
\label{cor;23.9.30.40}
For any $j> 0$ and $0<\gamma<2\sqrt{2}$,
there exists $C_j>0$ such that 
the following holds  
on $X\setminus \bigcup_{P}X_P(2\delta(\kappa))$:
\[
\Bigl|
 \nabla_{\infty}^j
 (s(h_{\infty},h_t))
 \Bigr|_{h_{\infty},g_{\phi}}
 \leq
 C_{j,\gamma}
 e^{-\gamma \rho_{\kappa} t}.
\] 
\hfill\qed
\end{cor}

\subsubsection{Refinement}

We obtain the following.
\begin{lem}
\label{lem;23.9.30.41}
There exists $C>0$ such that 
the following holds 
on $X\setminus \bigcup_{P}X^{(2)}_{P,\kappa}$:
\[
 |s(h_{\infty},h_t)^{\bot}|
\leq Ce^{-4(\kappa+5\delta(\kappa)) t}.
\] 
\end{lem}
\pf
It follows from Proposition \ref{prop;23.9.27.10}
and Corollary \ref{cor;23.9.30.40}.
\hfill\qed

\begin{cor}
\label{cor;23.10.1.3}
There exists $C>0$ such that 
the following holds 
on $X\setminus \bigcup_{P}X^{(2)}_{P,\kappa}$:
\begin{equation}
 \label{eq;23.10.1.1}
 |s_{\kappa,t}^{\bot}|
\leq C e^{-4(\kappa+5\delta(\kappa)) t}.
\end{equation}
There exists $C'>0$ such that
\begin{equation}
\label{eq;23.10.1.2}
 \Bigl|
 \Tr\bigl(s_{\kappa,t}F(\htilde_{\kappa,t})\bigr)
 \Bigr|
\leq C'e^{-8(\kappa+4\delta(\kappa)) t}.
\end{equation}
There exists $C''>0$ such that
\begin{equation}
 \label{eq;23.10.1.3}
 \|\delbar_E(s_{\kappa,t})\|_{L^2,\htilde_t,g_X}
\leq Ce^{-4(\kappa+4\delta(\kappa)) t},
\quad
 \|[\theta,s_{\kappa,t}]\|_{L^2,\htilde_t,g_X}
 \leq Ce^{-4(\kappa+4\delta(\kappa)) t}.
\end{equation}
\end{cor}
\pf
We obtain (\ref{eq;23.10.1.1}) from Lemma \ref{lem;23.9.30.41}.
We obtain the others by repeating the argument
in \S\ref{subsection;23.9.30.50}.
\hfill\qed

\vspace{.1in}
We refine Lemma \ref{lem;23.10.1.2}.

\begin{lem}
For any $j>0$,
there exists $C_j>0$ such that 
the following holds on 
$X\setminus \bigcup X^{(2)}_{P,\kappa}$:
\begin{equation}
\label{eq;23.10.1.10}
 |\nabla_{\infty}^j(s_{\kappa,t})|_{h_{\infty},g_{\phi}}
 \leq
 C_je^{-4(\kappa+4\delta(\kappa)) t}.
\end{equation}
For any $j\geq 0$,
there exists $C'_j\geq 0$ such that 
the following holds on 
$X\setminus \bigcup X^{(2)}_{P,\kappa}$:
\begin{equation}
\label{eq;24.7.6.20}
 \bigl|
 \nabla_{\infty}^j([\theta,s_{\kappa,t}])
 \bigr|_{h_{\infty},g_{\phi}}
+\bigl|
 \nabla_{\infty}^j([\theta^{\dagger}_{h_t},s_{\kappa,t}])
 \bigr|_{h_{\infty},g_{\phi}}
 \leq
 C_je^{-4(\kappa+4\delta(\kappa)) t}.
\end{equation}
\end{lem}
\pf
By using Corollary \ref{cor;23.10.1.3},
we obtain that there exists $C>0$ such that
\[
 \Bigl|
 R(h_t)
 \Bigr|_{h_{\infty},g_{\phi}}
 \leq
 Ct^2e^{-4(\kappa+5\delta(\kappa))t}
\]
on $X\setminus \bigcup X^{(2)}_{P,\kappa}$.
Then, we obtain the estimate (\ref{eq;23.10.1.10})
by using Lemma \ref{lem;23.9.28.110}
and Lemma \ref{lem;23.9.28.111}.
We also obtain (\ref{eq;24.7.6.20})
from Corollary \ref{cor;24.7.6.11}.
\hfill\qed

\subsubsection{$C^0$-estimate}

\begin{prop}
\label{prop;23.10.1.20}
There exists $C>0$ such that 
the following holds on $X$:
\begin{equation}
\label{eq;23.10.1.20}
 |s_{\kappa,t}-1|_{\htilde_{\kappa,t}}
 \leq
 Ce^{-4(\kappa+4\delta(\kappa)) t}.
\end{equation}
\end{prop}
\pf
We already have the following estimates on
$X\setminus \bigcup X_{P,\kappa}^{(2)}$:
\[
 |s_{\kappa,t}^{\bot}|_{h_{\infty}}
\leq Ce^{-4(\kappa+5\delta(\kappa))t},
\quad\quad
 |\nabla_{\infty}(s_{\kappa,t})|_{h_{\infty}}
\leq Ce^{-4(\kappa+4\delta(\kappa))t}.
\]
Let $Q\in X\setminus \bigcup X_{P,\kappa}^{(2)}$.
Let $v_i$ $(i=1,2)$ be a holomorphic frame of $E_{Q,i}$ around $Q$
such that $h_{\infty}(v_i,v_i)=1$.
We have
\[
 \bigl|
 h_t(v_1,v_2)
 \bigr|
 \leq Ce^{-4(\kappa+4\delta(\kappa))t}.
\]
We also have
\[
\Bigl|
h_t(v_1,v_1)\cdot h_t(v_2,v_2)-1
\Bigr|
\leq Ce^{-4(\kappa +4\delta(\kappa))t},
\quad\quad
|d h_{t}(v_i,v_i)|\leq
Ce^{-4(\kappa +4\delta(\kappa))t}.
\]
Note that the inclusion
$X\setminus \bigcup X_{P,\kappa}^{(2)}\lrarr X$
induces the surjection of the fundamental groups.
Hence, there exists a path
$\gamma:[0,1]\to X\setminus \bigcup X_{P,\kappa}^{(2)}$
such that
(i) $\gamma(0)=\gamma(1)=Q$,
(ii) the two points of $\pi^{-1}(Q)$ are exchanged
by the monodromy along $\gamma$.
By considering the above estimates along the path,
we obtain
\[
\Bigl|
h_t(v_1,v_1)- h_t(v_2,v_2)
\Bigr|\leq
Ce^{-4(\kappa+4\delta(\kappa))t}.
\]
Then, we obtain
\[
\Bigl|
 h_t(v_i,v_i)-1
 \Bigr|
 \leq Ce^{-4(\kappa+4\delta(\kappa)) t}.
\]
Hence,
(\ref{eq;23.10.1.20}) holds for some $C>0$ on
$X\setminus \bigcup_PX_{P,\kappa}^{(2)}$.

On $X_{P,\kappa}^{(2)}$, 
the function
$\Tr(s_{\kappa,t}-\id)\geq 0$ is subharmonic.
Hence, (\ref{eq;23.10.1.20})
also holds on $\bigcup_P X_{P,\kappa}^{(2)}$.
\hfill\qed

\subsubsection{Estimates for higher derivatives}
\label{subsection;23.10.16.2}

We consider the connection
$\nabla_{\htilde_{\kappa,t},g_X}$
of $\End(E)\otimes T^{\ast}X$
induced by the Chern connection of $\htilde_{\kappa,t}$
and the Levi-Civita connection of $g_X$.

\begin{lem}
\label{lem;23.10.1.41}
For any $j>0$,
there exists $C_j>0$ such that 
\begin{equation}
\label{eq;23.10.1.31}
\Bigl|
 \nabla_{\htilde_{\kappa,t},g_X}^j
 \bigl(s_{\kappa,t}\bigr)
 \Bigr|_{\htilde_{\kappa,t},g_X}
 \leq
 C_jt^je^{-4(\kappa+3\delta(\kappa)) t}.
\end{equation}
For any $j\geq 0$,
there exists $C'_j>0$ such that 
\begin{equation}
\label{eq;24.7.6.21}
\Bigl|
 \nabla_{\htilde_{\kappa,t},g_X}^j
 \bigl([\theta,s_{\kappa,t}]\bigr)
 \Bigr|_{\htilde_{\kappa,t},g_X}
+\Bigl|
 \nabla_{\htilde_{\kappa,t},g_X}^j
 \bigl([\theta^{\dagger}_h,s_{\kappa,t}]\bigr)
 \Bigr|_{\htilde_{\kappa,t},g_X}
 \leq
 C_jt^je^{-4(\kappa+3\delta(\kappa)) t}.
\end{equation}
\end{lem}
\pf
By Proposition \ref{prop;23.10.1.20},
there exists $C>0$ such that
\[
\left|
 t^2[\theta,\theta^{\dagger}_{\htilde_{\kappa,t}}]
 -t^2[\theta,\theta^{\dagger}_{h_t}]
 \right|_{\htilde_{\kappa,t},g_X}
 \leq
 Ct^2e^{-4(\kappa+4\delta(\kappa)t)}.
\]
We obtain
\[
\Bigl|
 \delbar_E
 \bigl(
 s_{\kappa,t}^{-1}
 \del_{E,\htilde_{\kappa,t}} s_{\kappa,t}
 \bigr)
 \Bigr|_{\htilde_{\kappa,t},g_X}
 \leq
 Ct^2e^{-4(\kappa+4\delta(\kappa))t}.
\]
Then, we obtain (\ref{eq;23.10.1.31})
from Lemma \ref{lem;23.9.28.110},
Lemma \ref{lem;23.9.28.111}
and Corollary \ref{cor;23.10.1.3}.
We also obtain (\ref{eq;24.7.6.21})
by using Corollary \ref{cor;24.7.6.11}.
\hfill\qed

\vspace{.1in}
We obtain Theorem \ref{thm;23.10.1.40} from
Proposition \ref{prop;23.10.1.20}
and Lemma \ref{lem;23.10.1.41}.
\hfill\qed

\subsection{An auxiliary Hermitian product}

\subsubsection{Induced $1$-forms on the spectral curve}
\label{subsection;23.10.19.13}

We obtain the following $\End(E)$-valued $1$-form
on $X\setminus D_{\theta}$:
\[
 \Psi_t
 =(\del_{E,h_{\infty}}+t\theta^{\dagger}_{h_{\infty}})
-(\del_{E,h_t}+t\theta^{\dagger}_{h_t}).
\]
It satisfies
$(\delbar_E+t\ad\theta)\Psi_t=0$.
Hence, there exist
holomorphic $1$-forms
$\eta_t$ on
$\Sigma_{\theta|X\setminus D_{\theta}}$
such that
$F_{\eta_t}=(\Psi_t^{\circ})^{1,0}$.

\begin{lem}
\label{lem;23.10.19.1}
$\eta_t$ are holomorphic $1$-forms on $\Sigma_{\theta}$.
\end{lem}
\pf
Let $P\in D_{\theta}$.
We consider
the following $\End(E)$-valued $1$-form
on $X(\kappa_0)\setminus\{P\}$:
\[
\Psi'_{P,t}=
(\del_{E,h_{\infty}}+t\theta^{\dagger}_{h_{\infty}})-
 (\del_{E,h_{P,t}}+t\theta^{\dagger}_{h_{P,t}}).
\]
It is a section of
\[
 \End(E_{|X_P(\kappa_0)})^{\asym}\otimes\Omega^{1,0}
 \oplus
 \End(E_{|X_P}(\kappa_0))^{\sym}\otimes\Omega^{0,1}.
\]
Hence, we obtain $((\Psi'_{P,t})^{\circ})^{1,0}=0$.
We consider the following
$\End(E)$-valued $1$-form on $X_P(\kappa_0)$:
\[
  \Psi_{P,t}=
 (\del_{E,h_{P,t}}+t\theta^{\dagger}_{h_{P,t}}) 
-(\del_{E,h_{t}}+t\theta^{\dagger}_{h_{t}}).
\]
It satisfies
$(\delbar_E+t\ad\theta)\Psi_{P,t}=0$.
Because $\Psi_{P,t}$ is given on $X_P(\kappa_0)$,
there exist holomorphic $1$-forms $\eta_{P,t}$ 
on $\Sigma_{\theta|X_P(\kappa)}$
such that
$F_{\eta_{P,t}}=((\Psi_{P,t})^{1,0})^{\circ}$
on $X_P(\kappa_0)\setminus\{P\}$.
(See \cite{Mochizuki-Asymptotic-Hitchin-metric}.)
Because $\eta_t=\eta_{P,t}$
on $\Sigma_{\theta|X_P(\kappa_0)\setminus\{P\}}$,
we obtain the claim of Lemma \ref{lem;23.10.19.1}.
\hfill\qed

\begin{lem}
\label{lem;23.10.24.1}
$\eta_t=O(e^{-4\kappa t})$
for any $\kappa<\kappa_0$.
\end{lem}
\pf
Let
$W_{\kappa}=X\setminus \bigcup_PX^{(2)}_{P,\kappa}$.
Because
\[
 \Psi_t=
 s(h_{\infty},h_t)^{-1}
 \del_{h_{\infty}}
 s(h_{\infty},h_t)
+s(h_{\infty},h_t)^{-1}
 \bigl[
 \theta^{\dagger}_{h_{\infty}},
  s(h_{\infty},h_t)
 \bigr],
\]
we obtain
$\eta_t=O(e^{-4\kappa t})$
on $\Sigma_{\theta|W_{\kappa}}$.
Because $\eta_t$ are holomorphic,
it implies the estimate on $\Sigma_{\theta}$.
\hfill\qed

\begin{lem}
\label{lem;23.10.19.100}
If $(E,\theta)$ has a globally defined
non-degenerate symmetric pairing $C$,
then $\eta_t=0$.
\end{lem}
\pf
Because $\Psi_t$
is a section of
$\End(E)^{\asym}\otimes\Omega^{1,0}
\oplus
\End(E)^{\sym}\otimes\Omega^{0,1}$,
we obtain $(\Psi_t^{\circ})^{1,0}=0$.
It implies the claim of the lemma.
\hfill\qed

\subsubsection{Canonical choice of a root of $\phi$
and a root of the coordinate $z_P$}
\label{subsection;23.10.19.12}

We obtain the quadratic differential
$\pi^{\ast}\phi$ on $\Sigma_{\theta}$.
Let $a_{T^{\ast}X}$ denote the tautological $1$-form
on $T^{\ast}X$.
Let $\phi^{1/2}$ denote the pull back of $a_{T^{\ast}X}$
by the inclusion $\Sigma_{\theta}\to T^{\ast}X$.

\begin{lem}
We obtain $(\phi^{1/2})^2=\pi^{\ast}(\phi)$.
\hfill\qed
\end{lem}

For each $P\in D_{\theta}$,
let $z_P$ be a holomorphic coordinate on $X_P(\kappa_0)$
such that
$\phi_{|X_P(\kappa_0)}=(\frac{3}{2})^2z_P(dz_P)^2$.
It induces a holomorphic coordinate system
$(z_P,\xi_P)$ on $T^{\ast}X_P$
such that $\xi_P\,dz_P$
is the tautological $1$-form.
Because
\[
\Sigma_{\theta|X_P(\kappa_0)}
=\left\{
 (z_P,\xi_P)\,\left|\,
 \xi_P^2-\left(\frac{3}{2}\right)^2z_P=0
 \right.
\right\},
\]
$\xi_P$ is a holomorphic local coordinate
on $\Sigma_{\theta|X_P(\kappa_0)}$.
We have
\[
 \phi^{1/2}
=\frac{2^3}{3^2}
 \xi_P^2d\xi_P. 
\]

\subsubsection{Induced pairings on the cohomology group
of $\Sigma_{\theta}$}
\label{subsection;23.10.19.10}

For each $P\in D_{\theta}$,
let $\Ptilde\in \Sigma_{\theta}$
denote the pre-image of $P$.
We set $\Dtilde_{\theta}=
\{\Ptilde\,|\,P\in D_{\theta}\}$,
which can be regarded as
the divisor on $\Sigma_{\theta}$.
For each $\Ptilde\in \Dtilde_{\theta}$,
there exists the residue morphism
\[
 \Res^{(2)}_{\Ptilde}:
 H^0\bigl(
 \Sigma_{\theta},
 K_{\Sigma_{\theta}}^2(2\Dtilde)
 \bigr)
 \lrarr
 \cnum.
\]

Let $\nu_1,\nu_2\in H^0(\Sigma_{\theta},K_{\Sigma_{\theta}})$.
We obtain
\[
  \nu_1\cdot\nu_2\cdot\eta_t\cdot
  (\phi^{1/2})^{-1}
  \in H^0(\Sigma_{\theta},K_{\Sigma_{\theta}}^2(2D)).
\]
We define the pairing
$\langle\cdot,\cdot\rangle^{\aux}_{P,h_t}$
on $H^0(\Sigma_{\theta},K_{\Sigma_{\theta}})$
as follows:
\[
 \langle\nu_1,\nu_2\rangle^{\aux}_{P,h_t}
=-4\pi\Res^{(2)}_{\Ptilde}
 \Bigl(
 \nu_1\cdot\nu_2\cdot\eta_t\cdot
 (\phi^{1/2})^{-1}
 \Bigr).
\]
We set
\[
\langle\cdot,\cdot\rangle_{h_t}^{\aux}
=\sum_{P\in D_{\theta}}
\langle\cdot,\cdot\rangle^{\aux}_{P,h_t}.
\]

\subsubsection{Induced Hermitian pairings on
$T_{(E,t\theta)}\nbigm_H'$}

We identify
$T_{(E,t\theta)}\nbigm_H'
=H^1(X,\End(E,t\theta))$.
As recalled in \S\ref{subsection;24.7.6.10},
we have the decomposition
\[
 T_{(E,t\theta)}\nbigm_H'
 =(T_{(E,t\theta)}\nbigm_H')^{\ver}
 \oplus
 (T_{(E,t\theta)}\nbigm_H')^{\hor},
\]
and there exist the natural isomorphisms
\[
 \iota^{\ver}_t:
 H^1(\Sigma_{\theta},\nbigo_{\Sigma_{\theta}})
 \simeq
 (T_{(E,t\theta)}\nbigm_H')^{\ver},
 \quad
 \iota^{\hor}_t:
 H^0(\Sigma_{\theta},K_{\Sigma_{\theta}})
 \simeq
 (T_{(E,t\theta)}\nbigm_H')^{\hor}.
\]
(See \cite{Mochizuki-Asymptotic-Hitchin-metric} for more details.)
We define $g^{\aux}_{|(E,t\theta)}$
by using $\langle\cdot,\cdot\rangle^{\aux}_{h_t}$
in \S\ref{subsection;23.10.19.10} as follows.
\begin{itemize}
 \item $g^{\aux}_{|(E,t\theta)}(v_1,v_2)=0$
       if
       $v_i\in
       (T_{(E,t\theta)}\nbigm_H')^{\ver}$ $(i=1,2)$
       or 
       $v_i\in
       (T_{(E,t\theta)}\nbigm_H')^{\hor}$ $(i=1,2)$.
 \item Let $\nu\in H^0(\Sigma_{\theta},K_{\Sigma_{\theta}})$
       and
       $\tau\in H^1(\Sigma_{\theta},\nbigo_{\Sigma_{\theta}})
       =H^0(\Sigma_{\theta}^{\dagger},K_{\Sigma_{\theta}^{\dagger}})$.
       Then,
\[
       g^{\aux}_{|(E,t\theta)}(\iota^{\hor}_t(\nu),
       \iota^{\ver}_t(\tau))
       =
       \langle
       \nu,\overline{\tau}
       \rangle^{\aux}_{h_t}
\]
\end{itemize}

We obtain the following lemma from
Lemma \ref{lem;23.10.24.1}.
\begin{lem}
$g^{\aux}_{|(E,t\theta)}(\nu,\tau)
=O(e^{-4\kappa t})\|\nu\|_{L^2}\cdot\|\tau\|_{L^2}$
for any $0<\kappa<\kappa_0$. 
\hfill\qed
\end{lem}

\begin{cor}
\label{cor;23.10.19.11}
$g^{\aux}_{|(E,t\theta)}=O(e^{-4\kappa t})$
with respect to $g_{\semiflat|(E,t\theta)}$
for any $0<\kappa<\kappa_0$.
As a result,
there exists $t_0>0$ such that
 $g_{\semiflat|(E,t\theta)}
+g^{\aux}_{|(E,t\theta)}$
are positive definite Hermitian products 
on $T_{(E,t\theta)}\nbigm_H'$
for any $t\geq t_0$.
\hfill\qed
\end{cor}

\subsection{Comparison of the Hitchin metric and the semi-flat metric}

\subsubsection{Main theorem}

Let $t_0$ be as in Corollary \ref{cor;23.10.19.11}.
For any $t\geq t_0$,
let $\sigma_t$ be the automorphism of
$T_{(E,t\theta)}\nbigm_H$
determined by
$g_{H|(E,t\theta)}
=\bigl(
g_{\semiflat|(E,t\theta)}
+g^{\aux}_{|(E,t\theta)}
\bigr)
\cdot \sigma_t$.

\begin{thm}
\label{thm;23.10.8.10}
For any $0<\kappa<\kappa_0$,
there exists $C(\kappa)>0$ such that
\[
 |\sigma_t-\id|_{g_{\semiflat}}
\leq C(\kappa)e^{-8\kappa t}.
\]
\end{thm}

Let $\nu\in H^0(\Sigma_{\theta},K_{\Sigma_{\theta}})$.
It determines the horizontal class
\[
 \iota_t^{\hor}(\nu)
 \in H^1(X,\Def(E,t\theta))^{\hor}
 \subset
 H^1(X,\Def(E,t\theta)).
\]
(See \cite{Mochizuki-Asymptotic-Hitchin-metric}.)
Let $\ttH_t(\nu)$ denote the harmonic representative
of $\iota_t^{\hor}(\nu)$
with respect to $h_t$.
We shall prove the following proposition in
\S\ref{subsection;23.10.3.1}.
\begin{prop}
\label{prop;23.9.27.1}
For any $0<\kappa<\kappa_0$,
\[
 (\ttH_t(\nu),\ttH_t(\nu))_{L^2,h_t}
 =\|\nu\|_{L^2}^2
 +O\bigl(e^{-8\kappa t}\|\nu\|_{L^2}^2\bigr).
\]
\end{prop}

Let $\tau\in
H^0(\Sigma_{\theta}^{\dagger},K_{\Sigma_{\theta}^{\dagger}})
\simeq H^1(\Sigma_{\theta},\nbigo_{\Sigma_{\theta}})$.
We also regard $\tau$ as a harmonic $(0,1)$-form
on $\Sigma_{\theta}$.
We obtain
$\iota^{\ver}_t(\tau)
\in H^1(X,(\End(E),t\theta))$.
There exists a harmonic representative
$\ttV_t(\tau)$
of $(\End(E),t\theta,h_t)$.
We shall prove the following proposition in
\S\ref{subsection;23.10.16.3}.
\begin{prop}
\label{prop;23.9.27.2}
For any $0<\kappa<\kappa_0$,
\[
 (\ttV_t(\tau),\ttV_t(\tau))_{L^2,h_t}
=(\tau,\tau)_{L^2}
+O\bigl(e^{-8\kappa t}\|\tau\|_{L^2}^2\bigr).
\]
\end{prop}
We shall prove the following proposition
in \S\ref{subsection;23.10.16.4}.

\begin{prop}
\label{prop;23.9.27.4}
$(\ttH_t(\nu),\ttV_t(\tau))_{L^2,h_t}
=\langle \nu,\overline{\tau}\rangle^{\aux}_{h_t}
+O\bigl(e^{-8\kappa t}\cdot
 \|\nu\|_{L^2}\cdot
 \|\tau\|_{L^2}\bigr)$.
\end{prop}

We obtain Theorem \ref{thm;23.10.8.10}
from Proposition \ref{prop;23.9.27.1},
Proposition \ref{prop;23.9.27.2}
and Proposition \ref{prop;23.9.27.4}.

\subsubsection{Preliminary}

Let $P\in D_{\theta}$.
Let $s_{P,t}$ be
determined by $h_{P,t}=h_t\cdot s_{P,t}$ on $X_P(\kappa_0)$.
(Note that we use a different notation from
\S\ref{subsection;23.10.16.5}.)
We set
\[
\Psi_{P,t}:=
 \bigl(\del_{h_{P,t}}
 +t\theta^{\dagger}_{h_{P,t}}\bigr)
-\bigl(\del_{h_t}+t\theta^{\dagger}_{h_t}\bigr)
 =s_{P,t}^{-1}\del_{h_t}s_{P,t}
 +s_{P,t}^{-1}[t\theta^{\dagger}_{h_t},s_{P,t}].
\]
Because
$[\delbar_E+t\theta,
\del_{h_{P,t}}
+t\theta^{\dagger}_{h_{P,t}}]=0$
and 
$[\delbar_E+t\theta,
\del_{h_t}+t\theta^{\dagger}_{h_t}]=0$,
we obtain $(\delbar_E+t\ad\theta)\Psi_{P,t}=0$.

Let $f_P$ be the endomorphism determined by
$\theta=f_P\,dz_P$.
Let $M_{P,i,t}$ be the endomorphisms of $E_{|X_P(\kappa_0)}$
defined by
\[
 \Psi_{P,t}=M_{P,1,t}\,dz_P
 +M_{P,2,t}\,d\zbar_P.
\]
Because $(\delbar_E+t\ad\theta)\Psi_{P,t}=0$,
we obtain
\begin{equation}
\label{eq;23.10.8.1}
 \delbar M_{P,1,t}-[tf_P,M_{P,2,t}\,d\zbar_P]=0.
\end{equation}

\begin{lem}
\label{lem;23.10.16.10}
We have $\delbar M_{P,1,t}^{\circ}=0$.
There exist  $C_j>0$ such that
       the following holds on
       $X^{(1)}_{P,\kappa}\setminus X^{(2)}_{P,\kappa}$:
\[
       |\nabla_{\infty}^jM_{P,1,t}|_{h_{\infty},g_{\phi}}
       \leq C_je^{-4(\kappa+\delta(\kappa)) t}.
\]       
\end{lem}
\pf
The first claim follows from (\ref{eq;23.10.8.1}).
The second follows from Theorem \ref{thm;23.10.1.40}.
\hfill\qed

\vspace{.1in}
Let $\eta_t$ be the $1$-forms on $\Sigma_{\theta}$
as in \S\ref{subsection;23.10.19.13}.
Let $(z_P,\xi_P)$ be the coordinate as in
\S\ref{subsection;23.10.19.12}.
Let $\upsilon_{P,t}$ be the holomorphic function on
$\Sigma_{\theta|X_P(\kappa_0)\setminus\{P\}}$
determined by
\[
 \eta_t=\upsilon_{P,t}\cdot \pi^{\ast}(dz_P).
\]
We note that $\upsilon_{P,t}$
has at most pole of order $1$ at $\pi^{-1}(P)$.
By the construction, we have the following lemma.
\begin{lem}
$M_{P,1,t}^{\circ}
=F_{\upsilon_{P,t}}$.
\hfill\qed
\end{lem}

\subsubsection{Some $1$-forms around $P\in D_{\theta}$}

Let $\nu\in H^0(\Sigma_{\theta},K_{\Sigma_{\theta}})$.
We obtain the holomorphic section $F_{\nu}$ of $\End(E)\otimes K_X$
on $X\setminus D_{\theta}$.

Let $P\in D_{\theta}$.
There exists a holomorphic function $\alpha_P$ on
$\Sigma_{\theta|X_P(\kappa)}$
such that $d\alpha_P=\nu$
and that $\alpha_P=0$ at $\pi^{-1}(P)$.
We obtain the induced section $F_{\alpha_P}$
of $\End(E_P)$ on $X_P(\kappa_0)$.
Let $\alpha_0$ and $\alpha_1$ be the holomorphic functions
on $X_P(\kappa_0)$
determined by
\[
 \alpha_P(\xi_P)=
 \alpha_0(z_P)
+\alpha_1(z_P)\xi_P,
\]
where $(z_P,\xi_P)$ be the holomorphic coordinate system
as in \S\ref{subsection;23.10.19.12}.
Then, the following holds:
\[
 F_{\alpha_P}=\alpha_0\id+\alpha_1f_P.
\]

\begin{lem}
There exists $C>0$
independently from $\nu\in H^0(\Sigma_{\theta},K_{\Sigma_{\theta}})$
and $t\geq 1$
such that the following holds:
\[
 \|F_{\alpha_P}\|^2_{h_{P,t},X_P(\kappa_0)}
 \leq C\|\nu\|_{L^2}^2.
\]
\end{lem}
\pf
We have
$|\alpha_P|\leq C_0\|\nu\|_{L^2}$.
By Theorem \ref{thm;23.10.1.40},
there exists $C_1>0$ such that
$|F_{\alpha_P}|_{h_{P,t}}^2\leq C_1\|\nu\|_{L^2}^2$
on $\del X_P^{(2)}$
for any $t\geq 1$.
Because $(\delbar_E+t\ad\theta)F_{\alpha_P}=0$,
we obtain
$-\del_z\del_{\zbar}\bigl|
 F_{\alpha_P}
 \bigr|_{h_t}^2
 \leq 0$.
Then, we obtain the claim of the lemma.
\hfill\qed

\vspace{.1in}

We obtain the following harmonic $1$-forms 
of $(\End(E_P),t\theta_P,h_{P,t})$
on $X_P(\kappa_0)$:
\[
\ttH_{P,t}(\nu):=
\bigl(
\del_{h_{P,t}}
+\ad(t\theta^{\dagger}_{h_{P,t}})
\bigr) F_{\alpha_P}.
\]
According to \cite[Lemma 4.14]{Mochizuki-Asymptotic-Hitchin-metric},
there exists $C>0$ such that
the following holds on $X^{(1)}_{P,\kappa}$
for any $t\geq 1$:
\[
|\ttH_{P,t}(\nu)|_{h_{P,t},g_X}\leq C(1+t)\|\nu\|_{L^2}.
\]
We remark the following.
\begin{lem}
The $(1,0)$-part $\ttH_{P,t}(\nu)^{1,0}$ is symmetric,
and the $(0,1)$-part $\ttH_{P,t}(\nu)^{0,1}$
is anti-symmetric with respect to $C_P$. 
We have
$\ttH_{P,t}(\nu)^{\circ}=F_{\nu}$
on $X_P(\kappa_0)$.
\hfill\qed
\end{lem}

On $X_P(\kappa_0)\setminus \{P\}$,
there exists
a $C^{\infty}$-section $\rho_{P,t}(\nu)$
of $\End(E)$,
which is anti-symmetric with respect to $C_P$,
such that
\[
 \ttH_{P,t}(\nu)-F_{\nu}
=(\delbar_E+t\ad\theta) \rho_{P,t}(\nu).
\]
On $X_{P,\kappa}^{(1)}\setminus X_{P,\kappa}^{(2)}$,
we have
\begin{equation}
\label{eq;24.7.6.30}
 \|\rho_{P,t}(\nu)\|_{L_{\ell}^p,h_{\infty},g_{\phi}}
 =O\bigl(e^{-4(\kappa+7\delta(\kappa)) t}
 \|\nu\|_{L^2}\bigr).
\end{equation}
We also note that $\rho_{P,t}(\nu)^{\circ}=0$.

We also obtain the following harmonic $1$-forms
of $(\End(E),t\ad\theta,h_t)$ on $X_P(\kappa_0)$:
\[
 \ttHtilde_{P,t}(\nu)
:=(\del_{E,h_{t}}+\ad t\theta^{\dagger}_{h_{t}})F_{\alpha_P}.
\]
We set
$\ttI_{P,t}(\nu)=\ttH_{P,t}(\nu)-\ttHtilde_{P,t}(\nu)$.

\begin{lem}
We have
\begin{equation}
\label{eq;23.10.7.1}
 \ttI_{P,t}(\nu)
=\bigl(\delbar+t\ad\theta\bigr)(-t^{-1}\alpha_1 M_{P,1,t}).
\end{equation}
\end{lem}
\pf
We have
\[
\ttI_{P,t}(\nu)=[\Psi_t,F_{\alpha_P}]
=[M_{P,1,t}dz_P,\alpha_1f_P]
+[M_{P,2,t}\,d\zbar_P,\alpha_1f_P]
=-[tf_P\,dz_P,t^{-1}\alpha_1M_{P,1,t}]
-[f_P,\alpha_1M_{P,2,t}d\zbar_P].
\]
We also have
\[
 \delbar(t^{-1}\alpha_1 M_{P,1,t})
 =t^{-1}\alpha_1\delbar M_{P,1,t}
 =t^{-1}\alpha_1[tf_P,M_{P,2,t}d\zbar_P]
 =[f_P,\alpha_1M_{P,2,t}d\zbar_P].
\]
Then, we obtain (\ref{eq;23.10.7.1}).
\hfill\qed

\begin{lem}
\label{lem;23.10.16.11}
We have $\ttI_{P,t}(\nu)^{\circ}=0$.
There exists $C>0$ such that
 $|\ttI_{P,t}(\nu)|_{h_{\infty},g_{\phi}}
 \leq Ce^{-4\kappa t}\|\nu\|_{L^2}$
on $X_{P,\kappa}^{(1)}\setminus X_{P,\kappa}^{(2)}$.
\end{lem}
\pf
The first claim is clear by the construction.
The second follows from
Lemma \ref{lem;23.10.16.10}.
\hfill\qed

\subsubsection{Approximations of harmonic $1$-forms in the horizontal direction}
\label{subsection;23.10.8.11}

Let $\nu\in H^0(\Sigma_{\theta},K_{\Sigma_{\theta}})$.
We recall a construction in \cite{Mochizuki-Asymptotic-Hitchin-metric}.
Let $\chi_{P,\kappa}:
X\to[0,1]$ be a $C^{\infty}$-function
such that
$\chi_{P,\kappa}=1$ on $X^{(2)}_{P,\kappa}$
and $\chi_{P,\kappa}=0$ outside
$X^{(1)}_{P,\kappa}$.
By patching
$\ttH_{P,t}(\nu)$ $(P\in D_{\theta})$
and $F_{\nu}$
by using $\rho_{P,t}(\nu)$,
we construct
a section $\ttH'_{\kappa,t}(\nu)$ of
$\End(E)\otimes(\Omega^{1,0}\oplus\Omega^{0,1})$
as follows:
\begin{itemize}
 \item 
On $X\setminus \bigcup X^{(1)}_{P,\kappa}$,
we set
$\ttH'_{\kappa,t}(\nu)=F_{\nu}$.
\item
     On $X^{(1)}_{P,\kappa}\setminus
     X^{(2)}_{P,\kappa}$,
we set
$\ttH'_{\kappa,t}(\nu)=F_{\nu}+
 (\delbar_E+t\ad\theta)(\chi_{P,\kappa} \rho_{P,t}(\nu))$.
\item
     On $X^{(2)}_{P,\kappa}$,
we set
$\ttH'_{\kappa,t}(\nu)=\ttH_{P,t}(\nu)$.
\end{itemize}

We obtain the following lemma by the construction.

\begin{lem}
\label{lem;23.10.8.20}
$\ttH'_{\kappa,t}(\nu)^{\circ}=F_{\nu}$.
\hfill\qed 
\end{lem}

\begin{lem}
We have $(\delbar_E+t\ad\theta)\ttH'_{\kappa,t}(\nu)=0$,
and 
$(\delbar_E+t\ad\theta)^{\ast}_{g_X,h_t}
 \ttH'_{\kappa,t}(\nu)
 =O\bigl(e^{-4\kappa t}\|\nu\|_{L^2}\bigr)$.
The cohomology class is $\iota^{\hor}_t(\nu)$.
\end{lem}
\pf
By the construction, we have
$(\delbar_E+t\ad\theta)\ttH'_{\kappa,t}(\nu)=0$.
On $X\setminus D_{\theta}$,
we have 
$\ttH'_{\kappa,t}(\nu)^{\circ}=F_{\nu}$.
On $X^{(2)}_{P,\kappa}$, $\ttH'_{\kappa,t}(\nu)^{1,0}$
is self-adjoint,
and $\ttH'_{\kappa,t}(\nu)^{0,1}$
is anti-self-adjoint with respect to $C_P$.
Hence, as studied in \cite{Mochizuki-Asymptotic-Hitchin-metric},
the cohomology class
of $\ttH'_{\kappa,t}(\nu)$
is $\iota^{\hor}_t(\nu)$.

We have
$(\delbar_E+t\ad\theta)^{\ast}_{g_X,h_t}
 \ttH'_{\kappa,t}(\nu)
 =O\bigl(e^{-4(\kappa+7\delta(\kappa)) t}\|\nu\|_{L^2}\bigr)$
on $X^{(1)}_{P,\kappa}\setminus X^{(2)}_{P,\kappa}$
by the construction. 
On $X^{(2)}_{P,\kappa}$,
we have
\[
 \ttH'_{\kappa,t}(\nu)
=(\del_{E,h_{P,t}}+\ad t\theta^{\dagger}_{h_{P,t}})F_{\alpha_P}
=(\del_{E,h_{t}}+\ad t\theta^{\dagger}_{h_{t}})F_{\alpha_P}
+\Bigl[
 s_{\kappa,t}^{-1}\del_{E,h_{t}}(s_{\kappa,t})
+s_{\kappa,t}^{-1}[t\theta^{\dagger}_{h_{t}},s_{\kappa,t}],\,
 F_{\alpha_P}
 \Bigr].
\]
By Theorem \ref{thm;23.10.1.40},
we obtain
\[
 (\del_{E,h_{t}}+t\ad\theta^{\dagger}_{h_t})
 \ttH'_{\kappa,t}(\nu)
 =(\del_{E,h_{t}}+t\ad\theta^{\dagger}_{h_t})
 \left(
 \Bigl[
 s_{\kappa,t}^{-1}\del_{E,h_{t}}(s_{\kappa,t})
+s_{\kappa,t}^{-1}[t\theta^{\dagger}_{h_{t}},s_{\kappa,t}],\,
 F_{\alpha_P}
 \Bigr]
 \right)
=O\bigl(e^{-4\kappa t}\|\nu\|_{L^2}\bigr).
\]
We obtain the claim of the lemma.
\hfill\qed

\vspace{.1in}
Note that
$\Tr\bigl(
 (\del_{E,h_t}+t\ad\theta^{\dagger}_{h_t})
 \ttH'_{\kappa,t}(\nu)
\bigr)
=\del\Tr(\ttH'_{\kappa,t}(\nu))=0$.
There exist a unique section
$\gamma_{\kappa,t}(\nu)$ of $\End(E)$
satisfying
$\Tr\gamma_{\kappa,t}(\nu)=0$
and
\[
(\delbar_E+t\ad\theta)^{\ast}_{g_X,h_t}
(\delbar_E+t\ad\theta)\gamma_{\kappa,t}(\nu)=
(\delbar_E+t\ad\theta)^{\ast}_{g_X,h_t}
 \ttH'_{\kappa,t}(\nu).
\]
The following lemma is clear.
\begin{lem}
$\ttH_t(\nu)
=\ttH'_{\kappa,t}(\nu)-\gamma_{\kappa,t}(\nu)$.
\hfill\qed
\end{lem}

By \cite[Proposition 2.39]{Mochizuki-Asymptotic-Hitchin-metric},
we have the following estimate:
\[
 \|\gamma_{\kappa,t}(\nu)\|_{L^2,h_t,g_X}
+\|(\delbar_E+t\theta)\gamma_{\kappa,t}(\nu)\|_{L^2,h_t}
\leq Ce^{-4(\kappa+\delta(\kappa)) t}.
\]

\subsubsection{Proof of Proposition \ref{prop;23.9.27.1}}
\label{subsection;23.10.3.1}

\begin{lem}
\label{lem;23.9.27.3}
$(\ttH_t(\nu),\ttH_t(\nu))_{L^2,h_t}
 =(\ttH'_{\kappa,t}(\nu),
 \ttH'_{\kappa,t}(\nu))_{L^2,h_t}
+O\bigl(e^{-8\kappa t}\|\nu\|^2_{L^2}\bigr)$.
\end{lem}
\pf
By the standard formula (see Lemma \ref{lem;23.10.7.20} below),
we obtain
\[
\bigl(\ttH_t(\nu),
(\delbar_E+t\ad\theta)\gamma_{\kappa,t}(\nu)
\bigr)_{L^2,h_t}
=\bigl((\delbar_E+t\ad\theta)^{\ast}_{h_t,g_X}\ttH_t(\nu),
\gamma_{\kappa,t}(\nu)
\bigr)_{L^2,h_t,g_X}
=0,
\]
\[
 \bigl(\ttH'_{\kappa,t}(\nu),
 (\delbar_E+t\ad\theta)\gamma_{\kappa,t}(\nu)\bigr)_{L^2,h_t}
 =\bigl(
 (\delbar_E+t\ad\theta)^{\ast}_{h_t,g_X}\ttH'_{\kappa,t}(\nu),
  \gamma_{\kappa,t}(\nu)
  \bigr)_{L^2,h_t,g_X}
=O\bigl(e^{-8\kappa t}\|\nu\|_{L^2}^2\bigr).
\]
Hence, we obtain 
\[
 \bigl(
 \ttH_t(\nu),\ttH_t(\nu)
  \bigr)_{L^2,g_X,h_t}
= \bigl(
 \ttH_t(\nu),
 \ttH'_{\kappa,t}(\nu)
  \bigr)_{L^2,g_X,h_t}
\\
=\bigl(
 \ttH'_{\kappa,t}(\nu),
 \ttH'_{\kappa,t}(\nu)
  \bigr)_{L^2,g_X,h_t}
+O\bigl(e^{-8\kappa t}\|\nu\|_{L^2}^2\bigr).
\]
\hfill\qed

\vspace{.1in}
We set
$W_{\kappa}=X\setminus \bigcup_{P}X^{(2)}_{P,\kappa}$.

\begin{lem}
\[
 \bigl(
 \ttH'_{\kappa,t}(\nu),
\ttH'_{\kappa,t}(\nu)
 \bigr)_{L^2,h_t,W_{\kappa}}
=2\sqrt{-1}\int_{\Sigma_{\theta|W_{\kappa}}}
 \nu\cdot\nubar
+O\bigl(e^{-8\kappa t}\|\nu\|_{L^2}^2\bigr).
\]
\end{lem}
\pf
It follows from
the estimate (\ref{eq;24.7.6.30}),
Lemma \ref{lem;23.10.8.20} and Lemma \ref{lem;23.9.25.1} below.
\hfill\qed

\vspace{.1in}

On $X_{P,\kappa}^{(2)}$,
we have
$\ttH_{\kappa,t}(\nu)=\ttH_{P,t}(\nu)$.

\begin{lem}
\begin{equation}
\label{eq;23.10.8.21}
\bigl(
\ttH_{P,t}(\nu),\ttH_{P,t}(\nu)
\bigr)_{L^2,h_t,X_{P,\kappa}^{(2)}}
= 
\bigl(
\ttHtilde_{P,t}(\nu),\ttHtilde_{P,t}(\nu)
\bigr)_{L^2,h_t,X_{P,\kappa}^{(2)}}
+O\bigl(e^{-8\kappa t}\|\nu\|_{L^2}^2\bigr).
\end{equation}
\end{lem}
\pf
Because
$(\delbar+t\ad\theta)\ttI_{P,t}(\nu)=0$,
we have
$(\del_{h_t}+t\ad\theta^{\dagger}_{h_t})^{\ast}_{h_t,g_X}
\ttI_{P,t}(\nu)=0$.
Hence, by Stokes formula (see Lemma \ref{lem;23.10.7.20} below),
and by Lemma \ref{lem;23.10.16.11} and
Lemma \ref{lem;23.9.25.1} below,
we obtain
\[
\bigl(
\ttHtilde_{P,t}(\nu),\ttI_{P,t}(\nu)
\bigr)_{L^2,h_t,X_{P,\kappa}^{(2)}}
=
2\sqrt{-1}
\int_{\del X^{(2)}_{P,\kappa}}
\Tr\Bigl(
 F_{\alpha_P}\cdot
 \bigl(
 \ttI_{P,t}(\nu)^{1,0}
 \bigr)^{\dagger}_{h_t}
 \Bigr)
=O\bigl(e^{-8\kappa t}\bigr)\|\nu\|_{L^2}^2.
\]
We also have
$\bigl(
\ttI_{P,t}(\nu),\ttI_{P,t}(\nu)
\bigr)_{L^2,h_t,X_{P,\kappa}^{(2)}}
=O\bigl(e^{-8\kappa t}\bigr)\|\nu\|_{L^2}^2$.
Then, we obtain (\ref{eq;23.10.8.21}).
\hfill\qed

\begin{lem}
\begin{equation}
\label{eq;23.10.6.2} 
\bigl(
\ttHtilde_{P,t}(\nu),\ttHtilde_{P,t}(\nu)
\bigr)_{L^2,h_t,X_{P,\kappa}^{(2)}}
=2\sqrt{-1}\int_{\Sigma_{\theta|X^{(2)}_{P,\kappa}}}
 \nu\nubar
+O\bigl(e^{-8\kappa t}\|\nu\|_{L^2}^2\bigr).
\end{equation}
\end{lem}
\pf
By using
Lemma \ref{lem;23.10.7.20} and
Lemma \ref{lem;23.9.25.1} below,
we obtain
\begin{multline}
 \bigl(
\ttHtilde_t(\nu),\ttHtilde_t(\nu)
 \bigr)_{L^2,h_t,X_{P,\kappa}^{(2)}}
=2\sqrt{-1}\int_{\del X^{(2)}_{P,\kappa}}
 \Tr\bigl(
 F_{\alpha}
 (\ttHtilde_t(\nu)^{1,0})^{\dagger}_{h_t}
 \bigr)
\\
 =2\sqrt{-1}\int_{\del\Sigma_{\theta|X^{(2)}_{P,\kappa}}}
 \alpha_P \nubar
+O\bigl(e^{-8\kappa t}\|\nu\|_{L^2}^2\bigr)
= 2\sqrt{-1}\int_{\Sigma_{\theta|X^{(2)}_{P,\kappa}}}
 \nu \nubar
+O\bigl(e^{-8\kappa t}\|\nu\|_{L^2}^2\bigr).
\end{multline}
We obtain (\ref{eq;23.10.6.2}).
\hfill\qed

\vspace{.1in}
In all, we obtain Proposition \ref{prop;23.9.27.1}.
\hfill\qed

\subsubsection{Approximation of harmonic $1$-forms in the vertical direction}
\label{subsection;23.10.8.12}

Let
$\tau\in
H^0(\Sigma_{\theta}^{\dagger},K_{\Sigma_{\theta}^{\dagger}})
\simeq H^1(\Sigma_{\theta},\nbigo_{\Sigma_{\theta}})$.
We recall the construction of approximations of
$\ttV_t(\tau)$ introduced in \cite{Mochizuki-Asymptotic-Hitchin-metric}.

Let $P\in D_{\theta}$.
There exist an anti-holomorphic function $\beta_P$
on $\Sigma_{\theta|X_P(\kappa_0)}$
such that $d\beta_P=\tau$.
We obtain the endomorphism
$F^{\dagger}_{\beta_P}$ of $E_P$
by using
the spectral curve $\Sigma^{\dagger}_{\theta}\to X^{\dagger}$
of $(E,\theta^{\dagger}_{h_t})$.
It equals
$(F_{\betabar_P})^{\dagger}_{h_t}$.

Let $\chi_{P,\kappa}$ be as in \S\ref{subsection;23.10.3.1}.
We put
\[
 \tau_{\kappa}^{\circ}
=\tau-\sum \delbar(\chi_{P,\kappa}\beta_P)
\]
on $\Sigma_{\theta}$.
We set
\[
 \ttV'_{\kappa,t}(\tau)
 =F_{\tau_{\kappa}^{\circ}}
+\sum (\delbar_E+t\ad\theta)(\chi_{P,\kappa} F^{\dagger}_{\beta_{P}}).
\]
We have
\[
 (\delbar_E+t\ad\theta)\ttV'_{\kappa,t}(\tau)=0,
 \quad\quad
 (\del_{E,h_t}+t\ad\theta^{\dagger}_{h_t})\ttV'_{\kappa,t}(\tau)
 =O\bigl(e^{-4(\kappa+7\delta(\kappa)) t}\|\tau\|_{L^2}\bigr).
\]
We note that
$\Tr\bigl(
 (\del_{E,h_t}+t\ad\theta^{\dagger}_{h_t})
 \ttV'_{\kappa,t}(\tau)
 \bigr)
 =\del\Tr \ttV'_{\kappa,t}(\tau)=0$
(see \cite[Lemma 4.28]{Mochizuki-Asymptotic-Hitchin-metric}).
There exist $\gamma_{\kappa,t}(\tau)$
satisfying
$\Tr\gamma_{\kappa,t}(\tau)=0$
and 
\[
 (\delbar_E+t\ad\theta)(\del_{E,h_t}+t\ad\theta^{\dagger}_{h_t})
 \gamma_{\kappa,t}(\tau)
=(\del_{E,h_t}+t\ad\theta^{\dagger}_{h_t})\ttV'_{\kappa,t}(\tau).
\]
By \cite[Proposition 2.39]{Mochizuki-Asymptotic-Hitchin-metric},
we have
\[
 \|\gamma_{\kappa,t}\|_{L^2,h_t,g_X}
+\|(\delbar_E+t\theta)\gamma_{\kappa,t}(\tau)\|_{L^2,h_t}
 =O\bigl(e^{-4\kappa t}\|\tau\|_{L^2}\bigr).
\]
We have
\[
 \ttV_{t}(\tau)
 =\ttV'_{\kappa,t}(\tau)
 -(\delbar_E+t\ad\theta)\gamma_{\kappa,t}(\tau).
\]

\subsubsection{Proof of Proposition \ref{prop;23.9.27.2}}
\label{subsection;23.10.16.3}

We obtain the following estimate
as in the case of Lemma \ref{lem;23.9.27.3}:
\[
\bigl(
 \ttV_t(\tau),\ttV_t(\tau)
 \bigr)_{L^2,h_t}
=\bigl(
 \ttV'_{\kappa,t}(\tau),\ttV'_{\kappa,t}(\tau)
 \bigr)_{L^2,h_t}
+O\bigl(e^{-8\kappa t}\|\tau\|_{L^2}^2\bigr).
\]

By using Lemma {\rm\ref{lem;23.9.26.1}} below,
we obtain the following:
\[
 \bigl(
 \ttV'_{\kappa,t}(\tau),
 \ttV'_{\kappa,t}(\tau)
 \bigr)_{L^2,h_t,W_{\kappa}}
=-2\sqrt{-1}\int_{\Sigma_{\theta|W_{\kappa}}}
 \tau\cdot\taubar
+O\bigl(e^{-8\kappa t}\|\tau\|_{L^2}^2\bigr).
\] 
By using the Stokes formula (Lemma \ref{lem;23.10.7.20}),
we obtain
\[
 \bigl(
 \ttV'_{\kappa,t}(\tau),
 \ttV'_{\kappa,t}(\tau)
 \bigr)_{L^2,h_t,X^{(2)}_{P,\kappa}}
=-2\sqrt{-1}\int_{\del X^{(2)}_{P,\kappa}}
\Tr\bigl(F^{\dagger}_{\beta_P}
 (\ttV'_{\kappa,t}(\tau)^{0,1})^{\dagger}_{h_t}
\bigr).
\]
By using Lemma \ref{lem;23.9.26.1} below,
we obtain
\[
\int_{\del X^{(2)}_{P,\kappa}}
\Tr\bigl(
 F^{\dagger}_{\beta_P}
 (\ttV'_{\kappa,t}(\tau)^{0,1})^{\dagger}_{h_t}
 \bigr)
=\int_{\del \Sigma_{\theta|X^{(2)}_{P,\kappa}}}
 \beta_P\taubar
 +O\bigl(e^{-8\kappa t}\|\tau\|_{L^2}^2\bigr)
=\int_{\Sigma_{\theta|X^{(2)}_{P,\kappa}}}
 \tau \taubar
 +O\bigl(e^{-8\kappa t}\|\tau\|_{L^2}^2\bigr).
\]
In all, we obtain
Proposition \ref{prop;23.9.27.2}.
\hfill\qed

\subsubsection{Proof of Proposition \ref{prop;23.9.27.4}}
\label{subsection;23.10.16.4}

As before, we obtain
\[
(\ttH_t(\nu),\ttV_t(\tau))_{L^2,h_t}
=(\ttH'_{\kappa,t}(\nu),\ttV'_{\kappa,t}(\tau))_{L^2,h_t}
+O\bigl(e^{-8\kappa t}\cdot
 \|\nu\|_{L^2}\cdot
 \|\tau\|_{L^2}\bigr).
\]

By the construction, we have
$\bigl(
 \ttH'_{\kappa,t}(\nu)^{0,1}
 \bigr)^{\circ}=0$
and 
$\bigl(
\ttV'_{\kappa,t}(\tau)^{1,0}
\bigr)^{\circ}=0$.
Hence, we obtain the following
by using Lemma \ref{lem;23.9.25.1} below:
\[
 \bigl(
 \ttH'_{\kappa,t}(\nu),
 \ttV'_{\kappa,t}(\tau)
 \bigr)_{L^2,h_t,W_{\kappa}}
=O\bigl(e^{-8\kappa t}\|\tau\|_{L^2}\|\nu\|_{L^2}\bigr).
\]

On $X^{(2)}_{P,\kappa}$,
we have
$\ttH'_{\kappa,t}(\nu)=\ttH_{P,t}(\kappa)$.
To study
$\bigl(
 \ttH_{P,t}(\nu),
 \ttV_{\kappa,t}'(\tau)
 \bigr)_{L^2,h_t,X^{(2)}_{P,\kappa}}$,
we prepare the following lemma. 

\begin{lem}
\label{lem;23.10.7.3}
 \begin{equation}
\label{eq;23.10.7.2}
 \bigl(
 \ttI_{P,t}(\nu),
 \ttV'_{\kappa,t}(\tau)
\bigr)_{L^2,h_t,X_{P,\kappa}^{(2)}}
=
\langle \nu,\overline{\tau}\rangle^{\aux}_{P,h_t}
+O\bigl(e^{-8\kappa t}
\|\nu\|_{L^2}\cdot\|\tau\|_{L^2}\bigr).
\end{equation}
\end{lem}
\pf
By (\ref{eq;23.10.7.1})
and Lemma \ref{lem;23.10.7.20} below,
we have
\[
 \bigl(
 \ttI_{P,t}(\nu),
 \ttV'_{\kappa,t}(\tau)
\bigr)_{L^2,h_t,X_{P,\kappa}^{(2)}} 
=2\sqrt{-1}\int_{\del X_{P,\kappa}^{(2)}}
 \Tr\Bigl(
 (t^{-1}\alpha_1M_{P,1,t})\cdot
 \del_{h_t}F_{\betabar_P}
 \Bigr).
\]
On $\del X^{(2)}_{P,\kappa}$,
we have
$\bigl(
 \del_{h_t}F_{\betabar_P}
 \bigr)^{\bot}
  =O\bigl(e^{-4\kappa t}\|\tau\|_{L^2}\bigr)$.
We also have
$(\del_{h_t}F_{\betabar_P})^{\circ}=F_{d\betabar_P}$.
Hence, 
we obtain
\begin{multline}
 \int_{\del X_{P,\kappa}^{(2)}}
 \Tr\Bigl(
 (t^{-1}\alpha_1M_{P,1,t})^{\circ}
 (\del_{h_t}F_{\betabar_P})^{\circ}
 \Bigr)
=\int_{\del X_{P,\kappa}^{(2)}}
 \Tr\Bigl(
 t^{-1}\alpha_1F_{\upsilon_{P,t}}
 F_{d\betabar_P}
 \Bigr)
 =\int_{\del \Sigma_{\theta|X_{P,\kappa}^{(2)}}}
 t^{-1}\alpha_1\upsilon_{P,t}d\betabar_P
\\
 =2\pi\sqrt{-1}
 \Bigl(
 (\del_{\xi_P}\alpha_P)
 \cdot (\xi_P\upsilon_{P,t})
 \cdot
 (\del_{\xi_P}\betabar_P)
\Bigr)_{|\xi_P=0}.
\end{multline}
Because
$\eta_t=
 \frac{2^3}{3^2}
 \xi_P\upsilon_{P,t}
 d\xi_P$
and
$\phi_P^{1/2}
=\frac{2^3}{3^2} \xi_P^2\,d\xi_P$,
the following holds on
$\Sigma_{\theta|X_P(\kappa_0)}$:
\[
 \nu\cdot\taubar\cdot
 \eta_t\cdot(\phi^{1/2})^{-1}
 =\Bigl(
 (\del_{\xi_P}\alpha_P)
 \cdot
 (\del_{\xi_P}\betabar_P)
 \cdot
 (\xi_P\upsilon_{P,t})
 \Bigr)
 \left(
 \frac{d\xi_P}{\xi_P}
 \right)^2.
\]
Then, we obtain
the estimate (\ref{eq;23.10.7.2}).
\hfill\qed

\vspace{.1in}
By Lemma \ref{lem;23.10.7.3},
we obtain the following:
\begin{equation}
\label{eq;23.10.8.22}
 \bigl(
 \ttH_{P,t}(\nu),
 \ttV_{\kappa,t}'(\tau)
 \bigr)_{L^2,h_t,X^{(2)}_{P,\kappa}}
=
 \langle \nu,\overline{\tau}\rangle^{\aux}_{P,h_t}
+\bigl(
 \ttHtilde_{P,t}(\nu),
 \ttV_{\kappa,t}'(\tau)
 \bigr)_{L^2,h_t,X^{(2)}_{P,\kappa}}
+O\bigl(e^{-8\kappa t}\|\nu\|_{L^2}\cdot\|\tau\|_{L^2}\bigr).
\end{equation}

By Lemma \ref{lem;23.10.7.20} below, we obtain
\[
\bigl(
 \ttHtilde_{P,t}(\nu),
 \ttV_{\kappa,t}'(\tau)
 \bigr)_{L^2,h_t,X^{(2)}_{P,\kappa}}
=2\sqrt{-1}\int_{\del X_{P,\kappa}^{(2)}}
 \Tr\Bigl(
 F_{\alpha}\cdot
 (\ttV_{\kappa,t}'(\nu)^{1,0})^{\dagger}_{h_t}
 \Bigr).
\]
We have
\[
 \Tr\Bigl(
 F_{\alpha}\cdot
 (\ttV_{\kappa,t}'(\nu)^{1,0})^{\dagger}_{h_t}
 \Bigr)
 =\Tr\Bigl(
 F_{\alpha}\cdot
 [\theta^{\dagger}_{h_t},F_{\betabar}]
 \Bigr)
 =\Tr\Bigl(
 \theta^{\dagger}[F_{\betabar},F_{\alpha}]
 \Bigr)=0.
\]
Thus, we obtain Proposition \ref{prop;23.9.27.4}.
\hfill\qed

\subsection{Appendix}

\subsubsection{The error terms for metrics}
\label{subsection;23.10.4.1}

Let $\gamma>0$.
Let $A^{(\ell)}_t$ $(\ell=1,2,\,\,t\geq 1)$ be $2\times 2$-matrices.
We assume that
\begin{itemize}
 \item $|(A^{(\ell)}_t)_{i,j}|=O(e^{-\gamma t})$ $(i\neq j)$.
 \item $(A^{(\ell)}_t)_{i,i}=\alpha^{(\ell)}_i+O(e^{-2\gamma t})$
       for some complex numbers $\alpha^{(\ell)}_i$.
\end{itemize}

Let $H_t$ $(t\geq 1)$ be Hermitian matrices.
We set $B_t=H_t-I_2$.
Assume the following condition.
\begin{itemize}
 \item $|B_t|=O(e^{-\gamma t})$.
\end{itemize}

\begin{lem}
\label{lem;23.9.25.1}
We have
\[
\Tr(A^{(1)}(A_t^{(2)})^{\dagger}_{H_t})
=\alpha^{(1)}_1\overline{\alpha^{(2)}_1}
+\alpha^{(1)}_2\overline{\alpha^{(2)}_2}
+O(e^{-2\gamma t}).
\]
In particular, we have
\[
  |A^{(1)}_t|_{H_t}^2=
\Tr(A^{(1)}_t (A^{(1)}_t)_{H_t}^{\dagger})
 =|\alpha^{(1)}_1|^2+|\alpha^{(1)}_2|^2
 +O(e^{-2\gamma t}).
\]
\end{lem}
\pf
We have
\[
 (A_t^{(2)})_{H_t}^{\dagger}
=\Hbar_t^{-1}
\lefttop{t}\Abar_t^{(2)}
\cdot \Hbar_t
=\lefttop{t}\Abar_t^{(2)}
 +[\lefttop{t}\Abar_t^{(2)},\Bbar_t]
 +O(e^{-2\gamma t}).
\]
We note that
$[\lefttop{t}\Abar_t^{(2)},\Bbar_t]_{i,i}
=O(e^{-2\gamma t})$ $(i=1,2)$
and 
$[\lefttop{t}\Abar_t^{(2)},\Bbar_t]_{i,j}
=O(e^{-\gamma t})$ $(i\neq j)$.
As a result,
we obtain the claim of the lemma.
\hfill\qed

\subsubsection{A variant}
\label{subsection;23.10.5.12}

Let $H_t$ be a family of Hermitian matrices
such that
$B_t=H_t-I_2=O(e^{-\gamma t})$.
We set
\[
 B_t=
 \left(
 \begin{array}{cc}
  a_t& b_t \\ \bbar_t & c_t
 \end{array}
 \right).
\]

We assume $\det(H_t)=1$.
We set
\[
 \Pi_1=
 \left(
\begin{array}{cc}
 1 & 0 \\ 0 & 0
\end{array}
 \right),
 \quad\quad
 \Pi_2=
 \left(
\begin{array}{cc}
 0 & 0 \\ 0 & 1
\end{array}
 \right).
\]
We have
\[
 (\Pi_1)^{\dagger}_{H_t}
=\Hbar_t^{-1}\Pi_1\Hbar_t
=\Pi_1+
[\Pi_1,\Bbar_t]
+O(e^{-2\gamma t})
=\Pi_1
+\left(
 \begin{array}{cc}
  0 & \bbar_t \\ -b_t & 0
 \end{array}
 \right)
 +O(e^{-2\gamma t}).
\]
We also have
\[
  (\Pi_2)^{\dagger}_{H_t}
=\Pi_2
+\left(
 \begin{array}{cc}
  0 & -\bbar_t \\ b_t & 0
 \end{array}
 \right)
 +O(e^{-2\gamma t}).
\]

Let $\alpha_i,\beta_i\in\cnum$ $(i=1,2)$.
We set
\[
 G_t=
 \alpha_1\Pi_1+\alpha_2\Pi_2
 +\beta_1 (\Pi_1)^{\dagger}_{H_t}
 +\beta_2 (\Pi_2)^{\dagger}_{H_t}
 =
 \left(
 \begin{array}{cc}
  \alpha_1+\beta_1 & (\beta_1-\beta_2)\bbar_t \\
 (\beta_2-\beta_1)b_t & \alpha_2+\beta_2
 \end{array}
 \right)
+O(e^{-2\gamma t}).
\]
\begin{lem}
\label{lem;23.9.26.1}
We obtain
\[
 \bigl|
 G_t
 \bigr|^2_{H_t}
=|\alpha_1+\beta_1|^2+|\alpha_2+\beta_2|^2
+O(e^{-2\gamma t}).
\]
\end{lem}
\pf
It follows from Lemma \ref{lem;23.9.25.1}.
\hfill\qed

\subsubsection{$L^2$-product}

Let $(E,\delbar_E,\theta)$ be a Higgs bundle on
a compact Riemann surface $X$.
Let $W\subset X$ be an open subset.
Let $h$ be a harmonic metric.
Let $\rho_i=\rho_i^{1,0}+\rho_i^{0,1}$ $(i=1,2)$
be $C^{\infty}$-sections of
$\End(E)\otimes\Omega^{1,0}\oplus\End(E)\otimes\Omega^{0,1}$.
We set
\[
 (\rho_1,\rho_2)_{L^2,h,W}
 =2\sqrt{-1}\int_W
 \Tr\Bigl(
 \rho_1^{1,0}
 (\rho_2^{1,0})^{\dagger}
-\rho_1^{0,1}
 (\rho_2^{0,1})^{\dagger}
 \Bigr).
\]
Let $\rho_i$ $(i=1,2)$
be $C^{\infty}$-sections 
of $\End(E)$.
Let $g_X$ be a K\"ahler metric of $X$.
Let $\dvol_{g_X}$ be the volume form.
For any $W\subset X$,
we set
\[
 (\rho_1,\rho_2)_{L^2,h,g_X,W}
 =\int_W
 \Tr\bigl(
 \rho_1\cdot(\rho_2)^{\dagger}_h
 \bigr)
 \dvol_{g_X}.
\]

If $W=X$,
$(\rho_1,\rho_2)_{L^2,h,X}$
and $(\rho_1,\rho_2)_{L^2,h,g_X,X}$
are denoted by
$(\rho_1,\rho_2)_{L^2,h}$
and $(\rho_1,\rho_2)_{L^2,h,g_X}$.
We set
$\|\rho\|_{L^2,h}=\sqrt{(\rho,\rho)_{L^2,h}}$
or
$\|\rho\|_{L^2,h}=\sqrt{(\rho,\rho)_{L^2,h,g_X}}$.

\subsubsection{Stokes formula}

Let $(E,\delbar_E,\theta,h)$ be a harmonic bundle on $X$.
Let $a$ be a $C^{\infty}$-section of $\End(E)$
and $b=b^{1,0}+b^{0,1}$ be a $C^{\infty}$ $\End(E)$-valued $1$-form.
\begin{lem}
\label{lem;23.10.7.20}
The following holds:
\begin{equation}
\label{eq;23.10.7.10}
 \Tr\Bigl(
 \del_{E,h} a\wedge (b^{1,0})^{\dagger}_h
-[\theta^{\dagger}_h,a] \wedge (b^{0,1})^{\dagger}_h
 \Bigr)
 =\del \Tr(a(b^{1,0})^{\dagger})
 -\Tr\bigl(a\cdot(\delbar_E+\ad\theta)b\bigr).
\end{equation}
\begin{equation}
\label{eq;23.10.7.11}
 \Tr\Bigl(
 [\theta,a]\wedge (b^{1,0})^{\dagger}_h
-\delbar a\wedge (b^{0,1})^{\dagger}_h
 \Bigr)
 =-\delbar \Tr(a(b^{0,1})_h^{\dagger})
 +\Tr\Bigl(
 a\bigl((\del_{E,h}+\ad\theta^{\dagger}_h)b\bigr)^{\dagger}_h
 \Bigr).
\end{equation}
 \end{lem}
\pf
We have
\[
 \Tr(\del a\wedge (b^{1,0})^{\dagger}_h)
=\del \Tr\bigl(a (b^{1,0})^{\dagger}_h\bigr)
-\Tr(a\wedge (\delbar b^{1,0})^{\dagger}_h)
\]
We have
\[
 \Tr\Bigl(
 [\theta_h^{\dagger},a]\wedge (b^{0,1})_h^{\dagger}
 \Bigr)
 =\Tr\Bigl(
 (\theta^{\dagger}_ha-a\theta^{\dagger}_h)
 \wedge (b^{0,1})^{\dagger}_h
 \Bigr)
=-\Tr\Bigl(
 a\bigl[
 \theta_h^{\dagger},
 (b^{0,1})^{\dagger}_h
 \bigr]
 \Bigr)
=\Tr\Bigl(
 a\bigl[
 \theta,
 b^{0,1}
 \bigr]^{\dagger}_h
 \Bigr).
\]
Thus, we obtain (\ref{eq;23.10.7.10}).
We obtain (\ref{eq;23.10.7.11}) similarly.
\hfill\qed

\begin{cor}
Let $W\subset X$ be a relatively compact open subset
with smooth boundary.
Then, the following holds.
\begin{equation}
\bigl(
(\del_{E,h}+\ad\theta^{\dagger}_h)a,b
\bigr)_{L^2,h,W}
=2\sqrt{-1}
\int_{\del W}
\Tr\bigl(a(b^{1,0})^{\dagger}_h\bigr)
+\bigl(
 a,(\del_{E,h}+\ad\theta^{\dagger}_h)^{\ast}_{h,g_X}b
\bigr)_{L^2,h,g_X,W},
\end{equation}
\begin{equation}
 \bigl(
  (\delbar_E+\ad\theta)a,b
\bigr)_{L^2,h,W}
=-2\sqrt{-1}
\int_{\del W}
\Tr\bigl(a(b^{0,1})^{\dagger}_h\bigr)
+\bigl(
 a,(\delbar_E+\ad\theta)^{\ast}_{h,g_X}b
\bigr)_{L^2,h,g_X,W}.
\end{equation}
\hfill\qed
\end{cor}

\end{document}